\documentclass[11pt,twoside]{article}

\def\aa{\alpha}
\def\bb{\beta}
\def\gg{\gamma}
\def\dd{\delta}
\def\ee{\epsilon}
\def\ff{\varphi}
\def\kk{\kappa}

\def\ll{\lambda}
\def\oo{\omega}
\def\ss{\sigma}

\def\zz{\zeta}

\def\cc{{\cal C}}
\def\cd{{\cal D}}
\def\ce{{\cal E}}
\def\cf{{\cal F}}

\def\cl{{\cal L}}
\def\cn{{\cal N}}
\def\co{{\cal O}}
\def\cp{{\cal P}}
\def\cs{{\cal S}}
\def\ct{{\cal T}}

\def\CC{\mathbbm{C}}

\def\FF{\Phi}

\def\II{\mathbbm{I}}

\def\NN{\mathbbm{N}}

\def\OO{\Omega}
\def\PP{\mathbbm{P}}

\def\RR{\mathbbm{R}}

\def\SS{\Sigma}

\def\ZZ{\mathbbm{Z}}
\def\pp{\partial}

\def\sym{Sp(n;\RR)}

\def\ra{\rightarrow}
\def\Ra{\Rightarrow}

\def\re{\preccurlyeq}  
\def\rp{\preceq}       

\def\Vol{\mbox{Vol}\,}
\def\de{\stackrel{\mbox{\scriptsize{def}}}{=}}

\def\sq{\Box}
\def\tr{\triangle}

\def\proof{\noindent {\it Proof. \;}}

\newcommand{\proofend}{\hspace*{\fill} $\Box$\\}
\newcommand{\diam}{\hspace*{\fill} $\Diamond$\\}

\newtheorem{theorem}{Theorem}[section]

\newtheorem{corollary}[theorem]{Corollary}

\newtheorem{proposition}[theorem]{Proposition}

\newtheorem{lemma}[theorem]{Lemma}

\newtheorem{conjecture}[theorem]{Conjecture}

\newtheorem{remark}[theorem]{Remark}

\newtheorem{example}[theorem]{Example}

\usepackage{amssymb,amsmath,bbm,epsfig,psfrag,epic,eepic,verbatim,latexsym}

\begin{document}

\begin{titlepage}
\title{On symplectic folding}    
$   $ \\
$   $ \\

\author{Felix Schlenk}
\date{10. January 1999}

\end{titlepage}
\maketitle

\begin{abstract}
\noindent
We study the rigidity and flexibility of symplectic embeddings of simple
shapes. It is first proved that under the condition $r_n^2 \le 2 r_1^2$ 
the symplectic ellipsoid $E(r_1, \dots ,r_n)$ with radii 
$r_1 \le \dots \le r_n$
does not embed in a ball of radius strictly smaller than $r_n$. We
then use symplectic folding to see that this condition is sharp and to
construct some nearly optimal embeddings of ellipsoids and polydiscs
into balls and cubes. 
It is finally shown that any connected symplectic
manifold of finite volume may be asymptotically filled with skinny
ellipsoids or polydiscs.

\end{abstract}

\tableofcontents

\setcounter{secnumdepth}{4}
\setcounter{equation}{0}

\section{Introduction}

Let $U$ be an open subset of $\RR^n$ which is diffeomorphic to a
ball, endow $U$ with the Euclidean volume form $\OO_0$, and let $(M,
\OO)$ be any connected $n$-dimensional volume manifold. 
Then $U$ embeds into $M$
via a volume preserving map if and only if 
$\Vol (U, \OO_0) \le \Vol (M, \OO)$.
(A proof of this ``folk-theorem'' is given below.)

Let $\oo_0 = \sum_{i=1}^n dx_i \wedge dy_i$ be the standard symplectic
form on $\RR^{2n}$ and equip any open subset $U$ of $\RR^{2n}$ with this
form. An embedding $\ff \colon U \hookrightarrow \RR^{2n}$ is called
symplectic, if $\ff^* \oo_0 = \oo_0$. In particular, every symplectic
embedding preserves the volume and the orientation. In dimension two,
the converse holds true. In higher dimensions, however, strong symplectic
rigidity phenomena appear. A spectacular example for this is Gromov's
Nonsqueezing Theorem \cite{G1}, which states that a ball $B^{2n}(r)$ of
radius $r$ symplectically embeds in the standard symplectic cylinder
$B^2(R) \times \RR^{2n-2}$ if and only if $r \le R$.
This and many other rigidity results for symplectic maps could later be
explained via symplectic capacities which arose from the variational
study of periodic orbits of Hamiltonian systems (see \cite{HZ} and the references therein). 
 
On the other hand, the flexibility of symplectic codimension 2
embeddings of open manifolds \cite[p.\ 335]{G2} implies that given any
symplectic ball $B^{2n-2}$ in $\RR^{2n-2}$ and a symplectic manifold
$(M^{2n}, \oo)$, there exists an $\ee >0$ such that $B^{2n-2} \times
B^2(\ee)$ symplectically embeds in $M$ (see \cite[p.\ 579]{FHW} for
details).

\smallskip
The aim of this work is to investigate the zone of transition between
rigidity and flexibility in symplectic topology. 
Unfortunately, symplectic capacities can be computed only for very
special sets, and there is still not much known about what one can do
with a symplectic map. We thus look at a model situation.
Let
\[
E(a_1, \dots , a_n) = \left\{ (z_1, \dots , z_n ) \in \CC^n \; \bigg| \,
\sum_{i=1}^n \frac{\pi |z_i|^2}{a_i}<1 \right\}
\]
be the open symplectic ellipsoid with radii $\sqrt{a_i / \pi}$, and
write $D(a)$ for the open disc of area $a$ and $P(a_1, \dots, a_n)$ for
the polydisc $D(a_1) \times \dots \times D(a_n)$. 
Since a permutation of the symplectic coordinate planes is a (linear)
symplectic map, we may assume $a_i \le a_j$ for $i<j$. Finally, denote
the ball $E^{2n}(a, \dots, a)$ by $B^{2n}(a)$ and the ``$n$-cube'' 
$P^{2n}(a, \dots, a)$ by $C^{2n}(a)$.
We call any of these sets a simple shape.
We ask:
\\
\\
\noindent
{\it ``Given a simple shape $S$, what is the smallest ball $B$ and what
is the smallest cube $C$ such that $S$ symplectically fits into $B$ and $C$?'' }
\\
\\
Observe that embedding $S$ into a minimal ball amounts to minimizing its
diameter, while embedding $S$ into a minimal cube amounts to minimizing
its symplectic width.

\smallskip
Our main  rigidity result states that
for ``round'' ellipsoids the identity provides already the optimal embedding.
\\
\\
{\bf Theorem 1} \;
{\it
Let $a_n \le 2a_1$ and $a<a_n$. Then $E(a_1, \dots, a_n)$ does not embed
symplectically in $B^{2n}(a)$.
}
\\

\noindent
An ordinary symplectic capacity only shows that if $a < a_1$, there is
no symplectic embedding of $E(a_1, \dots, a_n)$ into $B^{2n}(a)$. Our
proof uses the first $n$ Ekeland-Hofer capacities.
For $n = 2$, Theorem 1 was proved in \cite{FHW} as an early application
of symplectic homology, but the argument given here is much simpler and
works in all dimensions.

\smallskip
Our first flexibility result states that Theorem 1 is sharp.
\\
\\
{\bf Theorem 2A} \;
{\it
Given any $\ee >0$ and $a>2 \pi$, there exists a symplectic embedding
\[
E^{2n}(\pi, \dots, \pi, a) \hookrightarrow B^{2n} \left( \frac{a}{2} +
\pi + \ee \right).
\]
}

\noindent
Lalonde and McDuff observed in \cite{LM1} that their technique of
symplectic folding can be used to prove Theorem 2A for $n = 2$. 
The symplectic folding construction considers a
\mbox{4-ellipsoid} as a fibration of discs of varying size over a disc and
applies the flexibility of volume preserving maps to both the base and
the fibres. It is therefore purely four dimensional in nature. We refine
the method in such a way that it will nevertheless be sufficient to
prove the result for arbitrary dimension.
\begin{figure}[h] 
 \begin{center}
  \psfrag{1}{$1$}
  \psfrag{2}{$2$}
  \psfrag{3}{$3$}
  \psfrag{4}{$4$}
  \psfrag{5}{$5$}
  \psfrag{6}{$6$}
  \psfrag{8}{$8$}
  \psfrag{12}{$12$}
  \psfrag{15}{$15$}
  \psfrag{20}{$20$}
  \psfrag{24}{$24$}
  \psfrag{a/pi}{$\frac{a}{\pi}$}
  \psfrag{A/pi}{$\frac{A}{\pi}$}
  \psfrag{seb/pi}{$\frac{s_{EB}}{\pi}$}
  \psfrag{leb/pi}{$\frac{l_{EB}}{\pi}$}
  \psfrag{inclusion}{$\mbox{inclusion}$}
  \psfrag{vol}{$\mbox{volume condition}$}
  \psfrag{cEH}{$c_{EH}$}
  \psfrag{folding}{$\mbox{folding once}$}
  \leavevmode\epsfbox{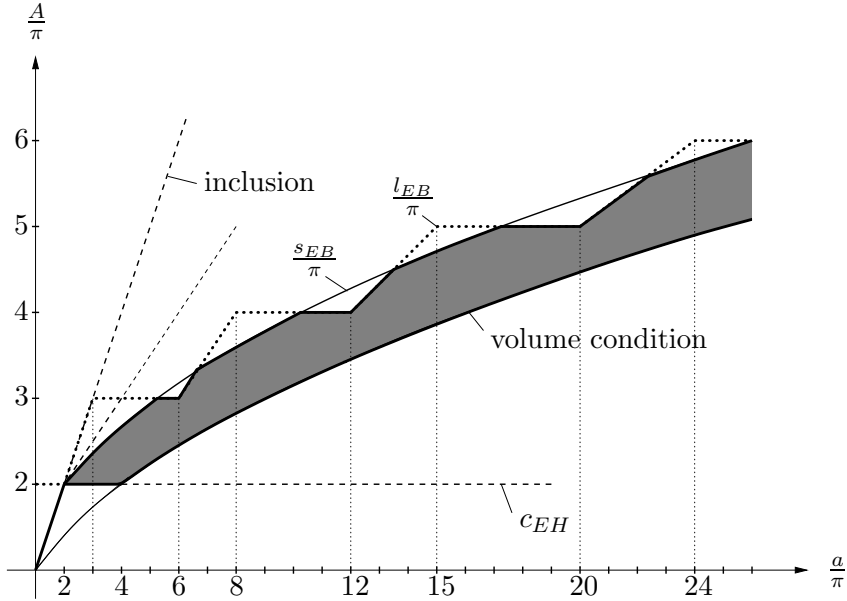}
 \end{center}
 \caption{What is known about $E(\pi,a) \hookrightarrow B^4(A)$} 
 \label{figure18}
\end{figure}
%
%

Theorem 1 and Theorem 2A shed some light on the power of Ekeland-Hofer
capacities: As soon as these invariants cease to imply that there is no
better embedding than the identity, there is indeed a better embedding.

For embeddings of ellipsoids into cubes, the same procedure yields a
similarly sharp result, but for embeddings of polydiscs into balls and
cubes the result is less satisfactory. In four dimensions, the precise
result is as follows.
\\
\\
{\bf Theorem 2B} \;
{\it
Let $\ee$ be any positive number.
\begin{itemize}
\item[(i)]
Let $a>\pi$. 
Then there is no symplectic embedding of $E(\pi,a)$ into
$C^4(\pi)$, but $E(\pi,a)$ symplectically embeds in
$C^4(\frac{a+\pi}{2}+\ee)$.
\item[(ii)]
Let $a>2\pi$. 
Then $P(\pi,a)$ symplectically embeds in $B^4(\frac{a}{2}+ 2\pi + \ee)$
as well as in $C^4(\frac{a}{2}+ \pi + \ee)$.
\end{itemize}
}
\bigskip
\noindent
{\bf Question 1} \,
{\it
Does $P(\pi,2 \pi)$ symplectically embed in $B^4(A)$ for some $A< 3 \pi$
or in $C^4(A)$ for some $A< 2 \pi$?
}
\\
\\
Both, Theorem 2A and Theorem 2B as well as its higher dimensional
version can be substantially improved by multiple folding. Let us
discuss the result in case of embeddings of 4-ellipsoids into 4-balls
(cf.\ Figure \ref{figure18}). Let $s_{EB}(a)$ be the function describing
the best embeddings obtainable by symplectic folding. It turns out that
\[ \limsup_{\ee \ra 0^+} \frac{s_{EB}(2 \pi + \ee) - 2 \pi}{\ee} =
\frac{3}{7}. \]
{\bf Question 2} \,
{\it
Let $f_{EB}(a) = \inf \{ A \, | \, E(\pi,a) \mbox{ symplectically embeds
in } B^4(A) \}$. How does $f_{EB}$ look like near $2 \pi$? In
particular, 
\[ \limsup_{\ee \ra 0^+} \frac{f_{EB}(2 \pi + \ee) - 2 \pi}{\ee} <
\frac{3}{7} \, ? \] 
}
Moreover, as $a \ra \infty$ the image of $E(\pi,a)$ fills up an
arbitrarily large percentage of the volume of $B^4(s_{EB}(a))$.
This can also be seen via a Lagrangian folding method, which was
developed by Traynor in \cite{T} and yielded the best previously
known results for
the above embedding problem (see the curve $l_{EB}$ in Figure
\ref{figure18}). Symplectic folding, however, may be used to prove that 
any connected symplectic manifold
$(M, \oo)$ of finite volume can be asymptotically filled by skinny
ellipsoids and polydiscs: For $a > \pi$ set
\[
p_a^E (M^{2n}, \oo) = \sup_\aa \frac{\Vol ( E^{2n}(\aa \pi, \dots, \aa \pi,
\aa a))}{\Vol (M, \oo)},
\]
where the supremum is taken over all $\aa$ for which $E^{2n}(\aa \pi, \dots,
\aa \pi, \aa a)$ symplectically embeds in $(M, \oo)$, and define 
$p_a^P (M, \oo)$ in a similar way.
\\
\\
{\bf Theorem 3 } 
{\it
$\lim_{a \ra \infty} p_a^E (M, \oo)$ \,\! and \,\!
$\lim_{a \ra \infty} p_a^P (M, \oo)$ \,\!
exist and equal 1. }
\\
\\
This result exhibits that in the limit symplectic rigidity disappears.
We finally give estimates of the convergence speed from below.

\smallskip
Appendix A  provides computer programs necessary to compute the optimal
embeddings of ellipsoids into a 4-ball and a 4-cube obtainable by our
methods, and in Appendix B we
give an overview on known results on the Gromov width of closed
symplectic manifolds.
\\
\\
%
%
%
{\bf Acknowledgement.} I greatly thank Dusa McDuff for her fine criticism
on an earlier more complicated attempt towards \mbox{Theorem 2A}, which gave worse
estimates, and for having explained to me the main point of the folding
construction.

I also thank Leonid Polterovich for suggesting to me to look closer at
Lagrangian folding.
%
\\

\section{Rigidity}  \label{rigidity}

Throughout this paper, if there is no explicit mention to the contrary, all maps will be assumed to be symplectic. In
dimension two this just means that they preserve the orientation and the
area.

Denote by $\co (n)$ the set of bounded domains in $\RR^{2n}$ endowed with the
standard symplectic structure $\oo_0 = \sum_{i=1}^n dx_i \wedge dy_i$.
Given $U \in \co (n)$, write $|U|$ for the volume of $U$ with respect to
the Euclidean volume form $\Omega_0 = \frac{1}{n!} \oo_0^n$.
Let $\cd (n)$ be the group of symplectomorphisms of $\RR ^{2n}$ and
$\cd_c(n)$ respectively $\sym$ the subgroups of compactly supported
respectively linear symplectomorphisms of $\RR ^{2n}$. 
Define the following relations on $\co (n)$: 
\medskip
\[
\begin{array}{lcl}
U \le_1 V  &  \Longleftrightarrow  &  \mbox{There exists a $\ff \in \sym$ with
$\ff (U) \subset V$.}  \\
U \le_2 V &  \Longleftrightarrow   &  \mbox{There exists a $\ff \in \cd
(n)$ with $\ff (U) \subset V$.} \\
U \le_3 V & \Longleftrightarrow    &  \mbox{There exists a symplectic embedding $\ff \colon U \hookrightarrow V$.}
\end{array}
\]
Of course, $\le_1 \; \Ra \; \le_2 \; \Ra \; \le_3$, but all the
relations are different:
That $\le_1$ and $\le_2$ are different  
is well known (see (\ref{equation2}) below and Traynor's theorem stated
at the beginning of \mbox{section \ref{flexibility}}).
The construction of sets $U$ and $V \in \co(n)$ with $U \le_3 V$ 
but $U \not\le_2 V$ relies on the following simple observation. Suppose
that $U$ and $V$ not only fulfill $U \le_3 V$ but are symplectomorphic,
whence, in particular, $|U| = |V|$. Thus, if $U \le_2 V$ and $\ff$ is a
map realizing $U \le_2 V$, no point of $\CC^n \setminus U$ can be
mapped to $V$, and we conclude that $\ff(\partial U) = \partial V$.
In particular, the characteristic foliations on $\partial U$ and
$\partial V$ are isomorphic, and if $\partial U$ is of contact type,
then so is $\partial V$ (see \cite{HZ} for basic notions in Hamiltonian
dynamics).  

Let now $U = B^{2n}(\pi)$, let
\[
SD = D(\pi) \setminus \{ (x,y) \; | \, x \ge 0, y=0 \}
\]
be the slit disc and set $V = B^{2n}(\pi) \cap (SD \times \dots \times SD)$. 
Traynor proved in \cite{T} that for $n \le 2$, $V$ is symplectomorphic to
$B^{2n}(\pi)$. But $\partial U$ and $\partial V$ are not even
diffeomorphic. 
For $n \ge 2$ very different examples were found in \cite{ElH} and
\cite{C}.
Theorem 1.1 in \cite{ElH} and its proof show that 
there exist $U,V \in \co (n)$ with smooth convex
boundaries such that $U$ and $V$ are symplectomorphic and
$C^\infty$-close to $B^{2n}(\pi)$, but the characteristic foliation of
$\partial U$ contains an isolated closed orbit while the one of
$\partial V$ does not. 
And Corollary A in \cite{C} and its proof imply that
given any $U \in \co (n)$, $n \ge 2$, with
smooth boundary $\partial U$ of contact type, there exists a
symplectomorphic and $C^0$-close $V \in \co (n)$ whose boundary is not
of contact type.

We in particular see that even for $U$ being a ball, 
$\le_3$ does not imply $\le_2$.
\\
\\
In order to detect some rigidity via the above relations we therefore
must pass to a small subcategory of sets: 

Let $\ce (n)$ be the collection
of symplectic ellipsoids described in the introduction
\[
\ce (n) = \{ E(a) = E(a_1, \dots , a_n) \}
\]
and write $\re_i$ for the restrictions of the relations $\le_i$
to $\ce (n)$. 

Notice again that
\[
\re_1 \;\; \Longrightarrow \;\; \re_2 \;\; \Longrightarrow \;\; \re_3\,.
\]
$\re_2$ and $\re_3$ are actually very similar: Since ellipsoids are
starlike, we may apply Alexander's trick to prove the extension after
restriction principle (see \cite{EH1} for details), which tells us that
given any embedding
$\ff \colon E(a) \hookrightarrow E(a')$ and any $\dd \in \, ]0,1[$ we can find a
$\psi \in \cd (n)$ which coincides with $\ff$ on $E(\dd a)$; hence
\begin{eqnarray} \label{equation1}
E(a) \re_3 E(a') \;\; \Longrightarrow \;\; E(\dd a) \re_2 E(a') 
\quad \mbox{ for all } \dd \in ]0,1[ \,. 
\end{eqnarray}
It is, however, not clear whether $\re_2$ and $\re_3$ are the same: While
\mbox{Theorem \ref{theoremrigidity}} proves this under an additional condition, the folding
construction of \mbox{section \ref{flexibility}} suggests that they are different in
general. But let us first prove a general and common rigidity property of these relations:

\begin{proposition} \label{proposition1}
The relations $\re_i$ are partial orderings on $\ce (n)\,$.
\end{proposition}

\proof
The relations are clearly reflexive and transitive, so we are left with
identitivity. Of course, the identitivity of $\re_3$ implies the one of
$\re_2$ which, in its turn, implies the one of $\re_1$. We still prefer
to give independent proofs which use tools whose difficulty is
about proportional to the depth of the results.

It is well known from linear symplectic algebra \cite[p.\ 40]{HZ} that 
\begin{equation} \label{equation2}
E(a) \re_1 E(a') \;\; \Longleftrightarrow \;\; a_i \le a_i' \quad \mbox{ for all } i,
\end{equation}
in particular $\re_1$ is identitive. 

Given $U \in \co (n)$ with smooth boundary $\pp U$, the spectrum $\ss
(U)$ of $U$ is defined to be the collection of the actions of closed
characteristics on $\pp U$. It is clearly
invariant under $\cd (n)$, and for an ellipsoid it is given by
\[
\ss (E(a_1, \dots, a_n)) = \{d_1(E) \le d_2(E) \le \ldots \} 
\de \{k a_i \;| \, k \in \NN, \, 1 \le i \le n \}.
\]
Let now $\ff$ be a map realizing $E(a) \re_2 E(a')$.
$E(a) \re_2 E(a') \re_2 E(a)$ gives in particular
$|E(a)| = |E(a')|$,  and we conclude as above that $\ff (\pp E(a)) = \pp
E(a')$. This implies $\ss (E(a)) = \ss (E(a'))$ and the claim for
$\re_2$ follows. 

To prove identitivity of $\re_3$ recall that Ekeland-Hofer capacities \cite{EH2} provide us with a
whole family of symplectic capacities for subsets of $\CC^n$. They are
invariant under $\cd_c (n)$, and for an ellipsoid $E$ they are given by 
the spectrum:
\begin{equation} \label{spectrum}
\{c_1(E) \le c_2(E) \le \ldots \} = \{d_1(E) \le d_2(E) \le \ldots \}.
\end{equation}
First
observe that in the proof of the extension after restriction principle the
generating Hamiltonian can be chosen to vanish outside a large ball, so
the extension can be assumed to be in $\cd_c (n)$. This shows that in the
definition of $\re_2$ we may replace $\cd (n)$ by $\cd_c (n)$ without
changing the relation, and that Ekeland-Hofer capacities may be applied
to $\re_2$. Next observe that for any $i \in \{1,2,3\}$ and $\aa >0$ 
\begin{equation} \label{equation3}
E(a) \re_i E(a') \;\;  \Longrightarrow  \;\;  E(\aa a) \re_i E(\aa a'), 
\end{equation}
just conjugate the
given map $\ff$ with the dilatation by $\aa^{-1}.$ Applying this and
(\ref{equation1}) we see that for any $\dd_1, \dd_2 \in \,]0,1[$ the assumed
relations
\[
E(a) \re_3 E(a') \re_3 E(a)
\]
imply
\[
E(\dd_2 \dd_1 a) \re_2 E(\dd_1 a') \re_2 E(a),
\]
and now the monotonicity of all the $c_i = d_i$ immediately gives
$a=a'$.
\proofend

\bigskip

It is well known (we refer again to the beginning of
\mbox{section \ref{flexibility}}) that $\re_2$ does not imply
$\re_1$ in general. 
However, a suitable pinching condition guarantees that
``linear'' and ``non linear'' coincide:

\begin{theorem}  \label{theoremrigidity}
Let $\kk \in \, ] \frac{\pi}{2}, \pi [$. Then the following 
statements are equivalent:  

\begin{itemize}

\item[(i)]                                                            
  $B^{2n}(\kk) \re_1 E(a) \re_1 E(a') \re_1 B^{2n}(\pi)$  
\item[(ii)]            
 $B^{2n}(\kk) \re_2 E(a) \re_2 E(a') \re_2 B^{2n}(\pi)$  
\item[(iii)]            
 $B^{2n}(\kk) \re_3 E(a) \re_3 E(a') \re_3 B^{2n}(\pi)\,.$  

\end{itemize}
\end{theorem}
Theorem 1 follows from \mbox{Theorem \ref{theoremrigidity}}, (\ref{equation2}) and
(\ref{equation3}). For $n=2$, \mbox{Theorem \ref{theoremrigidity}} 
was proved in \cite{FHW}. That proof uses a deep
result by McDuff, namely that the space of symplectic embeddings of a
ball into a larger ball is unknotted, and then applies the isotopy
invariance of symplectic homology. 
%
%
However, Ekeland-Hofer capacities provide an easy proof. The crucial
point is that as true capacities they have - very much in contrast to
symplectic homology - the monotonicity property. 
\\
\\
{\it Proof of Theorem \ref{theoremrigidity}}.\;
(ii) $\Ra$ (i): 
By assumption we have  $B^{2n}(\kk) \re_2 E(a) \re_2
B^{2n}(\pi)$, so the first Ekeland-Hofer capacity $c_1$ gives 
\begin{equation} \label{equation1'}
\kk \le a_1 \le \pi
\end{equation}
and $c_n$ gives 
\begin{equation} \label{equation2'}
\kk \le c_n(E(a)) \le \pi.
\end{equation}
(\ref{equation1'}) and  $\kk > \frac{\pi}{2}$ imply $2a_1 > \pi$, whence the
only elements
in $\ss (E(a))$ possibly smaller than $\pi$ are $ a_1, \dots , a_n$. It
follows therefore from (\ref{equation2'}) that
$a_n = c_n(E(a))$, whence $c_i(E(a))= a_i \;\,(1 \le i \le
n)$. 
Similarly we
find $c_i(E(a'))= a_i' \;\,(1 \le i \le n)$, and from $E(a) \re_2 E(a')$ we conclude $a_i \le
a_i'$.

(iii) $\Ra$ (i) follows now by a similar reasoning as in the proof of
the identitivity of $\re_3$: Starting from
\[
B^{2n}(\kk) \re_3 E(a) \re_3 E(a') \re_3 B^{2n}(\pi),
\]
(\ref{equation1}) shows that for any $\dd_1, \dd_2, \dd_3 \in \, ]0,1[$
\[
B^{2n}(\dd_3 \dd_2 \dd_1 \kk) \re_2 E(\dd_2 \dd_1 a) \re_2 E(\dd_1 a')
\re_2 B^{2n}(\pi)\,.
\]
Choosing $\dd_1, \dd_2, \dd_3$ so large that $\dd_3 \dd_2 \dd_1 \kk > \frac{\pi}{2}$
we may apply the already proved implication to see
\[
B^{2n}(\dd_3 \dd_2 \dd_1 \kk) \re_1 E(\dd_2 \dd_1 a) \re_1 E(\dd_1 a) \re_1
B^{2n}(\pi),
\]
and since $\dd_1, \dd_2, \dd_3$ may be chosen arbitrarily
close to 1, (\ref{equation2}) shows that we are done.
\proofend

\section{Flexibility}  \label{flexibility}

As it was pointed out in the introduction, 
the flexibility of symplectic codimension 2 embeddings of open
manifolds implies that a condition as in \mbox{Theorem
1} is necessary for rigidity. 
An explicit necessary condition was first obtained by Traynor in
\cite{T}. Her construction may be extended in an obvious way 
(see \mbox{subsection \ref{lagrangianfolding}}, 
in particular Corollary \ref{corollarylagrangianfolding}
$(i)_E$) to prove
\\
\\
{\bf Theorem}
{\rm (Traynor, \cite[Theorem 6.4]{T})}
{\it For all $k \in \NN$ and $\ee > 0$ there exists a symplectic embedding}
\[ E \left(\frac{\pi}{k+1}, \pi, \dots , \pi , k \pi \right) 
\hookrightarrow B^{2n}(\pi + \ee). \]
However, neither this theorem nor any refined version yielded by the
Lagrangian method used in its proof can decide whether \mbox{Theorem 1} is
sharp (cf.\ \mbox{Figure \ref{figure18}}). Our first flexibility result
states that this is indeed the case:

\begin{theorem}  \label{theoremfolding}  
Let $a>2\pi$ and $\ee>0$. Then $E^{2n}(\pi, \dots, \pi , a)$ embeds
symplectically in $B^{2n}(\frac{a}{2} + \pi + \ee)$.
\end{theorem}

For $n=2$, this theorem together with \mbox{Theorem 1} gives a complete
answer to our question in the introduction, whereas for arbitrary $n$ it only
states that \mbox{Theorem 1} is sharp. We indeed cannot expect a much
better result
since (as is seen using Ekeland-Hofer capacities) $E^{2n}(\pi, 3 \pi, \dots,
3 \pi)$ does not embed in any ball of capacity strictly smaller than $3
\pi$. 
\\
\\
{\it Proof of Theorem \ref{theoremfolding}}. 
We will construct an embedding 
\[ \Phi \colon E(a,\pi) \hookrightarrow B^4 \left(\frac{a}{2} + \pi +
\ee \right) \]
satisfying
\begin{equation} \label{equation3a}
\pi | \Phi (z_1, z_2)|^2 < \frac{a}{2} + \ee + \frac{\pi^2 |z_1|^2}{a} +
\pi |z_2|^2 \quad \mbox{ for all } (z_1, z_2) \in E(a, \pi).
\end{equation}
The composition of the linear symplectomorphism 
\[ E^{2n}(\pi, \dots, \pi ,a) \ra E^{2n}(a, \pi, \dots,\pi) \]
with the restriction of 
$\Phi \times id_{2n-4}$ to $E^{2n}(a, \pi, \dots, \pi)$ 
is then the desired embedding.

\medskip

The great flexibility of 2-dimensional area preserving maps is basic for
the construction of $\FF$. We now make sure that we may describe such a
map by prescribing it on an exhausting and nested family of loops.
\\

\noindent
{\bf Definition\;}
A family $\cl$ of loops in a simply connected domain $U \subset \RR^2$
is called {\it admissible} if there is
a diffeomorphism $\bb \colon D(|U|) \setminus \{0\} \ra U \setminus \{ p
\}$ for some point $p \in U$ such that 
\begin{itemize}
\item[(i)]
concentric circles are mapped to elements of $\cl$
\item[(ii)]
in a neighbourhood of the origin $\bb$ is an orientation preserving isometry.
\end{itemize}

\begin{lemma} \label{lemmaarea}
Let $U$ and $V$ be bounded and simply connected domains in $\RR^2$ of
equal area and let $\cl_U$ respectively $\cl_V$ be admissible families of loops
in $U$ respectively $V$. Then there is a symplectomorphism between $U$ and $V$
mapping loops to loops.
\end{lemma}
{\bf Remark.}
The regularity condition (ii) imposed on the families taken into
consideration can be weakened. Some condition, however, is necessary as
is seen from taking $\cl_U$ a family of concentric circles and $\cl_V$ a
family of rectangles with round corners and width larger than a positive
constant.
\diam
\\
{\it Proof of Lemma \ref{lemmaarea}}.
We may assume that $(U,\cl_U) = (D(\pi R^2), \{r e^{i \phi} \})$, and
after reparametrizing the $r$-variable by a diffeomorphism of $]0, R[$
which is the identity near $0$ we may assume that $\bb$ maps the loop
$C(r)$ of radius $r$ to the loop $L(r)$ in $\cl_V$ which encloses the
area $\pi r^2$.
 
We now search for a family $h(r, \cdot)$ of diffeomorphisms of
$S^1$ such that the map $\aa$ given by $\aa (r e^{i \phi}) = \bb (r
e^{i h(r, \phi)})$ is a symplectomorphism. With other words, we look
for a smooth $h \colon ]0, R[ \times S^1 \ra S^1$ which is a diffeomorphism
for $r$ fixed and solves the initial value problem
\[
(\ast)  \;\;\;\;\;  \left\{
                       \begin{array}{lcl}   
                              \frac{\pp h}{\pp \phi} (r, \phi) & = &
                                 1/\det \bb' (r e^{i h(r, \phi)}) \\
                              h(r,0) & = & 0  
                       \end{array}   \right.          
\]             
View $\phi$ for a moment as a real variable. The existence and
uniqueness theorem for ordinary differential equations with parameter
yields a smooth map $h \colon ]0,R[ \, \times \RR \ra \RR$ satisfying
($\ast$). Thus, $h(r,\cdot)$ is a diffeomorphism of $\RR$, and it remains
to check that it is $2\pi$-periodic. But this holds since the map $\aa \colon
r e^{i \phi} \mapsto \bb (r e^{i h(r, \phi)})$ locally preserves the
volume and $\aa (C(r))$ is contained in the loop $L(r)$. 

Finally, $\aa$ is an isometry in a punctured neighbourhood of the
origin and thus extends to all of $D(\pi R^2)$.
\proofend

\begin{figure}[h] 
 \begin{center}
   \psfrag{a}{$a$}
   \psfrag{b}{$b$}
   \psfrag{u}{$u$}
   \psfrag{fibre}{fibre capacity}
   \leavevmode\epsfbox{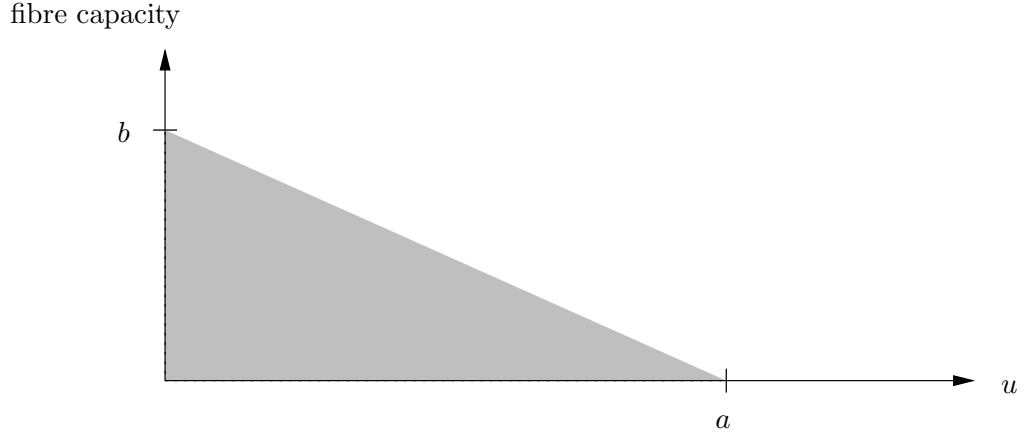}
 \end{center}
 \caption{The trapezoid $T(a,b)$} \label{figure1}
\end{figure}

While Traynor's construction relies mainly on considering a 4-ellipsoid
as a Lagrangian product of a rectangle and a triangle, we view it as a
trivial fibration over a symplectic disc with symplectic discs of
varying size as
fibres: More generally, define for $U \subset \CC$ open and $f \colon U \ra
\RR_{>0}$
\[
\cf (U,f) = \{ (z_1, z_2) \in \CC^2 \, | \, z_1 \in U, \, \pi |z_2|^2 <
f(z_1) \}.
\]
This is the trivial fibration over $U$ with fiber over $z_1$
the disc of capacity $f(z_1)$. For $\lambda \in \RR$ set
\[
U_\lambda = \{ z_1 \in U \, | \, f(z_1) \ge \lambda \}.
\]
Given two such fibrations $\cf (U, f)$ and $\cf (V,g)$, an embedding
$\psi \colon U \hookrightarrow V$ defines an embedding $\psi \times id
\colon \cf (U,f) \hookrightarrow \cf (V,g)$ if and only if $f(z_1) \le
g(\psi(z_1))$ for all
$z_1 \in U$, and under the assumption that all the sets $U_\lambda$ and 
$V_\lambda$ are connected, we see from Lemma \ref{lemmaarea} that inequalities
\[
\mbox{ area } U_\lambda < \mbox{ area } V_\lambda \quad \mbox{ for all }
\lambda
\]
are sufficient for the existence of an embedding $\cf (U,f)
\hookrightarrow \cf (V,g)$.
\\
\\
{\bf Example} 
{\rm (\cite[p.\ 54]{LM2})}
Let
$T(a,b) = \cf (R(a),g)$ with 
\[
R(a) = \{z_1 = (u,v) \, | \, 0 < u < a, \, 0 < v < 1 \}
\]
and $g(z_1) = g(u) = b-u$
be the trapezoid. 
We think of $T(a,b)$ as depicted in \mbox{Figure \ref{figure1}}.
\diam
%
%
%
\begin{lemma}  \label{lemma1}
For all $\ee > 0$,

\begin{itemize}

\item[ (i) ]  
$E(a,b)$ embeds in $T(a + \ee, b)$

\item[ (ii) ]
$T(a,b)$ embeds in $E(a + \ee, b)$.

\end{itemize}

\end{lemma}

\begin{figure}[h]
 \begin{center}
   \psfrag{D(a)}{$D(a)$}
   \psfrag{D(a+e)}{$D(a + \ee)$}
   \psfrag{z1}{$z_1$}
   \psfrag{alpha}{$\alpha$}
   \psfrag{o}{$\omega$}
   \psfrag{1}{$1$}
   \psfrag{v}{$v$}
   \psfrag{u}{$u$}
   \psfrag{a}{$a$}
   \psfrag{a+e}{$a+ \ee$}
   \leavevmode\epsfbox{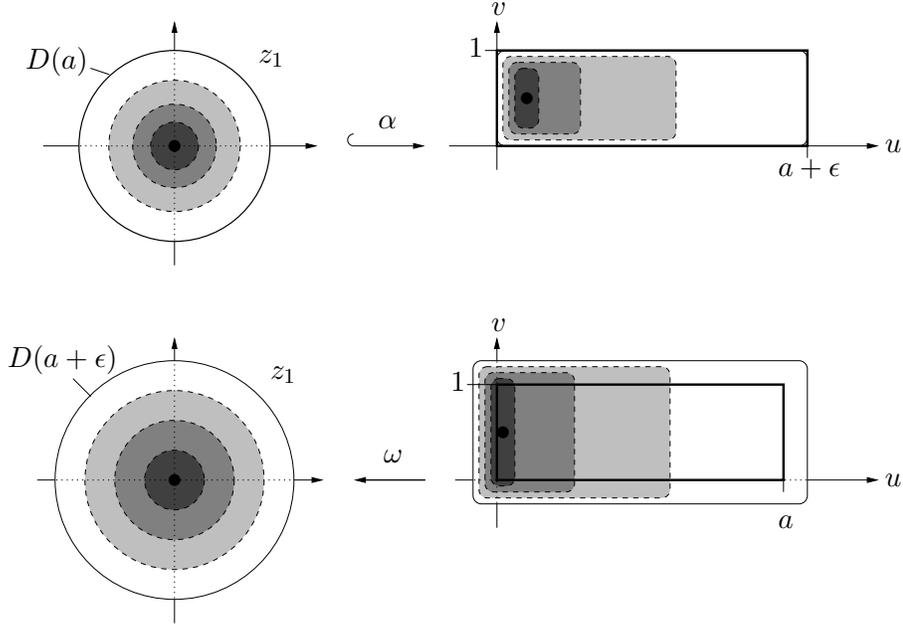}
 \end{center}
 \caption{The first and the last base deformation} \label{figure2}
\end{figure}

\proof
$E(a,b)$ is described by $U = D(a)$ and
$f(z_1) = b\,(1- \frac{\pi |z_1|^2}{a}).$ 
For (i) look at $\alpha$ and for (ii) at $\oo$ in \mbox{Figure \ref{figure2}}.
%
%
The symplectomorphism $\omega$ is defined on a round neighbourhood of $R(a)$.
\proofend

\mbox{Lemma \ref{lemma1}} and its proof readily imply that in order to
construct for any $a>2\pi$ and $\ee >0$ an embedding $\Phi$ satisfying
(\ref{equation3a}) it is enough to find for any $a>2\pi$ and $\ee >0$ an
embedding $\Psi \colon T(a,\pi) \hookrightarrow T(\frac{a}{2} + \pi +
\ee, \frac{a}{2} + \pi + \ee )$, $\,(u,v,z_2) \mapsto (u', v', z_2')$
satisfying

\begin{equation}  \label{equation4}
u'+\pi |z_2'|^2 < \frac{a}{2} + \ee + \frac{\pi u}{a} + \pi |z_2|^2 \quad
\mbox{ for all } (u,v,z_2) \in T(a, \pi ).
\end{equation}

\subsection{The folding construction}  \label{thefoldingconstruction}

The idea in the construction of an embedding $\Psi$ satisfying
(\ref{equation4}) is to separate the small fibres from the large
ones and then to fold the two parts on top of each other.
\\
\\
{\bf Step 1.} 
Following \cite[Lemma 2.1]{LM2} we first separate the
``low'' regions over $R(a)$ from the ``high'' ones: 

\begin{figure}[h] 
 \begin{center}
  \psfrag{v}{$v$}
  \psfrag{u}{$u$}
  \psfrag{1}{$1$}
  \psfrag{P1}{$P_1$}
  \psfrag{P2}{$P_2$} 
  \psfrag{L}{$L$}
  \psfrag{U}{$U$} 
  \psfrag{d}{$\dd$}
  \psfrag{a/2+d}{$\frac{a}{2}+ \dd$} 
  \psfrag{a/2+pi/2+11d}{$\frac{a+\pi}{2}+ 9 \dd$}
  \psfrag{a+pi/2+12d}{$a+ \frac{\pi}{2} + 10 \dd$} 
  \psfrag{f}{$f$}
  \psfrag{pi}{$\pi$}
  \psfrag{pi/2}{$\frac{\pi}{2}$}
  \leavevmode\epsfbox{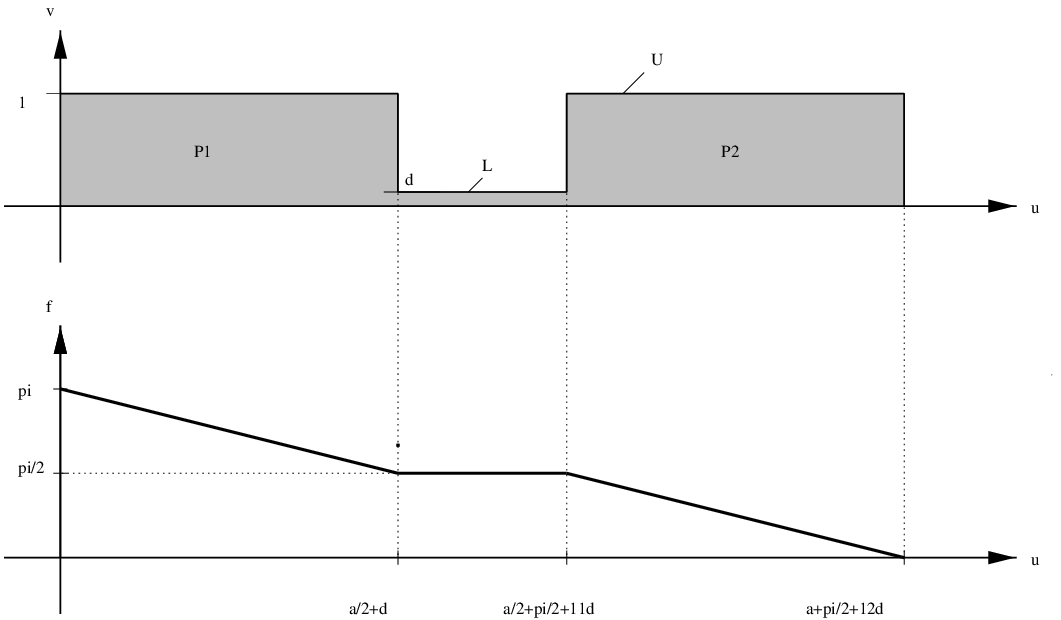}
 \end{center}
 \caption{Separating the low fibres from the large fibres} \label{figure3}
\end{figure}
%
%

Let $\dd >0$ be small. Let $\cf$ be described by $U$ and $f$ as in
\mbox{Figure \ref{figure3}} and write
\begin{gather*}
P_1 = U \cap \left\{ u \le \frac{a}{2}+ \dd \right\},                  \\ 
P_2 = U \cap \left\{ u \ge \frac{a+ \pi}{2} + 9 \dd \right\}           \\
L = U \setminus (P_1 \cup P_2).
\end{gather*}
It is clear from the discussion at the beginning of the proof that
there is an embedding $\bb \times id \colon T(a, \pi) \hookrightarrow \cf$
with

\begin{equation} \label{equation5}
\bb \, |_{ \left\{ u< \frac{a}{2} -\dd \right\}} = id \quad \mbox{ and } 
\quad \bb \, |_{ \left\{ u> \frac{a}{2} + \dd \right\}} = 
id + \left( \frac{\pi}{2} + 10 \dd ,0 \right).
\end{equation} 

\medskip
\noindent
{\bf Step 2.}
We next map the fibers into a convenient shape: 
\begin{figure}[h] 
 \begin{center}
  \psfrag{z2}{$z_2$}
  \psfrag{-dd}{$-\dd$}
  \psfrag{ss}{$\ss$}
  \psfrag{y}{$y$}         
  \psfrag{x}{$x$}
  \psfrag{-1/2}{$-\frac{1}{2}$}
  \psfrag{1/2}{$\frac{1}{2}$}
  \psfrag{-pi/2-2d}{$-\frac{\pi}{2} -2 \dd$}
  \psfrag{-pi-d}{$-\pi - \dd$}
  \psfrag{Ri}{$R_i$}
  \psfrag{Re}{$R_e$}
  \leavevmode\epsfbox{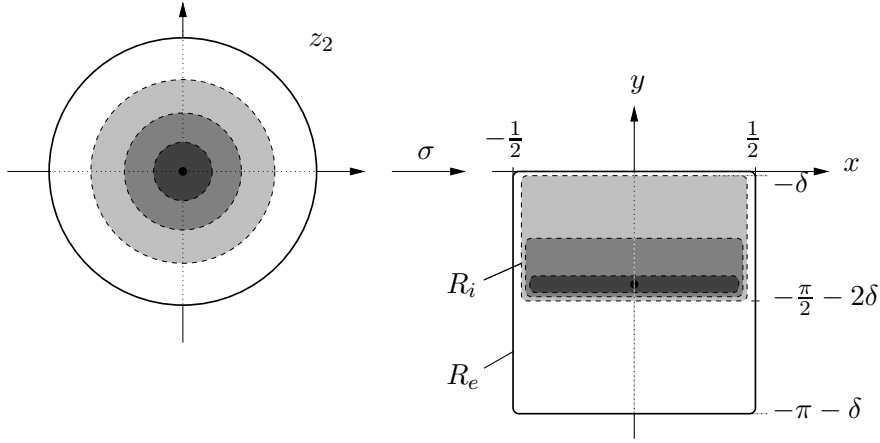}
 \end{center}
 \caption{Preparing the fibres} \label{figure4}
\end{figure}
%
%

Let $\ss$ be a symplectomorphism mapping $D(\pi)$ to $R_e$ and
$D(\frac{\pi}{2})$ to $R_i$ as specified in \mbox{Figure \ref{figure4}}.
We require that for $z_2 \in D(\frac{\pi}{2})$
\begin{equation*}
\pi |z_2|^2 + 2 \dd > y(\ss (z_2)) -
\left( - \frac{\pi}{2} - 2 \dd \right),
\end{equation*}
i.e.\
\begin{equation}  \label{equation6}
y(\ss (z_2)) < \pi |z_2|^2 - \frac{\pi}{2} \quad
\mbox{ for } z_2 \in D \left(\frac{\pi}{2} \right).
\end{equation}

Write for this bundle of round squares 
\[
(id \times \ss ) \cf = \cs = \cs (P_1) \coprod \cs (L) \coprod \cs (P_2).
\]

\medskip

In order to fold $\cs (P_2)$ over $\cs (P_1)$ we first move $\cs (P_2)$ along the
$y$-axis and then turn it in the $z_1$-direction over $\cs (P_1)$. 
\\
\\
{\bf Step 3.}
To move $\cs (P_2)$ along the $y$-axis we follow again
\cite[p.\ 355]{LM1}: 

Let $c \in C^\infty (\RR, \RR)$ with $c(\RR) = [ 0, \, 1-\dd ]$ and

\[
c(t) =
\left\{ \begin{array}{ll}
            0,                       &  t \le \frac{a}{2} + 2 \dd \;\mbox{ and } \;t \ge
                                        \frac{a + \pi}{2} + 8\dd \\
            1- \dd,                  &  \frac{a}{2} +3 \dd \le t \le \frac{a+\pi}{2} +
            7 \dd .
        \end{array}
   \right. 
\]
\begin{figure}[h] 
 \begin{center}
  \psfrag{t}{$t$}
  \psfrag{c}{$c(t)$}
  \psfrag{2}{$\frac{a}{2}+2\dd$}
  \psfrag{10}{$\frac{a+ \pi}{2} +8 \dd$}
  \psfrag{w}{$1-\dd$}
  \leavevmode\epsfbox{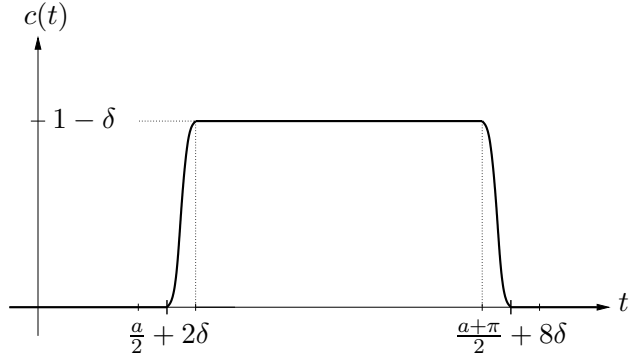}
 \end{center}
 \caption{The cut off $c$} \label{figure5}
\end{figure}
%
%
%
Put $I(t) = \int_0^t c(s) \, ds$ and define $\ff \in \CC^\infty (\RR^4,
\RR^4)$ by

\begin{eqnarray} \label{equation8}
\ff (u, x, v, y) = \left (u, x, v + c(u) \left (x + \frac{1}{2}
\right), y + I(u) \right).
\end{eqnarray}
We then find
\[
d \ff (u,x,v,y)  =  \left[
                       \begin{array}{cc}   
                              {\II}_2   &   0     \\
                              A            &   {\II}_2  
                       \end{array}   \right]          
\; \mbox{ with }              A =   \left[
                                        \begin{array}{cc}   
                                              *     &   c(u)    \\
                                              c(u)  &   0 
                                        \end{array}   \right],   
\]             
whence $\ff$ is a symplectomorphism. Moreover, with $I_\infty =
I(\frac{a+\pi}{2} + 8 \dd)$,
\begin{eqnarray}  \label{equation9}
\ff \, |_{\left\{ u \le \frac{a}{2} + 2 \dd \right\}} = id \quad 
\mbox{ and } \quad
\ff \, |_{\left\{ u \ge \frac{a + \pi}{2} + 8 \dd \right\}} = id + (0,0,0,I_\infty),
\end{eqnarray}
and assuming that $\dd < \frac{1}{10}$ we compute
\begin{eqnarray}  \label{equation10}
\frac{\pi}{2} + 2 \dd < I_\infty < \frac{\pi}{2} + 5 \dd.
\end{eqnarray}
The first inequality in (\ref{equation10}) implies
\begin{eqnarray} \label{equation11}
\ff (P_2 \times R_i) \cap (\RR^2 \times R_e) = \emptyset .
\end{eqnarray}
{\bf Remark.} $\ff$ is the crucial map of the construction; in fact, it
is the only truly symplectic, i.e.\ not 2-dimensional map. $\ff$ is just
the map which sends the lines $\{ v,x,y  \mbox{ constant}\}$ to the
characteristics of the hypersurface 
\[
(u,x,y) \mapsto \left( u,x, c(u) \left( x + \frac{1}{2} \right), y \right),
\]
which generates (the cut off of) the obvious flow separating $R_i$ from
$R_e$.
\diam
\\
\\
{\bf Step 4.} From (\ref{equation8}), \mbox{Figure \ref{figure3}} and
\mbox{Figure \ref{figure4}} we read off that the projection
of $\ff (\cs)$ onto the $(u,v)$-plane is contained in the union of $U$
with the open set bounded by the graph of $u \mapsto \dd + c(u)$ and
the $u$-axis. Observe that $\dd + c(u) \le 1$.

Define a local embedding $\gg$ of this union into the convex hull of $U$
as follows: On $P_1$ the map is the identity and
on $P_2$ it is the orientation preserving isometry between $P_2$ and
$P_1$ which maps the right edge of $P_2$ to the left edge of $P_1$. In
particular, we have for $z_1 = (u,v) \in P_2$
\begin{equation} \label{equation12}
u( \gg (z_1)) = a + \frac{\pi}{2} + 10 \dd -u .
\end{equation}
On the remaining domain $\gg$ looks as follows:
\begin{figure}[h] 
 \begin{center}
  \psfrag{1}{$1$}
  \psfrag{u}{$u$}
  \psfrag{u'}{$u'$}
  \psfrag{a}{$a$}
  \psfrag{b}{$b$}
  \psfrag{c}{$c$}
  \psfrag{d}{$d$}
  \psfrag{a'}{$a'$}
  \psfrag{b'}{$b'$}
  \psfrag{c'}{$c'$}
  \psfrag{d'}{$d'$}
  \psfrag{gg}{$\gg$}
  \leavevmode\epsfbox{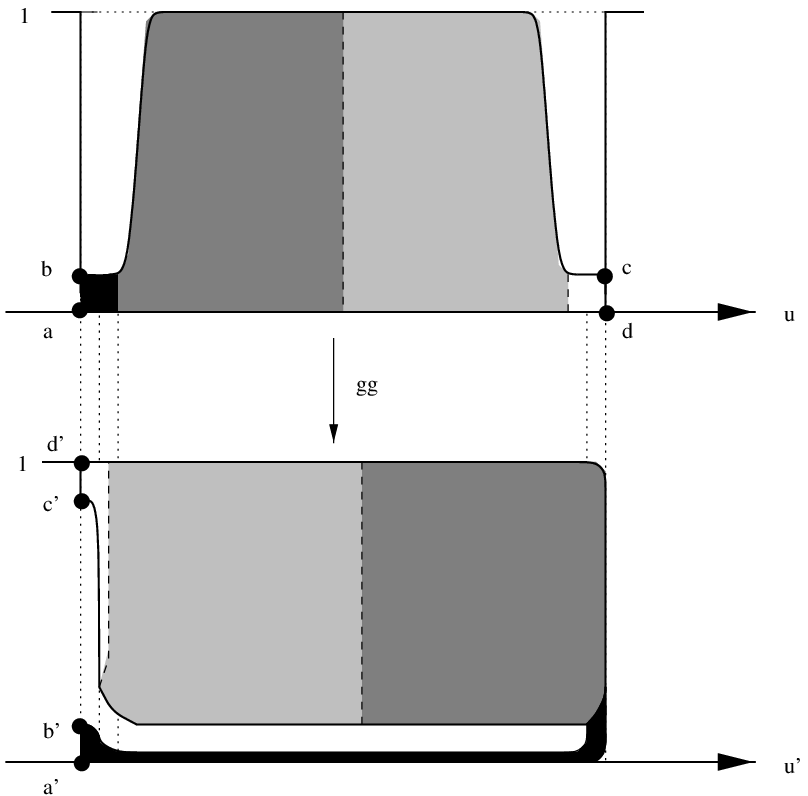}
 \end{center}
 \caption{Folding} \label{figure6}
\end{figure}
%
%
%
In a $\frac{\dd}{4}$-collar of the line from $a$ to $b$ the map is the
identity and on a $\frac{\dd}{4}$-collar of the line from $c$ to $d$ the
linear extension of the map on $P_2$, and we require
\begin{eqnarray*} 
 u'(\gg (u,v)) - \left( \frac{a}{2} + \dd \right) < \frac{\pi}{2} +8\dd -
       \left(u- \left( \frac{a}{2}+\dd \right)\right) +2 \dd,
\end{eqnarray*} 
i.e.
\begin{eqnarray}  \label{equation13}
u'(\gg (u,v)) < -u+ \frac{\pi}{2} + a + 12 \dd .
\end{eqnarray}
(\ref{equation11}) shows that $\gg \times id$ is one-to-one on $\ff (\cs
).$ 
\\
\\
{\bf Step 5.} We finally adjust the fibers: 
 
First of all observe that the projection of $\ff ( \cs )$ onto the
$z_2$-plane is contained in a tower shaped domain $\ct$ 
(cf.\ \mbox{Figure \ref{figure8}})
and that by the second inequality in (\ref{equation10}) we have $\ct
\subset \{ y< \frac{\pi}{2} + 4 \dd \}$. 
\begin{figure}[h] 
 \begin{center}
  \psfrag{x}{$x$}
  \psfrag{y}{$y$}
  \psfrag{-pi-d}{$-\pi-\dd$}
  \psfrag{-pi/2-2d}{$-\frac{\pi}{2}-2\dd$}
  \psfrag{pi/2+4d}{$\frac{\pi}{2}+4\dd$}
  \psfrag{t}{$\ct$}
  \leavevmode\epsfbox{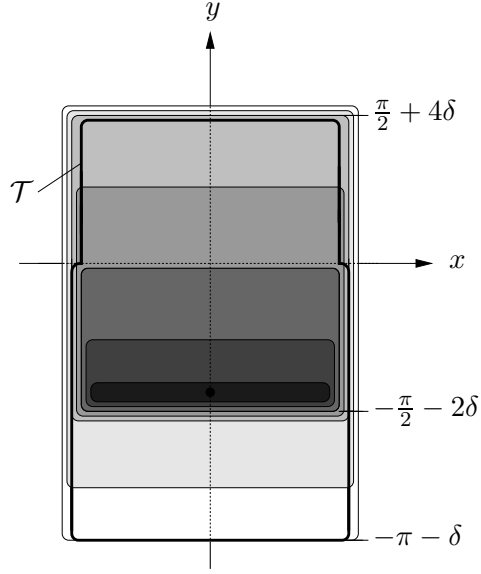}
 \end{center}
 \caption{Mapping the tower to a disc} \label{figure8}
\end{figure}
%
%

We define a symplectomorphism $\tau$ from a neighbourhood of $\ct$ to a
disc by prescribing the preimages of concentric circles as in \mbox{Figure
\ref{figure8}}: We require

\begin{align}
\bullet \;\; & \pi |\tau (z_2)|^2 < y + \frac{\pi}{2}+ 3\dd \quad 
\mbox{ for }  y \ge -\frac{\pi}{2} -2\dd          \label{equation14} \\
\bullet \;\; & \pi |\tau (z_2)|^2 <  \pi | \ss^{-1}(z_2)|^2 
+ \frac{\pi}{2} + 8 \dd \quad \mbox{ for } z_2 \in R_e.  \label{equation15} 
\end{align}

\begin{figure}[h] 
 \begin{center}
  \psfrag{u}{$u$}
  \psfrag{a}{$a$}
  \psfrag{A}{$A$}
  \psfrag{pi}{$\pi$}
  \leavevmode\epsfbox{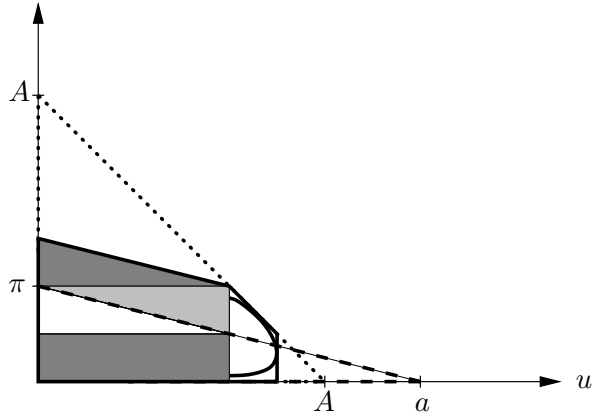}
 \end{center}
 \caption{Folding an ellipsoid into a ball} \label{figure9}
\end{figure}

This finishes the construction. We think of the result as depicted in \mbox{Figure
\ref{figure9}}.
%
%
\\

Let now $\ee >0$ arbitrary and choose $\dd = \min \{ \frac{1}{10}, \, \frac{\ee}{14} \}$.
It remains to check that
\[
\Psi \de (\gg \times \tau ) \circ \ff \circ (\bb \times \ss )
\]
satisfies (\ref{equation4}). So let $z= (z_1, z_2) = (u,v,x,y) \in T(a,
\pi)$ and write $\Psi (z) = (u',v',z_2')$. We have to show that
\begin{eqnarray}  \label{equation16}
u' - \frac{\pi u}{a} + \pi |z_2'|^2 - \pi |z_2|^2 < \frac{a}{2} + 14 \dd .
\end{eqnarray}
\\
{\bf Case 1.} $\bb(z_1) \in P_1$:

(\ref{equation5}) implies $u < \frac{a}{2} + \dd$,  and by
(\ref{equation9}) and step 4 we have $\ff = id$ and $\gg =
id$, whence (\ref{equation5}) and (\ref{equation15}) give
\begin{gather*}
u' = u'(\bb(u,v)) < u + 2 \dd,
\\
\pi |z_2'|^2 = \pi | \tau (\ss (z_2))|^2 < \pi |z_2|^2 + \frac{\pi}{2} +8 \dd .
\end{gather*}
Therefore 
\begin{eqnarray*}
u'-\frac{\pi u}{a}+\pi |z_2'|^2-\pi |z_2|^2  & < &
                                   u \left( 1-\frac{\pi}{a}
                                   \right) +2 \dd + \frac{\pi}{2} +8 \dd    \\ 
                                             & < &
                                   \frac{a}{2} \left( 1-\frac{\pi}{a}
                                   \right)  +\frac{\pi}{2}+11 \dd \\
                                             & = & \frac{a}{2} + 11 \dd .
\end{eqnarray*}

\noindent
{\bf Case 2.} $\bb(z_1) \in P_2$:

Step 2 shows $\ss (z_2) \in R_i$, by (\ref{equation9}) we have 
$\ff = id + (0,0,0,I_\infty)$, and  
(\ref{equation5}) implies $u > \frac{a}{2} - \dd$ and  $u(\bb(z_1)) +
2\dd \ge u + \frac{\pi}{2} + 10\dd$, whence by (\ref{equation12})
\[
u'=u'(\gg(\bb(z_1))) = a + \frac{\pi}{2} + 10 \dd - u(\bb(z_1)) \le
a-u+2 \dd.
\]
Moreover, from (\ref{equation14}), (\ref{equation6}) and (\ref{equation10}) we
see
\begin{eqnarray*}
\pi |z_2'|^2   & = & \pi |\tau (\ss (z_2) + (0, I_\infty))|^2   \\
               & < & y (\ss(z_2)) + I_\infty + \frac{\pi}{2} + 3\dd    \\
               & < &  \pi |z_2|^2 - \frac{\pi}{2} + \frac{\pi}{2} 
                          + 5 \dd + \frac{\pi}{2} + 3 \dd  \\
               & < & \pi |z_2|^2 + \frac{\pi}{2} + 8 \dd.
\end{eqnarray*}

Therefore
\begin{eqnarray*}
u'-\frac{\pi u}{a}+\pi |z_2'|^2-\pi |z_2|^2  & < & a - u \left( 1 +
                                              \frac{\pi}{a} \right) + 
                                              2\dd + \frac{\pi}{2} + 8 \dd  \\
               & < & a - \frac{a}{2} \left( 1 + \frac{\pi}{a} \right) +
                                              \frac{\pi}{2} + 12 \dd    \\
               & = & \frac{a}{2} + 12 \dd .                                                 
\end{eqnarray*}

\noindent
{\bf Case 3.} $\bb(z_1) \in L$:

By construction we have $\ss (z_2) \in R_i$, and 
using the definition of $\ff$, inequality (\ref{equation13}) implies
\[
u' < - u(\bb(u,v)) + \frac{\pi}{2} + a + 12 \dd.
\]
Next (\ref{equation14}), (\ref{equation6}) and the estimate $I(t) <
(1-\dd)(t- (\frac{a}{2} + 2 \dd))$
give
\begin{eqnarray*}
\pi |z_2'|^2  & < & \pi | \tau (x(\ss (z_2)), y(\ss (z_2)) + I(u(\bb
                             (u,v))))|^2  \\
              & < & y (\ss (z_2))
                            + I(u(\bb(u,v)) + \frac{\pi}{2} + 3\dd     \\   
              & < & \pi |z_2|^2 - \frac{\pi}{2} +
                            (1-\dd ) \left( u(\bb(u,v)) - \frac{a}{2} - 2 \dd \right)  
                            + \frac{\pi}{2} + 3\dd .
\end{eqnarray*}
Moreover, (\ref{equation5}) shows $\frac{a}{2} - \dd < u < \frac{a}{2} + \dd$,
whence $u(\bb(u,v)) > u > \frac{a}{2} - \dd$, and therefore
\begin{eqnarray*}
u'-\frac{\pi u}{a}+\pi |z_2'|^2-\pi |z_2|^2  & < & - u(\bb (u,v)) +
                         \frac{\pi}{2} + a + 12 \dd - \frac{\pi}{a}
                         \left( \frac{a}{2} -\dd \right)   \\
                                             &   & + u(\bb (u,v)) -
                         \frac{a}{2} -2 \dd - \dd \left( \frac{a}{2} - \dd \right) +
                         \frac{a}{2} \dd + 2 \dd^2 + 3 \dd        \\
                                             & = & \frac{a}{2} + 13 \dd
                         + \frac{\pi}{a} \dd + 3 \dd^2    \\
                                             & < & \frac{a}{2} + 14 \dd .
\end{eqnarray*}
\proofend

\subsection{Folding in four dimensions} \label{foldinginfourdimensions}

In four dimensions we may exploit the great flexibility of symplectic
maps which only depend on the fibre coordinates to provide rather
satisfactory embedding results for simple shapes.

We first discuss a modification of the folding construction described in
the previous section, then explain multiple folding and finally
calculate the optimal embeddings of ellipsoids and polydiscs into balls
and cubes which can be obtained by these methods. 

Not to disturb the exposition furthermore with $\dd$-terms we skip them
in the sequel. Since all sets considered will be bounded and all
constructions will involve only finitely many steps, we won't lose
control of them.

\subsubsection{The folding construction in four dimensions}
\label{thefoldingconstructioninfourdimensions}

\begin{figure}[h] 
 \begin{center}
  \psfrag{z2}{$z_2$}
  \psfrag{ss}{$\ss$}
  \psfrag{-1/2}{$-\frac{1}{2}$}
  \psfrag{1/2}{$\frac{1}{2}$}
  \psfrag{pi/2}{$\frac{\pi}{2}$}
  \psfrag{x}{$x$}
  \psfrag{y}{$y$}
  \leavevmode\epsfbox{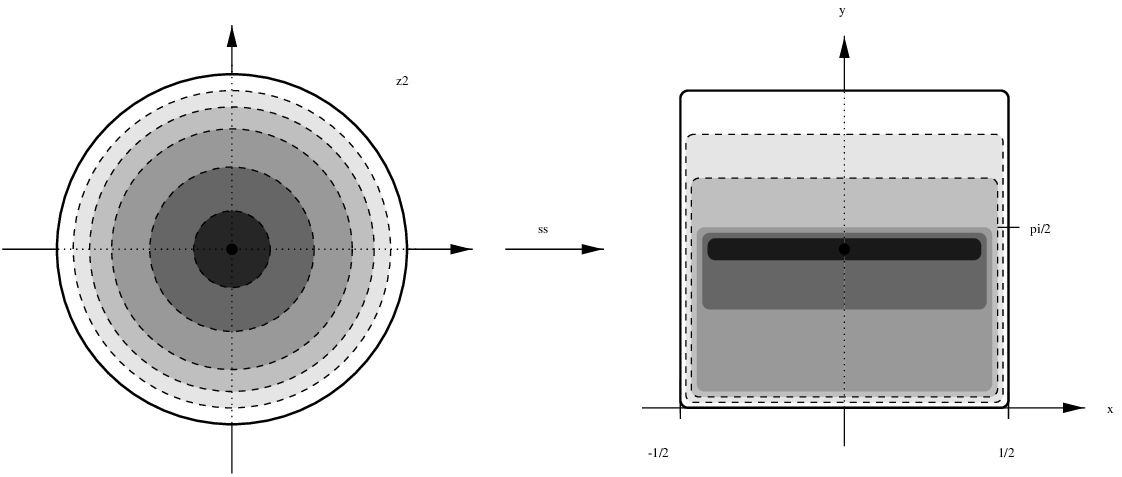}
 \end{center}
 \caption{The modified map $\ss$} \label{figure25}
\end{figure}

\begin{figure}[h] 
 \begin{center}
  \psfrag{u}{$u$}
  \psfrag{a}{$a$}
  \psfrag{A}{$A$}
  \psfrag{pi}{$\pi$}
  \leavevmode\epsfbox{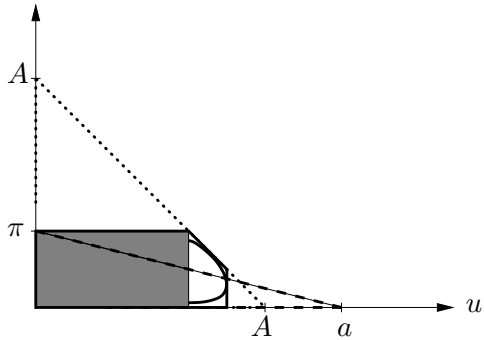}
 \end{center}
 \caption{Folding in four dimensions} \label{figure22}
\end{figure}

The map $\ss$ in step 2 of the folding construction given in the
previous section was dictated by the estimate (\ref{equation16})
necessary for the $n$-dimensional result. 
As a consequence, the map $\gg$ had to disjoin the $z_2$-projection of
$P_2$ from the one of $P_1$, and we ended up with the unused white
sandwiched triangle in \mbox{Figure \ref{figure9}}. In order to use this room
as well we modify the construction as follows:

Replace the map $\ss$ of step 2 by the map $\ss$ given by \mbox{Figure
\ref{figure25}}.
%
%
If we define $\ff$ as in (\ref{equation8}), the $z_2$-projection of the
image of $\ff$ will almost coincide with the image of $\ss$. Choose now
$\gg$ as in step 4 and define the final map $\tau$ on a neighbourhood of
the image of $\ff$ such that it restricts to $\ss^{-1}$ on the image of
$\ss$. If all the $\dd$'s were chosen appropriately, the composite map
$\Psi$ will be one-to-one, and the image $\Psi$ will be contained in
$T(a/2 + \pi + \ee, \, a/2 + \pi + \ee)$ for some small $\ee$. We think
of the result as depicted in \mbox{Figure \ref{figure22}}.
%
%

\subsubsection{Multiple folding}   \label{multiplefolding}

\begin{figure}[h] 
 \begin{center}
  \psfrag{F1}{$F_1$}
  \psfrag{F2}{$F_2$}
  \psfrag{F3}{$F_3$}
  \psfrag{F4}{$F_4$}
  \psfrag{S1}{$S_1$}
  \psfrag{S2}{$S_2$}
  \psfrag{S3}{$S_3$}
  \psfrag{bid}{$\bb \times id$}
  \psfrag{ids1}{$id \times \ss_1$}
  \psfrag{f1f2}{$\ff_2 \circ \ff_1$}
  \psfrag{g1id}{$\gg_1 \times id$}
  \psfrag{g2id}{$\gg_2 \times id$}
  \psfrag{g3id}{$\gg_3 \times id$}
  \psfrag{ids2}{$id \times \ss_2$}
  \psfrag{f3}{$\ff_3$}
  \leavevmode\epsfbox{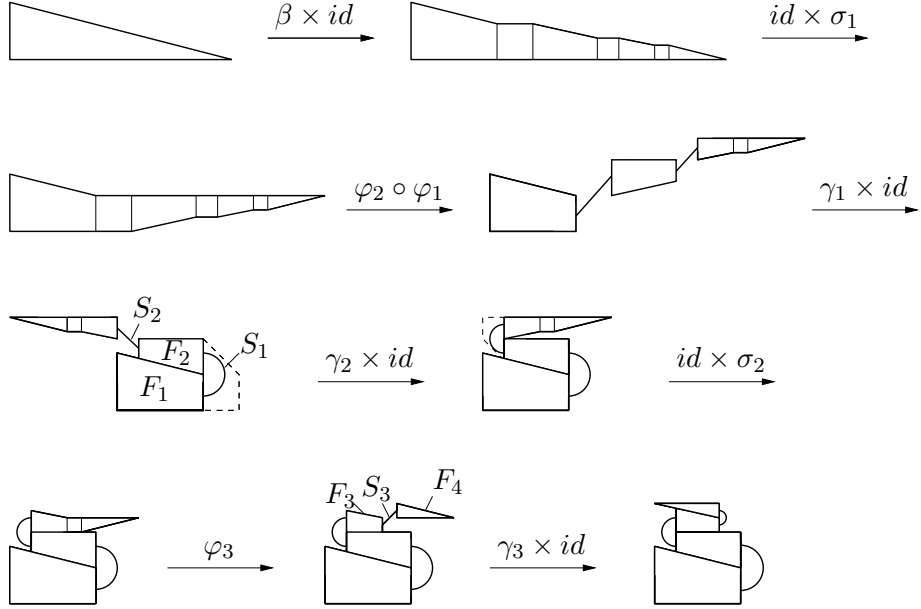}
 \end{center}
 \caption{Multiple folding} \label{figure26}
\end{figure}

\begin{figure}[h] 
 \begin{center}
  \psfrag{u}{$u$}
  \psfrag{v}{$v$}
  \psfrag{u1}{$u_1$}
  \psfrag{u2}{$u_2$}
  \psfrag{u3}{$u_3$}
  \psfrag{u4}{$u_4$}
  \psfrag{pi}{$\pi$}
  \psfrag{1}{$1$}
  \psfrag{l1}{$l_1$}
  \psfrag{l2}{$l_2$}
  \psfrag{l3}{$l_3$}  
  \psfrag{L1}{$L_1$}
  \psfrag{L2}{$L_2$}
  \psfrag{L3}{$L_3$}
  \psfrag{P1}{$P_1$}
  \psfrag{P2}{$P_2$}
  \psfrag{P3}{$P_3$}
  \psfrag{P4}{$P_4$}
  \psfrag{U}{$U$}
  \leavevmode\epsfbox{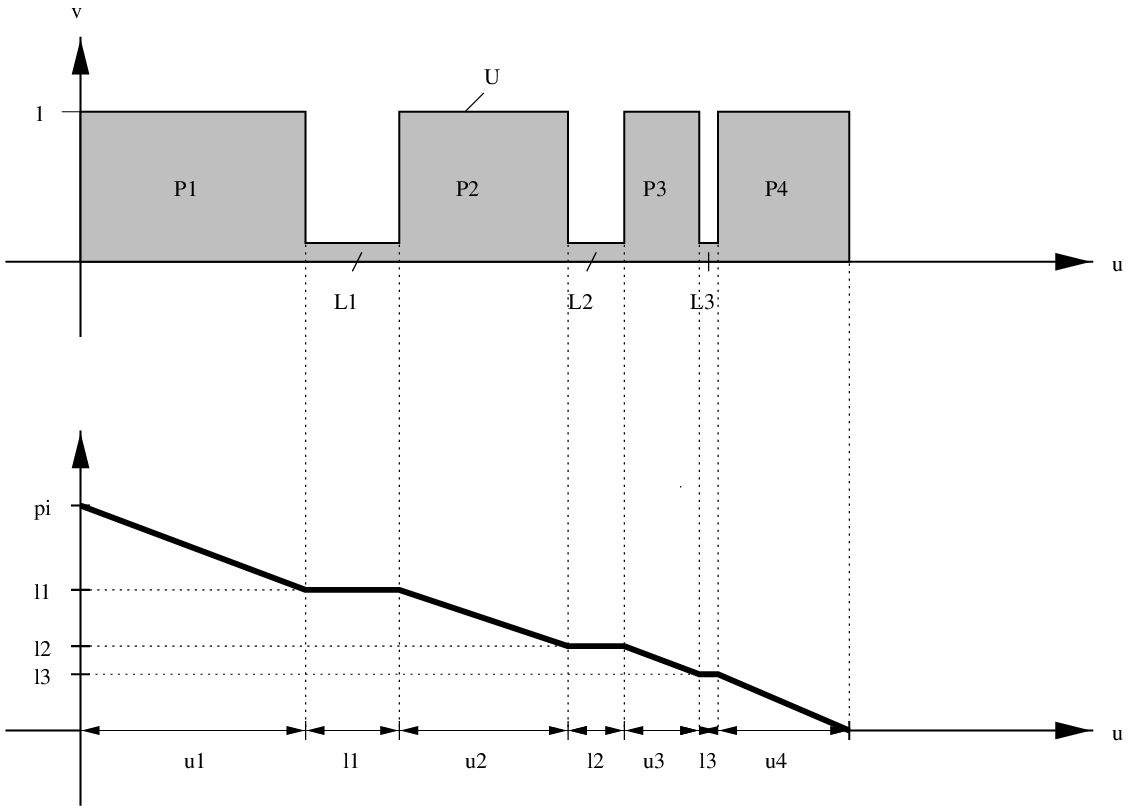}
 \end{center}
 \caption{$\bb$} \label{figure10}
\end{figure}

\begin{figure}[h] 
 \begin{center}
  \psfrag{x}{$x$}
  \psfrag{y}{$y$}
  \psfrag{-1/2}{$-\frac{1}{2}$}
  \psfrag{1/2}{$\frac{1}{2}$}
  \psfrag{ff1}{$\ff_1$}
  \psfrag{ff2}{$\ff_2$}
  \psfrag{l1}{$l_1$}
  \psfrag{2l1}{$2l_1$}
  \psfrag{l1-l2}{$l_1-l_2$}
  \psfrag{2l1+l2}{$2l_1+l_2$}
  \leavevmode\epsfbox{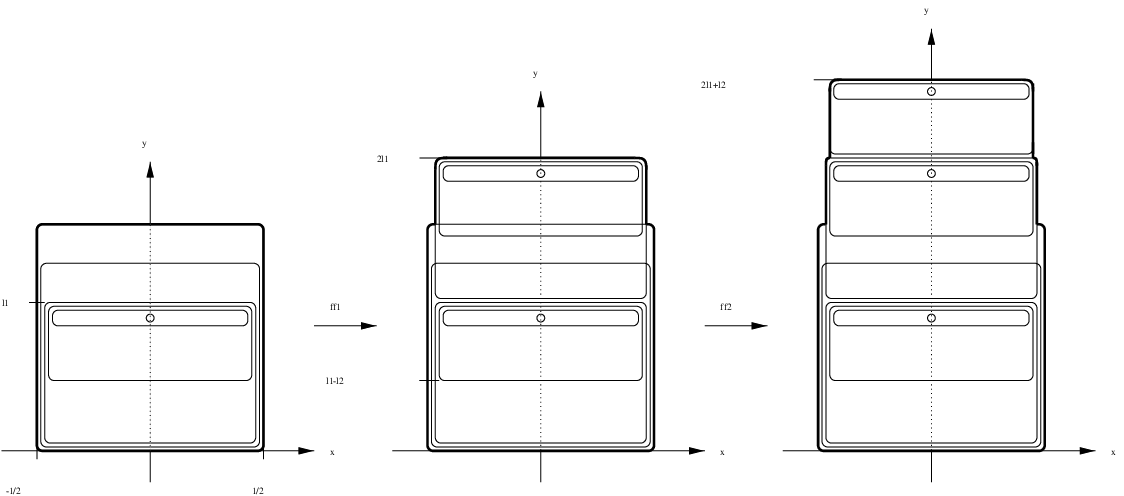}
 \end{center}
 \caption{The first and the second lift} \label{figure23}
\end{figure}

\begin{figure}[h] 
 \begin{center}
  \psfrag{gg}{$\gg_2$}
  \leavevmode\epsfbox{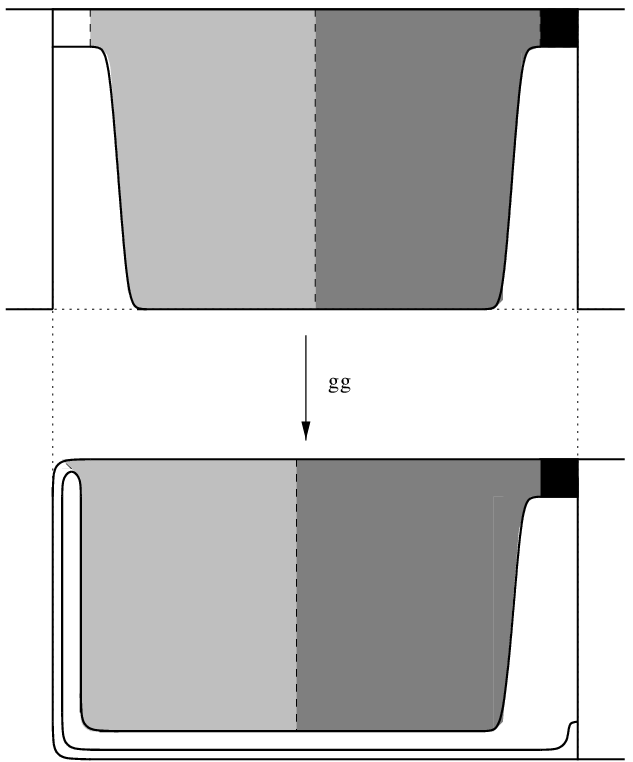}
 \end{center}
 \caption{Folding on the left}  \label{figure30}
\end{figure}

Neither \mbox{Theorem 2} nor Traynor's theorem stated at the beginning of
\mbox{section \ref{flexibility}} tells us if $E(\pi, 4 \pi)$ embeds in $B^4(a)$
for some $a \le 3 \pi$ (cf.\ \mbox{Figure \ref{figure18}}). Multiple
folding, which is explained in this subsection, will provide better embeddings. 
To understand the general construction it is enough to look at a
3-fold:
The folding map $\Psi$ is the composition of maps explained in \mbox{Figure
\ref{figure26}}.
%
%
Here are the details:
Pick reasonable $u_1, \dots, u_4 \in \RR_{>0}$ with $\sum_{j=1}^4 u_i
=a$ and put 
\begin{equation} \label{equation17}
l_i = \pi - \frac{\pi}{a} \sum_{j=1}^i u_j, \qquad i=1,2,3.
\end{equation}
{\bf Step 1.} 
Define $\bb \colon R(a) \ra U$ by \mbox{Figure \ref{figure10}}.
\\
%
%
\\
{\bf Step 2.} 
For $l_1 = \pi /2$ the map $\ss_1$ is given by \mbox{Figure \ref{figure25}},
and in general it is defined to be the symplectomorphism from $D(\pi)$
to the left round rectangle in \mbox{Figure \ref{figure23}}. 
\\
\\
{\bf Step 3.} 
Choose cut offs $c_i$ over $L_i$, $i=1,2$, put $I_i(t) =
\int_0^tc_i(s)\,ds$ and define $\ff_i$ on $\bb \times \ss_1 (T(a,\pi))$
by 
\[
\ff_i (u,x,v,y) = \left( u,x,v+c_i(u) \left(x + \frac{1}{2} \right), y
+ I_i(u) \right).
\]
The effect of $\ff_2 \circ \ff_1$ on the fibres is explained by \mbox{Figure
\ref{figure23}}.
%
%
\\
\\
{\bf Step 4.}
$\gg_1$ is essentially the map $\gg$ of the folding construction: On
$P_1$ it is the identity, for $u_1 \le u \le u_1+l_1$ it looks like the
map in \mbox{Figure \ref{figure6}}, and for $u>u_1+l_1$ it is an
isometry. Observe that by construction, the slope of the stairs $S_2$ is
$1$, while the one of the upper edge of the floor $F_1$ is less than
$1$. $S_2$ and $F_1$ are thus disjoint.
\\
\\
{\bf Step 5.}
$\gg_2 \times id$ is not really a global product map, but restricts to a
product on certain pieces of its domain: It fixes $F_1 \coprod S_1 \coprod
F_2$, and it is the product $\gg_2 \times id$ on the remaining domain where
$\gg_2$ restricts to an isometry on $u_1 \le 0$ and looks like the map
given by \mbox{Figure \ref{figure30}} on the $z_1$-projection of $S_2$.
%
%
\\
\\
For further reference, we summarize the result of the two preceding
steps in the
\\
\\
{\bf Folding Lemma.} 
{\it
Let $S$ be the stairs connecting two floors of minimal respectively
maximal height $l$.
\begin{itemize}
\item[(i)]
If the floors have been folded on top of each other by folding on the
right, $S$ is contained in a trapezoid with horizontal lower edge of
length $l$ and left respectively right edge of length $2l$ respectively $l$.
\item[(ii)]
If the floors have been folded on top of each other by folding on the
left, $S$ is contained in a trapezoid with horizontal upper edge of
length $l$ and left respectively right edge of length $l$ respectively $2l$.
\end{itemize}
}
\noindent
The remaining three maps are restrictions to the relevant parts of
already considered maps.
\\
\\
{\bf Step 6.}
On $\{ y > 2l_1 \}$ the map $\ss_2$ is the automorphism whose image is
described by the same scheme as the image of $\ss_1$, and $id \times
\ss_2$ restricts to the identity everywhere else.
\\
\\
{\bf Step 7.}
On $\{ y > 2l_1 \}$ the map $\ff_3$ restricts to the usual lift, and it
is the identity everywhere else.
\\
\\
{\bf Step 8.}
Finally, $\gg_3 \times id$ turns $F_4$ over $F_3$. It is an
isometry on $F_4$, looks like the map given by \mbox{Figure \ref{figure6}} on
$S_3$ and restricts to the identity everywhere else.
\\
\\
This finishes the multiple folding construction. 
\begin{figure}[h] 
 \begin{center}
  \psfrag{u}{$u$}
  \psfrag{u1}{$u_1$}
  \psfrag{a}{$a$}
  \psfrag{A}{$A$}
  \psfrag{pi}{$\pi$}
  \psfrag{u1+l1}{$u_1 \! +l_1$}
  \psfrag{l1}{$l_1$}
  \psfrag{l2}{$l_2$}
  \leavevmode\epsfbox{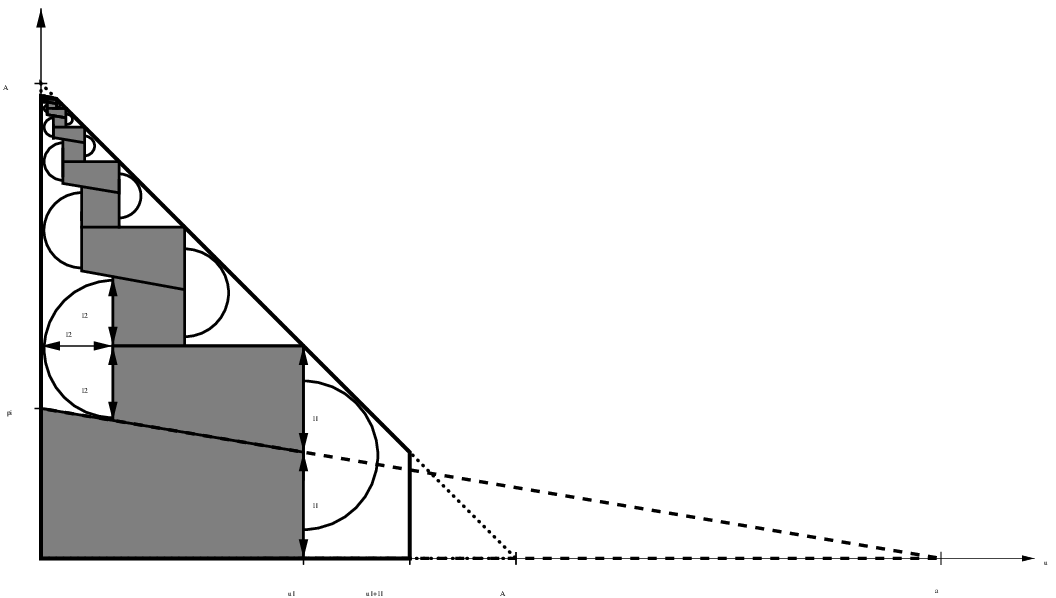}
 \end{center}
 \caption{A 12-fold} \label{figure14}
\end{figure}
%

\subsubsection{Embeddings into balls} \label{embeddingsintoballs}

In this subsection we use multiple folding to construct good embeddings of
ellipsoids into balls, and we also look at embeddings of polydiscs into
balls.

\paragraph{Embedding ellipsoids into balls}
\label{embeddingellipsoidsintoballs}

We now choose the $u_j$'s optimal.

Fix $u_1 >0$. As proposed in \mbox{Figure \ref{figure27}}, 
we assume that the second floor $F_2$ touches
the boundary of \,$T(A,A)$ and that all the other $u_j$'s are chosen
maximal. In other words, $A$ is given by
\begin{equation} \label{equation18}
A(a,u_1) = u_1 + 2\,l_1 = 2 \pi+ \left(1 - \frac{2 \pi}{a} \right) u_1,
\end{equation}
and we proceed as follows: If the remaining length $r_1 = a-u_1$ is
smaller than $u_1$, i.e.\ $u_1 \ge a/2$, we are done; otherwise we try
to fold a second time. By the Folding Lemma, this is possible if and
only if $l_1 < u_1$, i.e.
\begin{equation} \label{equation19}
u_1 > \frac{a \pi}{a + \pi}.
\end{equation}
If (\ref{equation19}) does not hold, the embedding attempt fails; if
(\ref{equation19}) holds, the Folding Lemma and the maximality of $u_2$
imply $u_2 = u_1 - l_2$, whence by (\ref{equation17})
\[
u_2 = \frac{a+ \pi}{a - \pi} u_1 - \frac{a \pi}{a- \pi}.
\]
If the upper left corner of $F_3$ lies outside $T(A,A)$, the embedding
attempt fails, otherwise we go on.

In general, assume that we folded already $j$ times and that $j$ is
even. If the length of the remainder $r_j = r_{j-1} -u_j$ is smaller
than $u_j$, we are done; if not, we try to fold again: The Folding Lemma
and the maximality of $u_{j+1}$ imply $u_{j+1} + 2l_{j+1} = u_j$, and
substituting $l_{j+1} = l_j - u_{j+1} \pi/a$ we get
\[
u_{j+1} = \frac{a}{a-2\pi} (u_j - 2l_j).
\]
If $u_j \le 2l_j$, the embedding attempt fails, otherwise we go on:
If the length of the new remainder $r_{j+1} = r_j - u_{j+1}$ is smaller than
$u_{j+1} + l_j$, we are done; otherwise we try to fold again: The
Folding Lemma and the maximality of $u_{j+2}$ imply 
$u_{j+2} + l_{j+2} = u_{j+1} + l_j$, whence by (\ref{equation17})
\[
u_{j+2} = \frac{a+\pi}{a-\pi} u_{j+1}.
\]
The embedding attempt fails here if and only if the upper left corner of
the floor $F_{j+3}$ lies outside $T(A,A)$; if this does not happen, we
may go on as before.

\medskip
First of all note that whenever the above embedding attempt succeeds, it
indeed describes an embedding of $E(\pi, a)$ into
$T(A(a,u_1),A(a,u_1))$. In fact, it is enough to define the fiber
adjusting map $\tau$ on a small neighbourhood of the resulting tower
$\ct$ in such a way that for any $z_2 = (x,y), \,z_2' = (x',y') \in \ct$ we
have
\[
y \le y' \;\; \Longrightarrow \;\; |\tau (z_2)|^2 < |\tau (z_2')|^2.
\]

(\ref{equation18}) shows that we have to look for the smallest $u_1$ for
which the above embedding attempt succeeds. Call it $u_0 = u_0(a)$. As
we have seen above, $u_0$ lies in the interval
\begin{equation} \label{equation23}
I_a = \left[ \frac{a \pi}{a+\pi},  \frac{a}{2} \right].
\end{equation}
Moreover, if the embedding attempt succeeds for $u_1$, the same clearly
holds true for any $u_1' > u_1$. Hence, given $u_1 \in I_a$, the
corresponding embedding attempt succeeds if and only if $u_1 \ge
u_0$. \mbox{Appendix A1} provides a computer program calculating $u_0$, and the
result $s_{EB}(a) = 2 \pi + (1- 2\pi /a)u_0$
is discussed and compared with
the one yielded by Lagrangian folding in \mbox{subsection 
\ref{symplecticversuslagrangianfolding}}.
\\
\\
{\bf Remarks.}
{\bf 1.\ \!\!\!\!}
Simple geometric considerations show that our choices in the above algorithm
are optimal, i.e.\ $s_{EB}(a)$ provides the best estimate for an
embedding of $E(\pi, a)$ into a ball obtainable by multiple folding.

{\bf 2.\ \!\!\!}
Let $u_1 > u_0$ and let $N(u_1)$ be the number of folds needed in the
above embedding procedure determined by $u_1$. Then $N(u_1) \ra \infty$
as $u_1 \searrow u_0$, i.e.\ the best embeddings are obtained by folding
arbitrarily many times. This follows again from an easy geometric
reasoning.

{\bf 3.\ \!\!\!}
Fix $N$ and let $A_N(a)$ be the function describing the optimal
embedding obtainable by folding $N$ times. Then $\{ A_N \}_{n \in \NN}$ 
is a monoton
decreasing family of rational functions on $[2\pi, \infty[$. For instance, 
\[ 
A_1(a) = 2 \pi + (a-2\pi) \frac{1}{2}, 
\qquad 
A_2(a) = 2 \pi + (a-2\pi) \frac{a+\pi}{3a+\pi}  \]
and
\[ A_3(a) = 2 \pi + (a -2\pi) \frac{(a+\pi)(a+2\pi)}{4 (a^2+a\pi+\pi^2)}. \]
So, $A_1'(2\pi) = \frac{1}{2}$ and $A_2'(2\pi) =A_3'(2\pi) = \frac{3}{7}$.
One can show that $A_N'(2\pi) = \frac{3}{7}$ for all $N \ge 3$. Thus
\[ \limsup_{\ee \ra 0^+} \frac{s_{EB}(2 \pi + \ee) - 2 \pi}{\ee} =
\frac{3}{7}. \]
\hfill{\diam}

\paragraph{Embedding polydiscs into balls}  
               \label{embeddingpolydiscsintoballs}

\begin{proposition}  \label{propositionpolydiscsintoballs}
Let $a>2\pi$ and $\ee >0$. Then $P(\pi,a)$ embeds in
\mbox{$B^4(s_{PB}(a)+\ee)$}, where $s_{PB}$ is given by
\[
s_{PB}(a) = \frac{a-2\pi}{2k}+(k+2)\pi, \qquad 2(k^2-k+1) < a/\pi \le
2(k^2+k+1).
\]
\end{proposition}

\proof
\begin{figure}[h] 
 \begin{center}
  \psfrag{u}{$u$}
  \psfrag{u1}{$u_1$}
  \psfrag{A3}{$A_3$}
  \psfrag{10pi}{$10\pi$}
  \psfrag{pi}{$\pi$}
  \leavevmode\epsfbox{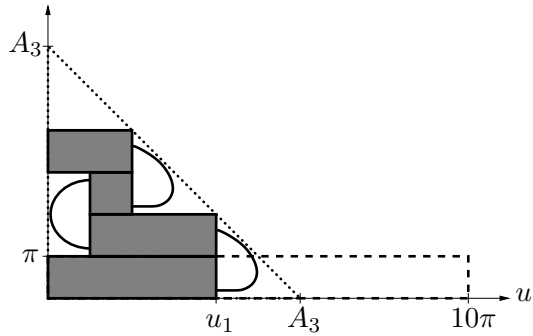}
 \end{center}
 \caption{The optimal embedding $P(\pi, 10 \pi) \hookrightarrow
  B^4(A)$} 
\label{figure16}
\end{figure}
%
%
Let $N=2k-1, \,k\in \NN$, be odd. From \mbox{Figure \ref{figure16}} we read
off that under the condition $u_1>N\pi$ the optimal embedding by folding
$N$ times is described by
\begin{eqnarray*}
a  & = & \pi + 2(u_1-\pi)+2(u_1-3\pi) + \dots + 2(u_1-N\pi) +\pi \\
   & = & 2\pi + 2ku_1 - 2k^2\pi
\end{eqnarray*}
and $A_N(a) = u_1 + 2\pi$; hence
\[
A_N(a) = \frac{a-2\pi}{2k}+(k+2)\pi,
\]
provided that $A_N(a) - 2\pi > (2k-1)\pi$. This condition translates to
$a> 2(k^2-k+1)\pi$, and the claim follows.
\proofend

\noindent
{\bf Remark.} $s_{PB}$ is the optimal result obtainable by multiple
folding. In fact, a simple geometric argument or a similar calculation
as in the proof shows that folding $2k$ times yields worse estimates.
\diam
\begin{figure}[h] 
 \begin{center}
  \psfrag{a/pi}{$\frac{a}{\pi}$}
  \psfrag{A/pi}{$\frac{A}{\pi}$}
  \psfrag{1}{$1$}
  \psfrag{2}{$2$}
  \psfrag{3}{$3$}
  \psfrag{5}{$5$}
  \psfrag{6}{$6$}
  \psfrag{7}{$7$}
  \psfrag{10}{$10$}
  \psfrag{cEH}{$c_{EH}$}
  \psfrag{sPB}{$\frac{s_{PB}}{\pi}$}
  \psfrag{inclusion}{\mbox{inclusion}}
  \psfrag{vol}{\mbox{volume condition}}
  \leavevmode\epsfbox{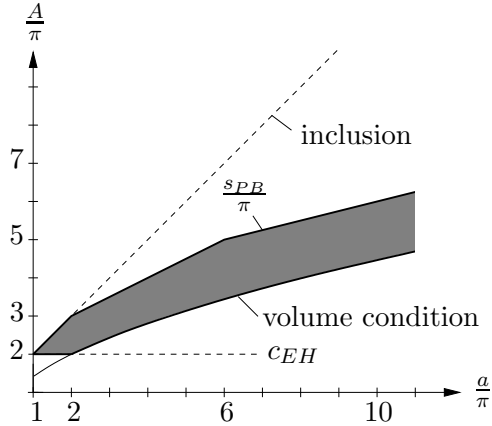}
 \end{center}
 \caption{What is known about $P(\pi, a) \hookrightarrow
  B^4(A)$} 
\label{figure17}
\end{figure}
%
%
\begin{remark}
{\rm
Let $d_{PB}(a) = s_{PB}(a) -\sqrt{2\pi a}$ be the difference between
$s_{PB}$ and the volume condition. $d_{PB}$ attains its local maxima at
$a_k = 2(k^2-k+1)\pi$, where $d_{PB}(a_k) = (2k+1)\pi - 2\pi
\sqrt{k^2-k+1}$. This is an increasing sequence converging to $2\pi$.
\diam
}
\end{remark}

\subsubsection{Embeddings into cubes}
\label{embeddingintocubes}

Given an open set $U$ in $\CC^n$, call the orthogonal projections of $U$
onto the $n$ symplectic coordinate planes the shadows of $U$. As pointed
out in \cite[p.\ 580]{FHW}, symplectic capacities measure to some extent
the areas of the shadows of a set. Of course, this can not be made
rigorous since the areas of shadows are no symplectic invariants, but
for sufficiently regular sets these areas indeed are symplectic
capacities: As remarked before, the capacities $a_1, \dots\!, a_n$ of the
ellipsoid $E(a_1, \dots, a_n)$ are symplectic capacities and, more
generally, given any bounded $U$ with connected smooth boundary
$\partial U$ of restricted contact type and with a shadow whose boundary 
is the shadow of a closed characteristic on $\partial U$ which
lies in a single symplectic coordinate direction, this shadow is a
capacity of $U$ \cite[Proposition 2]{EH2}. Moreover, the smallest shadow
of a polydisc and of a symplectic cylinder are capacities.

Instead of studying embeddings into minimal balls, i.e.\ to reduce the
diameter of a set, it is therefore a more symplectic enterprise to look
for minimal embeddings into a polydisc $C^{2n}(a)$,
i.e.\ to reduce the maximal shadow.

The Non-Squeezing Theorem states that the smallest shadow of simple sets
(like ellipsoids, polydiscs or cylinders) can not be reduced. We therefore
call obstructions to the reduction of the maximal shadow highest order
rigidity. (More generally, calling an ellipsoid or a polydisc given by
$a_1 \le \dots \le a_n$ $i$-reducible if there is an embedding into
$C^{2i}(a') \times \RR^{2n-2i}$ for some $a'<a_i$, one might explore
$i$-th order rigidity.)
\\
\\
The disadvantage of this approach to higher order rigidity is that for a
polydisc there are no good higher invariants available, in fact,
Ekeland-Hofer-capacities see only the smallest shadow \cite[Proposition
5]{EH2}:
\[
c_j(P(a_1, \dots, a_n)) = j a_1.
\]
Many of the polydisc-analogues of the rigidity results for ellipsoids
proved in section \ref{rigidity} are therefore either wrong or much harder to prove. It is
for instance not true that $P(a_1, \dots, a_n)$ embeds linearly in
$P(a_1', \dots, a_n')$ if and only if $a_i \le a_i'$ for all $i$, for a long enough
4-polydisc may be turned into the diagonal of a cube of smaller maximal
shadow:
\begin{lemma}
Let $r>1+ \sqrt{2}$. Then $P^{2n}(\pi, \dots, \pi, \pi r^2)$ embeds
linearly in $C^{2n}(a)$ for some $a< \pi r^2$.
\end{lemma}
\proof
It is clearly enough to prove the lemma for $n=2$. Consider the linear
symplectomorphism given by
\[
(z_1, z_2) \mapsto (z_1', z_2') = \frac{1}{\sqrt{2}} (z_1+z_2, z_1-z_2).
\]
For $(z_1, z_2) \in P(\pi, \pi r^2)$  we have for $i = 1,2$
\begin{eqnarray} \label{equation30}
| z_i'|^2 \le \frac{1}{2} (|z_1|^2 + |z_2|^2 + 2 |z_1| |z_2|) \le
\frac{1}{2} + \frac{r^2}{2} +r,
\end{eqnarray}
and the right hand side of (\ref{equation30}) is strictly smaller than
$r^2$ provided that $r>1 + \sqrt{2}$.
\proofend
\\
Similarly, we don't know how to prove the full analogue of \mbox{Proposition
\ref{proposition1}}: 

Let $\cp (n)$ be the collection of polydiscs
\[
\cp (n) = \{ P(a_1, \dots, a_n) \}
\]
and write $\rp_i$ for the restrictions of the relations $\le_i$ to $\cp
(n)$. Again
$\rp_2$ and $\rp_3$ are very similar, again all the relations $\rp_i$
are clearly reflexive and transitive, and again the identitivity of
$\rp_2$, which again implies the one of $\rp_1$, follows from the
equality of the spectra, which is implied by
the equality of the volumes. (Observe that, even though the boundary of a
polydisc is not smooth, its spectrum is still well defined.)
For $n$=2 the identitivity of $\rp_3$ is
seen by using any symplectic capacity, which determines the smallest
shadow, and the equality of the volumes; but for arbitrary $n$ we don't
know a proof.

\bigskip

While the lack of convenient invariants made it impossible to get good
rigidity results for embeddings into polydiscs, the folding construction
provides us with rather satisfactory flexibility results.

\paragraph{Embedding ellipsoids into cubes}
\label{embeddingellipsoidsintocubes}

\begin{figure}[h] 
 \begin{center}
  \psfrag{u}{$u$}
  \psfrag{u1}{$u_1$}
  \psfrag{7pi}{$7 \pi$}
  \psfrag{A}{$A$}
  \psfrag{pi}{$\pi$}
  \leavevmode\epsfbox{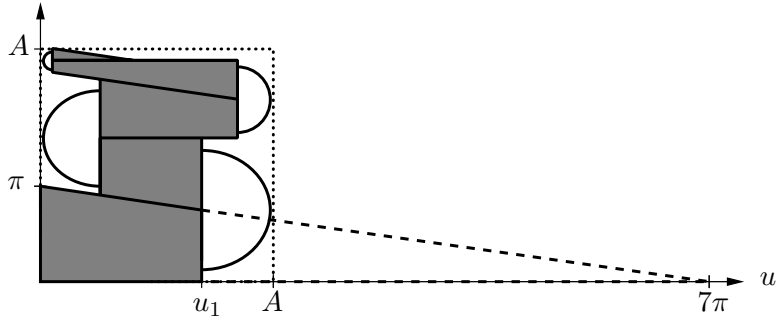}
 \end{center}
 \caption{The optimal embedding $E(\pi, 7 \pi) \hookrightarrow C^4(A)$} \label{figure18a}
\end{figure}

We again use the notation of section \ref{foldinginfourdimensions}, 
fold first at some
reasonable $u_1$ and then choose the subsequent $u_j$'s maximal (see
\mbox{Figure \ref{figure18a}}).
%
%
Let $w(a,u_1) = u_1 + l_1 = \pi+(1-\pi/a)u_1$ be the width of the image
and $h = h(a,u_1)$ its height.

Let's first see what we get by folding once: The only condition on $u_1$
is $a/2 \le u_1$, whence $h(a,u_1) = \pi < \pi + (1-\pi/a)u_1 =
w(a,u_1)$. The optimal choice of $u_1$ is thus $u_1 = a/2$. 

Suppose now that we fold at least twice. The only condition on $u_1$ is
then again $l_1 < u_1$, i.e.
\begin{equation*}  
u_1 > \frac{a \pi}{a+\pi}.
\end{equation*}
Observe that $h(a,u_1)$ diverges if $u_1$ approaches $a\pi / (a+\pi)$. 
Note also that $w$ is increasing in $u_1$ while $h$ is
decreasing. Thus, $w(a,u_1)$ and $h(a,u_1)$ intersect exactly once,
namely in the optimal $u_1$, which we call $u_0$. In particular, we see
that folding only once never yields an optimal embedding.
Write $s_{EC}(a) = \pi+(1-\pi /a)u_0$ for the resulting estimate. It is
computed in \mbox{Appendix A2}.
Again, it is easy to see that our choices in the above procedure are
optimal, i.e.\ $s_{EC}(a)$ provides the best estimate for an embedding
of $E(\pi,a)$ into a cube obtainable by symplectic folding.
\\
\\
{\bf Example.}
If we fold exactly twice, we have $h = 2l_1+l_2$, or, since $l_2$
satisfies $a = u_1+ u_2 +(a/\pi) l_2$ and $u_2 = u_1-l_2$, 
\[
h = 2\pi - \frac{2\pi}{a}u_1+\frac{\pi(a-2u_1)}{a-\pi}.
\]
Thus, provided that $l_2 + (a/\pi) l_2 \le w$, the equation $h=w$ yields
\begin{equation}  \label{equation32}
u_0 = \frac{a \pi (2a-\pi)}{a^2+2a\pi-\pi^2}.
\end{equation}
Indeed, $u_0$ satisfies (\ref{equation32}) whenever $a>\pi$. Finally,
$l_2+(a/\pi) l_2 \le w$ holds if and only if $\pi \le a \le 3\pi$. 
\diam

\begin{figure}[h] 
 \begin{center}
  \psfrag{1}{$1$}
  \psfrag{2}{$2$}
  \psfrag{3}{$3$}
  \psfrag{4}{$4$}
  \psfrag{5}{$5$}
  \psfrag{6}{$6$}
  \psfrag{7}{$7$}
  \psfrag{A/pi}{$\frac{A}{\pi}$}
  \psfrag{a/pi}{$\frac{a}{\pi}$}
  \psfrag{sec/pi}{$\frac{s_{EC}}{\pi}$}
  \psfrag{inclusion}{$\mbox{inclusion}$}
  \psfrag{vol}{$\mbox{volume condition}$}
  \psfrag{cEH}{$c_{EH}$}
  \leavevmode\epsfbox{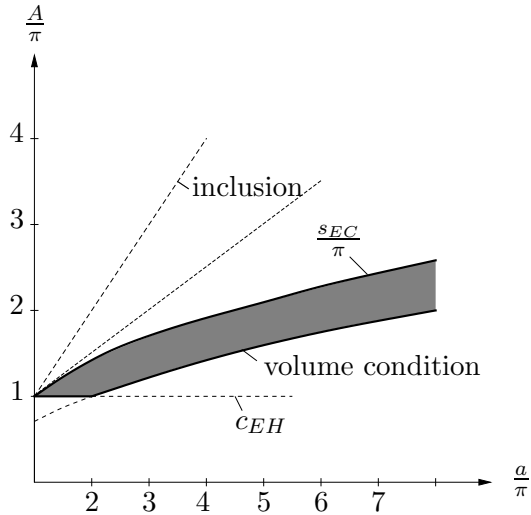}
 \end{center}
 \caption{What is known about $E(\pi,a) \hookrightarrow C^4(A)$} \label{figure19}
\end{figure}
%

In fact, (\ref{equation32}) also holds true for all $a$ for which the
optimal embedding of $E(\pi, a)$ obtainable by multiple folding is a
3-fold for which the height is still described by $h = 2l_1+l_2$, i.e.\
for which $u_4 \le u_3$. This happens for $3 < a/\pi < 4.2360\dots$,
whence
\[
s_{EC}(a) = \frac{a\pi (3a-\pi)}{a^2+2a\pi-\pi^2} \quad \mbox{ for } 1 \le
\frac{a}{\pi} \le 4.2360\dots.
\]
In general, $s_{EC}$ is a piecewise rational function. Its singularities
are those $a$ for which $u_{N(a)} = u_{N(a)+1}$, where we wrote $N(a)$
for the number of folds determined by $u_0(a)$.

\begin{remark}
\label{remarkdec}
{\rm
Let $d_{EC}(a) = s_{EC}(a) - \sqrt{\pi a/2}$ be the difference between
$s_{EC}$ and the volume condition. The set of local minima of $d_{EC}$
coincides with its singular set, i.e.\ with the singular set of
$s_{EC}$.
On the other hand, $d_{EC}$ attains its local maxima at those $a$ for
which the point of $F_{N(a)+1}$ touches the boundary of
$T(A,A)$. Computer calculations suggest that on this set, $d_{EC}$ is
increasing, but bounded by $(2/3) \pi$.
\diam
}
\end{remark}


\paragraph{Embedding polydiscs into cubes}
                       \label{embeddingpolydiscsintocubes}

\begin{proposition} \label{propositionpintocube}

Let $a>2 \pi$ and $\ee >0$. Then $P(\pi,a)$ embeds in \mbox{$C^4(s_{PC}(a)+\ee)$},
where $s_{PC}$ is given by
\[
s_{PC}(a) = 
\left\{ \begin{array}{cl}
           (N+1)\pi,             & (N-1)N +2 < \frac{a}{\pi} \le N^2+1 \\
           \frac{a+2N \pi}{N+1}, & N^2 +1    < \frac{a}{\pi} \le N(N+1)+2\,
           .
        \end{array}
   \right. 
\]
\end{proposition}

\proof
The optimal embedding by folding $N$ times is described by 
\[
2 u_1 + (N-1)(u_1 - \pi) = a, 
\]
whence $u_1 = \frac{a+ (N-1) \pi}{N+1}$; in fact, by the assumption on
$a$, the
only condition $u_1 > \pi$ for $N \ge 2$ is satisfied. Thus
$A_N(a) = \max \{ \frac{a+ 2N \pi}{N+1} , (N+1) \pi \}$, and the
proposition follows.
\proofend
\begin{figure}[h] 
 \begin{center}
  \psfrag{u}{$u$}
  \psfrag{u1}{$u_1$}
  \psfrag{u1+pi}{$u_1 \! + \pi$}
  \psfrag{a}{$a$}
  \psfrag{pi}{$\pi$}
  \leavevmode\epsfbox{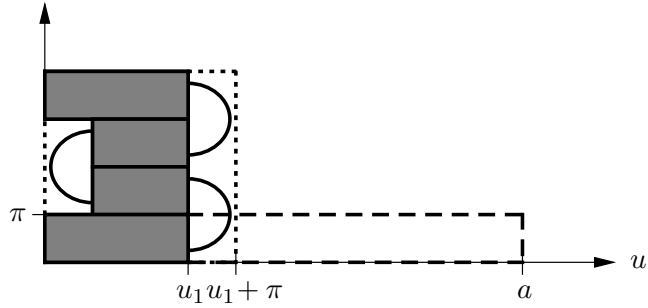}
 \end{center}
 \caption{Folding $P(\pi,a)$ three times} \label{figure20}
\end{figure}
%
%
       
\begin{remark}  \label{remarkdpc}
{\rm
The difference $d_{PC}(a) = s_{PC}(a) - \sqrt{\pi a}$ between $s_{PC}$
and the volume condition attains its local maxima at $a_N =
(N^2-N+2)\pi$, where $d_{PC}(a_N) = (N+1)\pi -
\sqrt{N^2-N+2}\,\pi$. This is an increasing sequence converging to
$(3/2)\pi$.
\diam
}
\end{remark}

Since for $a \le 2 \pi$ 
folding cannot reduce $P(\pi, a)$ and since we believe that for
small $a$ folding is essentially the only way to achieve a reduction
(see also \cite{LM3}), we state:

\begin{conjecture}
The polydisc-analogue of \mbox{Theorem 1'} holds. In particular,
\newline $P^{2n}(\pi,
\dots, \pi, a)$ embeds in $C^{2n}(A)$ for some $A<a$ if and only if $a>
2 \pi$.
\end{conjecture}

\subsection{Folding in higher dimensions}
\label{foldinginhigherdimensions}

Even though symplectic folding is an essentially four dimensional
process, we may still use it to get good embeddings in higher dimensions
as well. The point is that we may fold into different symplectic
directions of the fiber.
In view of the applications of higher dimensional folding in subsection
\ref{asymptoticpackings} and \ref{refinedasymptoticinvariants} we will
concentrate on embedding skinny
polydiscs into cubes and skinny ellipsoids into balls and cubes.

\bigskip

Given domains $U \subset \RR^{2n}$ and $V, W \subset \RR^n$ and given
$\aa >0$, we set 
\[
\aa U = \{ \aa z \in \RR^{2n} \,|\, z \in U \} 
\quad \mbox{ and } \quad
{\aa V} \times W = \aa (V \times W).
\]
As in the four dimensional case we may view an ellipsoid $E(a_1, \dots,
a_n)$ as fibered over the disc $D(a_n)$ with ellipsoids $\gg E(a_1,
\dots,a_{n-1})$ of varying size as fibres. By deforming the base $D(a_n)$ to
a rectangle as in \mbox{Figure \ref{figure2}} we may get rid of the
$y_1$-coordinate. It will be convenient to get rid of the other
$y_i$-coordinates too. 
Write $\RR^{2n}(x,y) = \RR^n(x) \times \RR^n(y)$ and set
\begin{gather*}
\tr (a_1, \dots, a_n) = \{ 0 < x_1, \dots, x_n \, | \, \sum_{i=1}^n
\frac{x_i}{a_i} <1 \} \subset \RR^n(x), 
\\
\sq (b_1, \dots, b_n) = \{ 0 < y_i < b_i, \; 1 \le i \le n \} \subset
\RR^n(y).
\end{gather*}
\begin{lemma}    \label{lemmafibres}
For all $\ee >0,$
\begin{itemize}
\item[(i)]
$E(a_1-\ee, \dots, a_n-\ee)$ embeds in $\tr (a_1, \dots, a_n) \times
\sq^n (1)$ in such a way that for all $\aa \in \; ]0,1]$, $\aa E(a_1-\ee,
\dots, a_n-\ee)$ is mapped into $(\aa +\ee) \tr(a_1, \dots,a_n) \times
\sq^n (1)$.
\item[(ii)]
$\tr (a_1-\ee, \dots, a_n-\ee) \times \sq^n (1)$ embeds in $E(a_1, \dots, a_n)$
in such a way that for all $\aa \in \; ]0,1]$, 
$\aa \tr (a_1-\ee, \dots, a_n-\ee) \times \sq^n (1)$ is mapped into
$(\aa +\ee) E(a_1, \dots, a_n)$.
\end{itemize}
\end{lemma}
\proof
By Lemma \ref{lemmaarea} we find embeddings $\aa_i \colon D(a_i-\ee)
\hookrightarrow \sq(a_i,1)$ satisfying
\[
x_i(\aa_i(z_i)) < \pi |z_i|^2 + \frac{\ee}{n} \frac{a_1}{\max(1,a_n)} \quad
\mbox{ for }  z_i \in D(a_i - \ee), \; 1 \le i \le n
\]
(cf.\ \mbox{Figure \ref{figure2}}). 
Given $(z_1, \dots, z_n) \in E(a_1 -\ee, \dots, a_n-\ee)$ we then find
\begin{eqnarray*}
\sum_{i=1}^n \frac{x_i(\aa_i(z_i))}{a_i} & < & \sum_{i=1}^n \frac{\pi
                                               |z_i|^2}{a_i} +
                                               \frac{1}{a_i}
                                               \frac{\ee}{n} \frac{a_1}{a_n}\\
               & < & \max_i \frac{a_i -\ee}{a_i} + \frac{\ee}{a_n} 
                 \;=\;  1 - \frac{\ee}{a_n} +\frac{\ee}{a_n} \;=\; 1,   
\end{eqnarray*}
and given $(z_1, \dots, z_n) \in \aa E(a_1 -\ee, \dots, a_n-\ee)$ we find
\begin{eqnarray*}
\sum_{i=1}^n \frac{x_i(\aa_i(z_i))}{a_i} & < & \sum_{i=1}^n \frac{\pi
                                               |z_i|^2}{a_i} +
                                               \frac{a_1}{a_i}
                                               \frac{\ee}{n}    
                \;<\;  \aa + \ee.   
\end{eqnarray*}

The proof of (ii), which uses products of maps $\oo_i$ as in \mbox{Figure 
\ref{figure2}}, is similar.
\proofend

Forgetting about all the $\ee$'s, we may thus view an ellipsoid as a
Lagrangian product of a simplex and a cube. 
In the setting of symplectic folding, however, we will still rather
think of $E(a_1, \dots, a_n)$ as fibered over the base $\sq (a_n, 1)$.
By Lemma \ref{lemmafibres}(i) we may assume that the fiber over
$(x_n,y_n)$ is $(1- x_1/a_n) \tr (a_1, \dots, a_{n-1}) \times
\sq^{n-1}(1)$.  

Similarly, by mapping the discs $D(a_i)$ symplectomorphically to the
rectangles $\sq(a_i, 1)$ and then looking at the Lagrangian instead of
the symplectic splitting, we may think of $P(a_1, \dots, a_n)$ as $\sq
(a_1, \dots, a_n) \times \sq^n(1)$.

\subsubsection{Embeddings of polydiscs}
\label{embeddingsofpolydiscs}

We fold a polydisc $P(a_1, \dots, a_n)$ by folding a four dimensional
factor $P(a_i, a_j)$ for some $ i \neq j \in \{ 1, \dots , n \}$ and
leaving the other factor alone.
An already folded polydisc may be folded again by restricting the
folding process to a component containing no stairs.
The choice of $i$ and $j$ is only restricted by the condition that the
new image should still be embedded.

\paragraph{Embedding polydiscs into cubes}
\label{highembeddingpolydiscsintocubes}
In view of an application in subsection \ref{asymptoticpackings} we are
particularly interested in embedding thin polydiscs into cubes.
So fix $P^{2n}(a,\pi, \dots, \pi)$ and let $A$ be reasonably large. 
As explained above,
we think of $P^{2n}(a,\pi,\dots,\pi)$ as $\sq^n(a, \pi, \dots, \pi) \times
\sq^n(1)$ and of $C^{2n}(A)$ as $\sq^n(A) \times
\sq^n(1)$. The base direction will thus be the $z_1$-direction. Folding
into the $z_i$-direction for some $i \in \{2, \dots, n\}$, we will always lift into the $x_i$-direction.

We describe the process for $n=3$: First, fill a $z_1$-$z_2$-layer as
well as possible by lifting $N$ times into the $x_2$-direction (cf.\
\mbox{Figure \ref{figure20}}). Then lift once into the $x_3$-direction and fill
a second $z_1$-$z_2$-layer $\dots$. 
If $u_1$ is chosen appropriately, we will fold $N$ times into the
$x_3$-direction and fill $N+1$ $z_1$-$z_2$-layers.

The following proposition generalizes Proposition
\ref{propositionpintocube} to arbitrary
dimension.
\begin{proposition}
\label{propositionhigherintocubes}
Let $a>2\pi$ and $\ee >0$. Then $P^{2n}(\pi, \dots, \pi,a)$ embeds in
\mbox{$C^{2n}(s_{PC}^{2n}(a)+\ee)$}, where $s_{PC}^{2n}$ is given by
\[
s_{PC}^{2n}(a) = 
\left\{ \begin{array}{cl}
           (N+1)\pi,    & (N-1)N^{n-1} < \frac{a}{\pi}-2 \le (N-1)(N+1)^{n-1} \\
           \frac{a-2\pi}{(N+1)^{n-1}}+2\pi, & (N-1)(N+1)^{n-1} < \frac{a}{\pi}-2 \le N(N+1)^{n-1}\,
           .
        \end{array}
   \right. 
\]
\end{proposition}
\proof
The optimal embedding by folding $N$ times in each $z_1$-$z_2$-layer is described by 
\[
2 u_1 + ((N+1)^{n-1}-2)(u_1-\pi) =a,
\]
whence
\[
u_1 = \frac{a + ((N+1)^{n-1}-2)\pi}{(N+1)^{n-1}}.
\]
Thus
\[
A_N(a) = 
\max \left\{ \frac{a+2((N+1)^{n-1}-1)\pi}{(N+1)^{n-1}}, (N+1)\pi \right\},
\]
and the proposition follows.
\proofend

\subsubsection{Embeddings of ellipsoids}
\label{embeddingsofellipsoids}

We will concentrate on embedding ellipsoids $E^{2n}(\pi, \dots, \pi, a)$
with $a$ very large.

\paragraph{Embedding ellipsoids into cubes}
\label{highembeddingellipsoidsintocubes}

Studying embeddings 
\newline
\nobreak{$E^{2n}(\pi, \dots, \pi, a) \hookrightarrow
C^{2n}(A)$} of skinny ellipsoids into minimal cubes, we face the problem
of filling the fibers $\sq^{n-1}(A) \times \sq^{n-1}(1)$ of the cube by
many small fibers $\gg \tr^{n-1}(\pi) \times \sq^{n-1}(1)$ of the
ellipsoid. Forget about the irrelevant $y$-factors. Since $a$ is very
large, $\gg$ decreases very slowly. We are thus essentially 
left with the problem of filling $n-1$-cubes by equal $n-1$-simplices. 
This is trivial for $n-1 = 1$ and $n-1 = 2$, 
but impossible for $n-1 \ge 3$. Indeed, only
$2^{m-1}$ $m$-simplices $\tr^m(\pi)$ fit into $\sq^m(\pi)$, whence we
only get
\begin{equation}   \label{equationbadlimes}
\lim_{a \ra \infty} \frac{|E^{2n}(\pi, \dots, \pi,a)|}{
|C^{2n}(s_{EC}^{2n}(a))|} \ge \frac{2^{n-2}}{(n-1)!}.
\end{equation}

\begin{figure}[h] 
 \begin{center}
  \psfrag{x2}{$x_2$}
  \psfrag{x3}{$x_3$}
  \psfrag{A}{$A$}
  \psfrag{m}{$\mu_1$}
  \psfrag{n}{$\nu_1$}
  \psfrag{d}{$\dd_1$}
  \psfrag{dots}{$\vdots$}
  \leavevmode\epsfbox{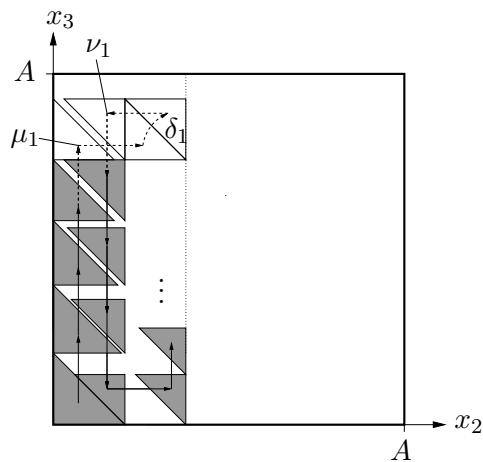}
 \end{center}
 \caption{Filling the cube fibres by the ellipsoid fibres} \label{figure31}
\end{figure}

We describe now the embedding process for $n-1=2$ in more detail 
\nobreak{(cf.\ \mbox{Figure \ref{figure31}})}. 
%
%
%
We first fill almost half of the ``first column'' of the cube fiber,
move the ellipsoid fibre out of this first column ($\mu_1$), deform it
to its complementary fiber ($\dd_1$), move this fiber back to the first
column ($\nu_1$), and fill almost all of the remaining room in the first
column. We then  pass to the second column and proceed as before.
The deformations $\dd_i$ are performed by applying 2-dimensional maps to
both symplectic directions of the ellipsoid fibers (see \mbox{Figure
\ref{figure52}} in \ref{highembeddingellipsoidsintoballs} and the text
belonging to it for more details). 
In order to guarantee that different stairs do not intersect, we arrange
the stairs arising from folding in such a way that the $z_1$-projections
of ``upward-stairs'' lie in $\{ 0 < y_1 < 1/2 \}$ while the
$z_1$-projections of ``downward-stairs'' lie in $\{1/2 < y_1 <1\}$, and
we arrange the stairs arising from moving in such a way that the
$z_1$-projections of the $\mu_i$- respectively $\nu_i$-stairs lie in 
$\{ 0 < y_1 < 1/4 \}$ respectively $\{ 1/4 < y_1 < 1/2 \}$ if $i$ is odd
and in $\{ 1/2 < y_1 < 3/4 \}$ respectively $\{ 3/4 < y_1 < 1 \}$ if $i$
is even (cf.\ \mbox{Figure \ref{figure6}}).
The $x_1$-intervals used for folding respectively moving will then be
double respectively four times as large as usual, but this will not
affect (\ref{equationbadlimes}).
\\
\\
{\bf Remark.}
We will prove in subsection \ref{asymptoticpackings} that the left hand
side of (\ref{equationbadlimes}) is 1 for any $n$.
\diam

\paragraph{Embedding ellipsoids into balls}
\label{highembeddingellipsoidsintoballs}

If we try to fill the fibers $\tr^{n-1}(A) \times \sq^{n-1}(1)$ of a
ball by many small fibers $\gg \tr^{n-1}(\pi) \times \sq^{n-1}(1)$ of a skinny
ellipsoid, we end up with a result for $s_{EB}^{2n}(a)$ as in
(\ref{equationbadlimes}). In the problem of embedding a skinny ellipsoid
into a minimal ball, however, both the fibers of the ellipsoid and the
fibers of the ball are balls. This may be used to prove
\begin{proposition}   
\label{propositionlim=1}
For any $n$,
\[
\lim_{a \ra \infty} \frac{|E^{2n}(\pi, \dots , \pi,a)|}{
|B^{2n}(s_{EB}^{2n}(a))|} = 1.
\]
\end{proposition}

\proof
\begin{figure}[h] 
 \begin{center}
  \psfrag{x1}{$x_1$}
  \psfrag{A}{$A$}
  \psfrag{A/l}{$A/l$}
  \psfrag{2A/l}{$2A/l$}
  \psfrag{A-A/l}{$A-A/l$}
  \psfrag{T1}{$T_1$}
  \psfrag{T2}{$T_2$}
  \psfrag{Tl}{$T_l$}
  \leavevmode\epsfbox{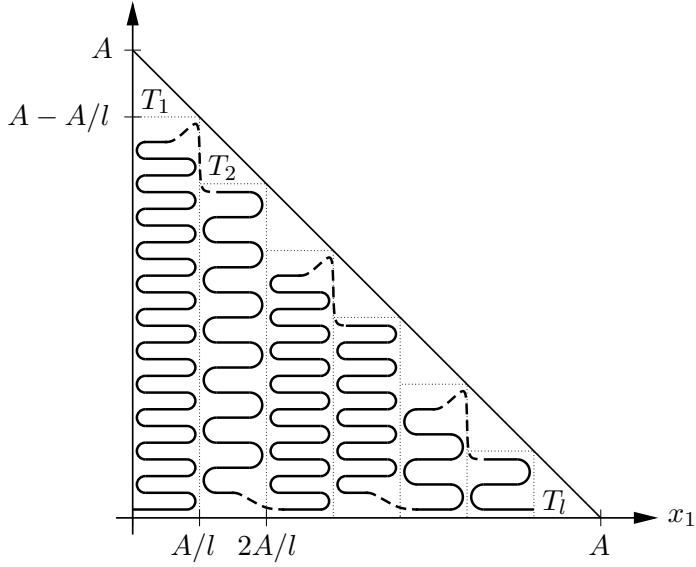}
 \end{center}
 \caption{Embedding a skinny ellipsoid into a ball} \label{figure50}
\end{figure}
%
The idea of the proof is very simple: 
Instead of packing a large simplex by small simplices, we will leave the
simplices alone and pack the cubes by small cubes, a trivial problem.

\begin{figure}[h] 
 \begin{center}
  \psfrag{y2}{$y_2$}
  \psfrag{y3}{$y_3$}
  \psfrag{1}{$1$}
  \psfrag{C}{$C$}
  \psfrag{C'}{$C'$}
  \psfrag{C''}{$C''$}
  \psfrag{1/k1}{$1/k_1$}
  \psfrag{1/k2}{$1/k_2$}
  \psfrag{N1/k1}{$N_1/k_1$}
  \psfrag{max}{$\max (1/k_1,1/k_2)$}
  \leavevmode\epsfbox{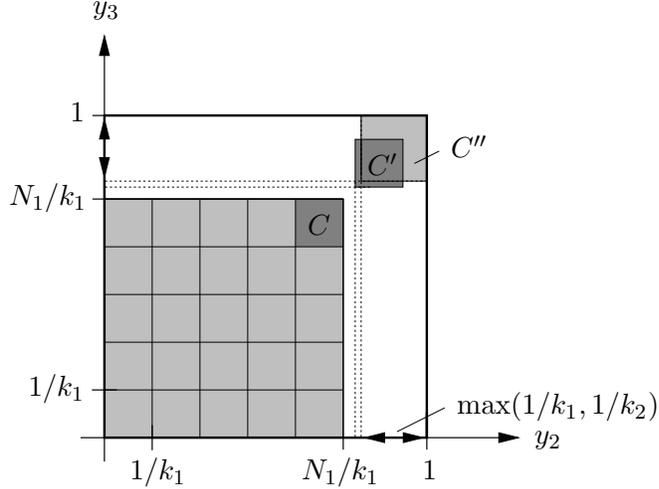}
 \end{center}
 \caption{Filling the $y$-factor of the fibers} \label{figure51}
\end{figure}

\begin{figure}[h] 
 \begin{center}
  \psfrag{xi}{$x_i$}
  \psfrag{yi}{$y_i$}
  \psfrag{1}{$1$}
  \psfrag{C}{$C$}
  \psfrag{C'}{$C'$}
  \psfrag{C''}{$C''$}
  \psfrag{k1pi}{$(1-\frac{b_1}{a}) k_1 \pi$}
  \psfrag{k2pi}{$(1-\frac{b_1}{a}) k_2 \pi$}
  \psfrag{N1/k1}{$N_1/k_1$}
  \psfrag{1/k1}{$1/k_1$}
  \psfrag{1/k2}{$1/k_2$}
  \psfrag{aa}{$\aa_i$}
  \psfrag{ee}{$\ee$}
  \leavevmode\epsfbox{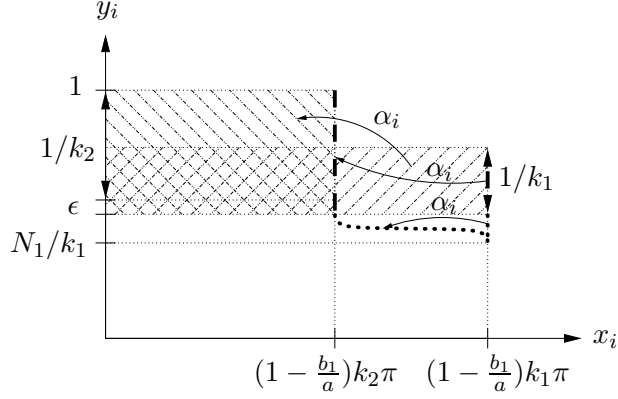}
 \end{center}
 \caption{Rescaling the fibers} \label{figure52}
\end{figure}

So pick a very large $l \in \NN$, write
\[
P_i = B^{2n}(A) \cap \left\{ \frac{(i-1)A}{l} < x_1 < \frac{iA}{l} \right\}, 
\qquad 1 \le i \le l,
\]
and set
\[
k_1 = \frac{A - A/l}{\pi},
\]
where $A$ is again a parameter which will be fixed later on.
After applying the diagonal map $\mbox{diag} \, [ k_1, \dots, k_1, 1/k_1, \dots,
1/k_1]$ to the fibers, the ellipsoid is contained in $\sq (a,1) \times
\tr^{n-1}(k_1 \pi) \times \sq^{n-1}(1/k_1).$
We will embed some part $\sq (b_1, 1) \times \tr^{n-1}(k_1 \pi) \times
\sq^{n-1}(1/k_1)$ of this set into $P_1$ by
fixing the simplices and moving the cubes along the $y_i$-directions $(2
\le i \le n)$ (see \mbox{Figure \ref{figure50}} and \mbox{Figure \ref{figure51}}).

%

We want to fill as much of $\sq^{n-1}(1)$ by cubes $\sq^{n-1}(1/k_1)$ as
possible. However, in order to use also the space in $P_2$ optimally, we
will have to deform the ellipsoid fibers before passing to $P_2$, and
for this we will have to use some space in $\sq^{n-1}(1)$.
Assume that we fold $N_1'$ times in each $z_1$-$z_2$-layer and by this
embed $\sq(b_1',1) \times \tr^{n-1}(k_1 \pi) \times \sq^{n-1}(1/k_1)$
into $P_1$. The maximal ellipsoid fiber over $P_2$ will then be 
\[
\left(1- \frac{b_1'}{a} \right) \tr^{n-1}(k_1 \pi) 
\times \sq^{n-1} \left( \frac{1}{k_1} \right).
\] 
We want to deform this fiber to a fiber 
\[
\left(1- \frac{b_1'}{a}\right) \tr^{n-1}(k_2' \pi) \times
\sq^{n-1}\left(\frac{1}{k_2'}\right)
\]
fitting into the minimal ball fiber $\tr^{n-1}(A-2A/l) \times
\sq^{n-1}(1)$ over $P_2$. We thus define $k_2'$ by 
$(1- b_1'/a) k_2' \pi = A - 2A/l$.
As we shall see below, the appropriate ellipsoid fiber deformation can
then be achieved in $\sq^{n-1}(1) \setminus \sq^{n-1}(1- \max (1/k_1,
1/k_2))$.

The optimal choice of $N_1'$ and $k_2'$ is the solution of the system
\begin{equation*} 
\left. \begin{array}{rcl}
            N_1 & = & \max \left\{ N \in \NN \, \big| \, N \mbox{ even},\,
            \frac{N}{k_1} < 1 - \max \big(\frac{1}{k_1},\frac{1}{k_2}\big) \right\}    \\     
            k_2 \pi & = & \Big(A- \frac{2A}{l}\Big)  \Big/
            \left(1 - \frac{b_1(N)}{a}\right) 
        \end{array}
   \right\}.
\end{equation*}
By folding $N_1$ times in each $z_1$-$z_2$-layer we fill nearly all
of $\sq^{n-1}(1- \max (1/k_1, 1/k_2))$ and indeed stay away from $\sq^{n-1}(1) \setminus
\sq^{n-1}(1- \max (1/k_1, 1/k_2))$ (cf.\ \mbox{Figure \ref{figure51}}).

The deformation of the ellipsoid fibres is achieved as follows:
We first move the cube $C$ along all $y_i$-directions, 
$i \ge 2$, by $1- \max (1/k_1, 1/k_2) - (N_1-1)/k_1 -\ee$ 
for some $\ee \in \; ]0, 1-
\max (1/k_1, 1/k_2) - N_1/k_1[$. This can be done whenever $A/l > n
\pi$. We then deform the translate $C'$ to $C''$.
%
%
This deformation is the restriction to $(1-b_1/a) \tr^{n-1}(k_1 \pi) \times
\sq^{n-1}(1)$ of a product of $n-1$ two-dimensional symplectic maps $\aa_i$ which are
explained in \mbox{Figure \ref{figure52}}: On $y_i \le N_1/k_1$, $\aa_i$ is the
identity, and on $y_i \ge 1-1/k_2 -\ee$ it is an affine map with linear
part 
\[
(x_i, y_i) \mapsto \left(\frac{k_2}{k_1} x_i,\frac{k_1}{k_2} y_i\right).
\] 

Assume that we can choose $A$ such that proceeding in this way, we
successively fill a large part of all the $P_i, \; 1 \le i \le l-1$,
and leave $P_l$ untouched, i.e.\ the embedding process ends exactly when
passing from $P_{l-1}$ to $P_l$ (cf.\ \mbox{Figure \ref{figure50}}). The
process is then described by the equations for the pairs $(N_i,
k_{i+1}), \; 1 \le i \le l-2$, 
\begin{equation}  \label{equationmany}
\left.  \begin{array}{rcl} 
N_i         & = &   \max \left\{ N \in \NN \, \big| \, N \mbox{ even},\,
                    \frac{N}{k_i} 
                    < 1 - \max \big(\frac{1}{k_i}, \frac{1}{k_{i+1}}\big) \right\}      \\ 
k_{i+1} \pi & = &   \Big(A - \frac{(i+1)A}{l}\Big) \Big/ 
              \left(1 - \frac{\sum_{j=1}^{i-1} b_j(N_j) + b_i(N)}{a}\right)  
        \end{array}   
\right\}, 
\end{equation} 
where $b_j(N_j)$ is the $x_1$-length of the part embedded into $P_j$,
and by
\begin{equation*}
N_{l-1} = \max \{ n \in \NN \, | \, N \mbox{ even}, \, N < k_{l-1} \}.
\end{equation*}
We finally observe that, in reality, the system (\ref{equationmany})
splits. Indeed, the second line in (\ref{equationmany}) readily implies
that $k_i < 2 k_{i+1}$ whenever $i \le l-2$. Thus, the first line in 
(\ref{equationmany}) reads
$N_i = \max \{ N \in \NN \, | \, N \mbox{ even},\, N/k_i < 1 - 1/k_i \}$,
and the embedding process is described by 
\addtocounter{equation}{1}
\begin{align*}   
N_i         & \;\;=\;\;  \max \{ N \in 2 \NN \, | \, N < k_i-1 \}     
                      \tag{\arabic{equation}.1}  \label{equationsystem1}                       \\ 
k_{i+1} \pi & \;\;=\;\;  \left(A - \frac{(i+1)A}{l}\right) \Bigg/ 
                          \left(1 - \frac{\sum_{j=1}^i b_j(N_j)}{a}\right)
                      \tag{\arabic{equation}.2}  \label{equationsystem2}                        \\
N_{l-1}     & \;\;=\;\;  \max \{ N \in 2 \NN \, | \, N < k_{l-1} \}.
                      \tag{\arabic{equation}.3}  \label{equationsystem3} 
\end{align*}
We now argue that such an $A$ indeed exists, and that it is the minimal
$A$ for which the above embedding process succeeds.

Observe first that such a minimal $A$, which we denote by $A_0$, indeed
exists, for clearly, if $A$ was chosen very large, the embedding process
will end at some $P_i$ with $i<l-1$, and if $A$ was chosen very small,
it won't succeed at all.

Suppose now that the embedding process for $A_0$ ends before passing
from $P_{l-1}$ to $P_l$. Pick $A'<A_0$ and write $k_i$ and $N_i$
respectively $k_i'$ and $N_i'$ for the embedding parameters belonging
to $A_0$ respectively $A'$. If $A_0-A'$ is small, $k_1-k_1'$ is small
too; 
thus, by (\ref{equationsystem1}), $N_1 = N_1'$ whenever $A_0-A'$ is
small enough.
But then, $b_1(N_1) - b_1'(N_1)$ is small, whence
(\ref{equationsystem2}) shows that $k_2 - k_2'$ is small.
Arguing by induction, we assume that $N_j = N_j'$ and that
$b_j(N_j)-b_j'(N_j)$ and $k_{j+1}-k_{j+1}'$ are small for $j \le
i$. Then, by (\ref{equationsystem1}) or (\ref{equationsystem3}), and
after choosing $A_0-A'$ even smaller if necessary, we may assume that 
$N_{i+1} = N_{i+1}'$. If $i+2 \le l-1$, $b_{j+1}(N_{j+1}) -
b_{j+1}'(N_{j+1})$ 
is then small too, whence (\ref{equationsystem2}) shows that $k_{i+2} -
k_{i+2}'$ is small.

We hence may assume that all differences $b_i - b_i'$ are arbitrarily
small. But then the embedding process for $A'$ will succeed as well, a
contradiction.  

\medskip
Recall that $A_0 = A_0(a,l)$ still depends on $l$. The best embedding
result provided by the above procedure is thus 
\[
s_{EB}^{2n}(a) = \min_{l \in \NN} \{ A_0 (a,l) \}.
\]
Set
\[
q(a,l) = 1 - \frac{|E^{2n}(\pi, \dots, \pi, a)|}{|B^{2n}(A_0(a,l))|}
\]
and 
\[
q(a) = 1 - \frac{|E^{2n}(\pi, \dots, \pi, a)|}{|B^{2n}(s_{EB}^{2n}(a))|}.
\]
In order to prove the proposition, we have to show that 
\begin{equation}   \label{equationratio}   
\lim_{a \ra \infty} q(a) = 0.
\end{equation}

Given any $a$ and $l$, the region in $B^{2n}(A_0(a,l))$ which is not
covered by the image of $E^{2n}(\pi, \dots, \pi, a)$ is the disjoint
union of four types of regions $R_h(a,l)$, $1 \le h \le 4$. 
\begin{itemize}
\item[]
$R_1(a,l)$ is the union of the ``triangles'' $T_i(a,l)$ (see \mbox{Figure
\ref{figure50}}).
\item[]
$R_2(a,l)$ is the space needed for folding (see \mbox{Figure \ref{figure56}}).
\item[]
$R_3(a,l)$ is the union of the space needed to deform the ellipsoid
fibers and the space caused by the fact that the $N_i$
have to be integers (see \mbox{Figure \ref{figure51}}).
\item[]
$R_4(a,l)$ is the image of the difference set of the embedded set
and $E^{2n}(\pi, \dots, \pi, a)$ (see \mbox{Figure \ref{figure57}}).
\end{itemize}
Detailed descriptions of these sets are given below.

Let $\ee >0$ be small. We will find $a_\ee$ and $l_\ee$ such that             
\refstepcounter{equation}
\label{e:ee}
\makeatletter
\protected@edef\nummer{\theequation}
\makeatother
\begin{equation}
  \frac{|R_h(a, l_\ee)|}{|B^{2n}(A_0(a, l_\ee))|} < \ee \quad \mbox{ for all } a
  \ge a_\ee,                                            \tag{\nummer.h}
\end{equation}
$1 \le h \le 4$.
Since the sets $R_h(a,l)$ are disjoint and $q(a) \le q(a,l)$,
(\ref{e:ee}.h), $1 \le h \le 4$, imply (\ref{equationratio}).

Set $R_{h,i}(a,l) = R_h(a,l) \cap P_i(a,l)$.
We first of all observe that the ratio 
$
|R_1(a, l)| / |B^{2n}(A_0(a, l))|
$
depends only on $l$ and can be made arbitrarily small by taking $l$
large.
We thus find $l_1$ such that
\begin{equation*}     
\frac{|R_1(a, l)|}{|B^{2n}(A_0(a, l))|} < \ee \quad \mbox{ for all $a$ and
$l \ge l_1$}.
\end{equation*}
Moreover, notice that given $\zz >0$ we can choose $l_1$ such that for all
$a$ and $l \ge l_1$
\begin{equation}  \label{e:b}   
\frac{|R_{1,i}(a,l)|}{|P_i(a,l)|} < \zz \quad \mbox{ whenever $i$ is not
too near to $l-1$.}
\end{equation}
Here and in the sequel, ``$i$ too near to $l-1$'' stands for
``$1-i/(l-1)$ smaller  than a constant which can be made arbitrarily small by
taking first $l$ and then also $a$ large''.

Next, our construction clearly shows that given $\zz$ as above and 
$l$ being fixed we may find $a_1$ such that for $a \ge a_1$ and for all
$i \in \{ 1, \dots , l-1 \}$
\begin{equation}   \label{equationtworatios}   
\frac{|R_{2,i}(a,l)|}{|P_i(a,l)|} < \zz  \quad \mbox{ and  } \quad
\frac{|R_{3,i}(a,l)|}{|P_i(a,l)|} < \zz.
\end{equation}
In particular, given any $l_\ee \ge l_1$, we find $a_\ee$ 
such that
(\ref{e:ee}.1), (\ref{e:ee}.2) and (\ref{e:ee}.3) hold true.
\begin{figure}[h] 
 \begin{center}
  \psfrag{pi}{$\pi$}
  \psfrag{a}{$a$}
  \leavevmode\epsfbox{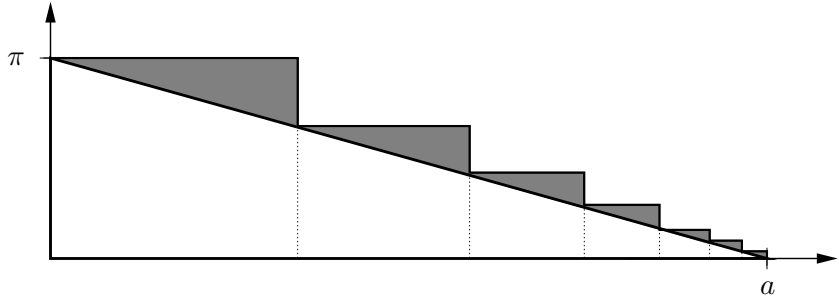}
 \end{center}
 \caption{$Y(a,8) \subset X(a,8)$} \label{figure57}
\end{figure}
%
%
%

Recall that the embedding $\ff_{a,l} \colon E^{2n}(\pi, \dots, \pi,a)
\hookrightarrow B^{2n}(A_0(a,l))$ is defined on a larger domain with
piecewise constant fibres. Set 
\begin{align*}
X_i(a,l) \;&= \;\ff_{a,l}^{-1}(P_i(a,l)),                       \\
Y_i(a,l) \;&= \;X_i(a,l) \setminus E^{2n}(\pi, \dots, \pi, a),  \\
Z_i(a,l) \;&= \;X_i(a,l) \cap E^{2n}(\pi, \dots, \pi, a)
\end{align*}
and
$X(a,l) = \coprod_{i=1}^{l-1} X_i(a,l)$,
$Y(a,l) = \coprod_{i=1}^{l-1} Y_i(a,l)$,
$Z(a,l) = \coprod_{i=1}^{l-1} Z_i(a,l)$
(cf.\ \mbox{Figure \ref{figure57}}),
and recall that we denoted the $u$-width of $X_i(a,l)$ by
$b_i(a,l)$. Assume now that $\zz$ is small. Then (\ref{e:b}) and
(\ref{equationtworatios}) show that for $a \ge a_\ee$ and $i$ not too
near to $l_\ee-1$, $|X_i(a,l_\ee)|/|P_i(a,l_\ee)|$ is near to $1$. Thus,
a simple volume comparison shows that if $l_\ee$ is large,
$b_i(a,l_\ee)/a$ and hence also $|R_{4,i}(a,l_\ee)|/|P_i(a,l_\ee)| =
|Y_i(a,l_\ee)|/|P_i(a,l_\ee)|$ is
small for these $a$ and $i$. In particular, we may choose $l_\ee$ and
$a_\ee$ such that (\ref{e:ee}.4) holds true too.

This completes the 
proof of Proposition \ref{propositionlim=1}.
For later purposes, we state that given $\zz >0$, we may find $l_0$ and
$a_0$ such that for all $a \ge a_0$ and $i$ not too near to $l_0-1$
\begin{equation}   \label{e:c}   
\frac{|R_{h,i}(a,l_0)|}{|P_i(a,l_0)|} < \zz, \qquad
1 \le h \le 4.
\end{equation}
\proofend

\medskip
The above proof gives no information about the convergence speed in
(\ref{equationratio}). The remainder of this paragraph is devoted to the
proof of 
\begin{proposition}  \label{p:c}
Given $\ee >0$ there is a constant $C(n,\ee)$ such that for all $a$
\[
1 - \frac{|E^{2n}(\pi, \dots, \pi,a)|}{|B^{2n}(s_{EB}^{2n}(a))|}
                                 < C(n, \ee) a^{-\frac{1}{2n} + \ee}.
\]
\end{proposition}
\proof
The proposition follows from the existence of a pair $(a_0, l_0)$
such that for $a \in I_k(a_0) =  [4^{kn} a_0, 4^{(k+1)n} a_0[$, $k \in \NN_0$,
\begin{equation}   \label{e:qq}  
(2-\ee) q(4^na, 2^{k+1}l_0) < q(a,2^k l_0).
\end{equation}
Indeed, choose $C(n, \ee)$ so large that $C(n, \ee) a^{-\frac{1}{2n} + \ee} 
> q(a)$ for $a<a_0$ and
\begin{equation}   \label{e:199}  
C(n,\ee) a^{-\frac{1}{2n}} > q(a, l_0) 
\quad \mbox{ for } a \in I_0(a_0).
\end{equation}
Then, if $a \in I_k(a_0) \; \mbox{ for some } k \in
\NN$, 
\begin{eqnarray*}
q(a) \; \le \; q(a, 2^k l_0) & \overset{(\ref{e:qq})}{<} &
(2-\ee)^{-k} q \left( \frac{a}{4^{kn}}, l_0 \right)           \\
                             & \overset{(\ref{e:199})}{<}    &
(2-\ee)^{-k} C(n,\ee) 2^k a^{-\ee} a^{-\frac{1}{2n} + \ee}    \\
                             & \le                           & 
(2-\ee)^{-k} C(n,\ee) 2^k 4^{-\ee kn} 
a_0^{-\ee} a^{-\frac{1}{2n} + \ee}                            \\
                             & <                             & 
(2-\ee)^{-k} 2^k 4^{-\ee kn} C(n,\ee)  
a^{-\frac{1}{2n} + \ee}                                       \\
                             & <                             & 
C(n,\ee) a^{-\frac{1}{2n} + \ee}.
\end{eqnarray*}

\smallskip
So let's prove (\ref{e:qq}). 
Fix $(a_0, l_0)$ and $\hat{a} \in I_0(a_0)$ and set
$a_k = 4^{kn} a_0$, $\hat{a}_k = 4^{kn} \hat{a}$, $l_k = 2^k l_0$ and
\[
\rho_k = \frac{A_0(\hat{a}_{k+1}, l_{k+1})}{A_0(\hat{a}_k, l_k)},
\]
$k \in \NN_0$.
Given a specified subset $S(a,l)$ of $B^{2n}(A_0(a,l))$ and a parameter
$p(a,l)$ belonging to the embedding $\ff_{a,l} \colon 
E^{2n}(\pi, \dots, \pi, a)
\hookrightarrow B^{2n}(A_0(a,l))$, we write
$_kS$ and $_kp$ instead of $S(\hat{a}_k,l_k)$ and
$p(\hat{a}_k,l_k)$. Moreover, we write $_kS'$ for the rescaled subset
$\frac{1}{\rho_k} S(\hat{a}_{k+1}, l_{k+1})$ of $\frac{1}{\rho_k}
B^{2n}(A_0(\hat{a}_{k+1}, l_{k+1}))$ and
$_kp'$ for the parameter belonging to the rescaled embedding
$\frac{1}{\rho_k} E^{2n}(\pi, \dots, \pi, \hat{a}_{k+1}) \hookrightarrow
\frac{1}{\rho_k} B^{2n}(A_0(\hat{a}_{k+1}, l_{k+1}))$.
Finally, write $\rho$, $S$, $S'$, $p$, $p'$ instead
of $\rho_0$, $_0S$, $_0S'$, $_0p$, $_0p'$,
set $E = E^{2n}(\pi, \dots, \pi, \hat{a})$, 
$E' = \frac{1}{\rho} E^{2n}(\pi, \dots, \pi, \hat{a}_1)$ and
$B = B^{2n}(A_0(\hat{a}, l_0)),$
and observe that $B = B'$.

We claim that we can find $(a_0, l_0)$ such that for all $k \in \NN_0$, 
$\hat{a}_k \in I_k(a_0)$ and $i$ not too near to $l_k-1$
\refstepcounter{equation}
\label{e:200} 
\makeatletter
\protected@edef\nummer{\theequation}
\makeatother
\begin{equation}   
(4-\ee) |_kR_{h,2i(-1)}'| < |_kR_{h,i}|,       \tag{\nummer.h.k}
\end{equation}
$1 \le h \le 4$.
We will first prove (\ref{e:200}.h.0)
and will then check that the
conditions valid for $(\hat{a}, l_0)$ which allowed us to conclude 
(\ref{e:200}.h.0) are
also valid for $(\hat{a}_k, l_k)$ provided that (\ref{e:200}.h.m) holds true for $m \le
k-1$. Arguing by induction, we thus see that (\ref{e:200}.h.k) holds true
for all  $k \in \NN_0$.

Set $\ee_1 = \ee/16$ and observe that for all $k \in \NN_0$ and $i$ not
too near to $l_k-1$
\begin{equation}   \label{e:pp}  
|_kP_{2i-1}'| > |_kP_{2i}'| > \left( \frac{1}{2} - \ee_1 \right)|_kP_i|.
\end{equation}
We conclude that for $k \in \NN_0$, $\hat{a}_k \in I_k(a_0)$ and $i$ not
too near to $l_k-1$
\refstepcounter{equation}
\label{e:201} 
\makeatletter
\protected@edef\nummer{\theequation}
\makeatother
\begin{equation}   
\left( 2- \frac{3 \ee}{4} \right) 
\frac{|R_{h,2i(-1)} (\hat{a}_{k+1}, l_{k+1})|}{|P_{2i(-1)}
(\hat{a}_{k+1}, l_{k+1})|}
<
\frac{|R_{h,i} (\hat{a}_k, l_k)|}{|P_i(\hat{a}_k, l_k)|}, 
                               \tag{\nummer.h.k}
\end{equation}
$1 \le h \le 4$. 
In particular, there is $(a_0,l_0)$ such that for all $\hat{a} \in
I_0(a_0)$, 
\[
(2-\ee) 
\frac{|R_h (\hat{a}_{k+1}, l_{k+1})|}{|B^{2n}(A_0(\hat{a}_{k+1}, l_{k+1}))|}
<
\frac{|R_h (\hat{a}_k, l_k)|}{|B^{2n}(A_0(\hat{a}_k, l_k))|},
\] 
$1 \le h \le 4$. 
Since $R_h(a,l)$ are disjoint, this implies
(\ref{e:qq}).

\bigskip
{\bf{(R1)}}
Let $R_1(a,l) = \coprod_{i=1}^l T_i(a,l)$ 
be the union of the ``triangles'' $T_i(a,l) \subset B^{2n}(A_0(a,l))$
(see \mbox{Figure \ref{figure50}}).
\begin{figure}[h] 
 \begin{center}
  \psfrag{T1}{$T_1$}
  \psfrag{T1'}{$T_1'$}
  \psfrag{T2'}{$T_2'$}
  \psfrag{Tl}{$T_{l_0}$}
  \psfrag{Tl1'}{$T_{2 l_0-1}'$}
  \psfrag{Tl'}{$T_{2 l_0}'$}
  \leavevmode\epsfbox{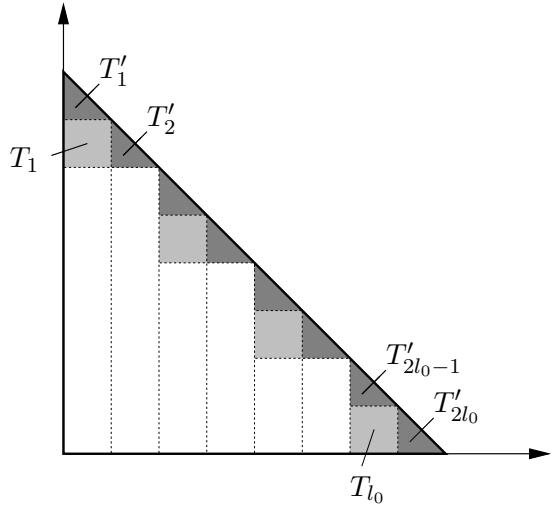}
 \end{center}
 \caption{$R_1$ and $R_1'$} \label{figure53}
\end{figure}
%
%
%
$R_{1,2i(-1)}'$ is a subset of $R_{1,i}$, and $|R_{1,i}|/|R_{1,2i(-1)}'| =
|T_i|/|T_{2i(-1)}'|$ depends only on
$l_0$ (see \mbox{Figure \ref{figure53}}). Clearly, 
$4 - |T_i|/|T_{2i(-1)}'|$ 
is small if $|T_i| / |P_i|$ is small enough.
By taking $l_0$ large, we may make $|T_i| / |P_i|$ arbitrarily small
for $i$ not too near to $l_0 -1$.
Thus, (\ref{e:200}.1.0) holds true whenever $l_0$ is large
enough. Observe finally that (\ref{e:200}.1.0) 
implies (\ref{e:200}.1.k), $k \in \NN$.

\bigskip
{\bf{(R2)}}
Recall that the $x_1$-length of the space needed for folding equals the
fiber capacity at the place where we fold. 
The staircases needed for folding are thus contained in $R_2(a,l) =
\coprod_{i=1}^{l-1} R_{2,i}(a,l)$, where $R_{2,i}(a,l)$ equals
\[
Q_i(a,l) \setminus 
\Bigg\{ \frac{(i-1)A}{l} + \pi \Bigg(1-\frac{\sum_{j=1}^{i-1} b_j}{a}\Bigg)
     < x_1 < 
\frac{iA}{l} - \pi \Bigg(1- \frac{\sum_{j=1}^{i-1} b_j}{a}\Bigg) \Bigg\}.
\]
Here, we put
\[
Q_i(a,l) = P_i(a,l) \setminus T_i(a,l).
\]
\begin{figure}[h] 
 \begin{center}
  \psfrag{x1}{$x_1$}
  \psfrag{A/2l}{$\frac{A_0}{2l_0}$}
  \psfrag{A/l}{$\frac{A_0}{l_0}$}
  \leavevmode\epsfbox{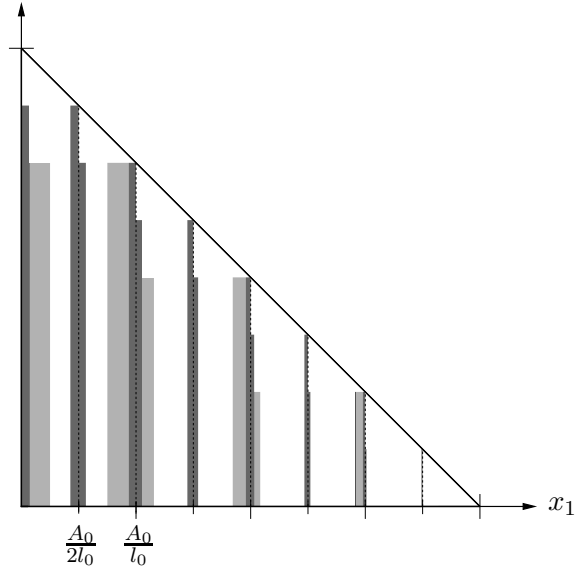}
 \end{center}
 \caption{$R_2$ and $R_2'$} \label{figure56}
\end{figure}
%
%

Observe that for $i$ not too near to $l_k-1$, $|_kQ_{2i-1}' \cap {_kT_i}|
/ | _kQ_{2i-1}'| \ra 0$ as $l_k \ra \infty$ (cf.\ \mbox{Figure \ref{figure53}}).
Hence, also $| _kR_{2,2i-1}' \cap {_kT_i}| / | _kR_{2,2i-1}'| \ra 0$ as
$l_k \ra \infty$. We may thus neglect $_kR_{2,2i-1}' \cap {_kT_i}$ and
prove
(\ref{e:200}.2.k)
with $_kR_{2,2i-1}'$ replaced by $_kR_{2,2i-1}'
\setminus {_kT_i}$ (which we denote again by $_kR_{2,2i-1}'$).

\smallskip
If $u_i = \sum_{j=1}^{i-1} b_j$ respectively $u_i' = \sum_{j=1}^{i-1}
b_j'$ is the $x_1$-coordinate at which the image of $E$
respectively $E'$ enters $P_i$, then the volume embedded into 
$
\coprod_{j=1}^{i-1} P_j
$
is
\begin{equation}   \label{equationvolumes}   
\frac{\pi^{n-1}}{\hat{a}^{n-1}n!} \, 
\big[ \hat{a}^n - (\hat{a}-u_i)^n \big]
\quad \mbox{ resp.\ } \quad
\frac{\pi^{n-1}}{\hat{a}_1^{n-1}n!} \Bigg[
\bigg(\frac{\hat{a}_1}{\rho} \bigg)^n -
\bigg(\frac{\hat{a}_1}{\rho}
-u_i' \bigg)^n \Bigg],
\end{equation}
and the fiber capacity at $u_i$ respectively $u_i'$ is
\begin{equation}   \label{equationcapacities}   
c_i = \frac{\pi}{\hat{a}} (\hat{a} -u_i) \quad \mbox{ resp.\ } \quad
c_i' = \frac{\pi}{\hat{a}_1} \bigg(\frac{\hat{a}_1}{\rho} -u_i'\bigg).
\end{equation}
Thus, $c_1 = \rho c_1'$.
We claim that
\begin{equation}   \label{equationclaim}  
\begin{array}{l} 
c_i > (1-\ee_1) \rho c_{2i(-1)}' 
\\
\;\;\;\;\;\;\;\;\;\;\;\;\;
\mbox{ whenever $\hat{a}$ is large enough and $i$ is not too near to
$l_0-1$}.
\end{array}
\end{equation}
Since $c_{2i-1}' > c_{2i}'$, it suffices to show that 
\begin{equation*}  
c_i > (1-\ee_1) \rho c_{2i-1}' 
\quad \mbox{ for $\hat{a}$ large enough and $i$ not too near to $l_0-1$.}
\tag{\ref{equationclaim}'}
\end{equation*}
So assume that there is an $i$ violating the inequality in
(\ref{equationclaim}') and set
\[
i_0 = \min \{ 1 \le i \le l_0-1 \, | \, c_i \le (1-\ee_1)\rho c_{2i-1}' \}.
\]
Let $\zz >0$ be so small that 
\begin{equation}   \label{e:x}   
\zz < \ee_1
\end{equation}
and set
\[
z_i(a,l) = \frac{|Z_i(a,l)|}{|P_i(a,l)|}  
\quad \mbox{ and } \quad 
z(a,l) = \frac{|Z(a,l)|}{|B^{2n}(A_0(a,l))|}.  
\]
By the definition of $\rho$, $z$ and $z'$,
\begin{equation}   \label{equationrhon}   
\rho^n = 4^n \frac{z}{z'}.
\end{equation}
By (\ref{e:c}), for any large enough $l_0$ there is $a_0$ such that for
all $\hat{a} \in I_0(a_0)$ and $i$ not too near to $l_0-1$
\begin{equation}   \label{e:zi}   
z_i > 1- \zz.
\end{equation}
We have seen in ($R_1$) that for all $i \in \{ 1, \dots , l_0 \}$
\begin{equation}   \label{e:1}   
|R_{1,2i(-1)}'| < |R_{1,i}|.
\end{equation}
Moreover, if $\zz$ is small enough, we clearly have that for $i$ not too
near to $l_0-1$
\begin{equation}   \label{e:ss}   
c_i > c_{2i(-1)}'.
\end{equation}
This implies that for these $i$
\begin{equation}   \label{e:2}   
|R_{2,2i(-1)}'| < |R_{2,i}|.
\end{equation}
We now assume that $a_0$ is so large compared to $l_0$ that
\begin{equation}   \label{e:Gg}
A_0(a_0,l_0) > 12 l_0 \pi.
\end{equation}
Then, $A_0(\hat{a},l_0) > 12 l_0 \pi > 12 l_0 c_i$, i.e.\
\begin{equation}   \label{e:Hh}
\frac{A_0(\hat{a},l_0)}{l_0} > 12 c_i, \qquad 1\le i \le l_0-1.
\end{equation}
\begin{equation}   \label{e:3}   
|R_{3,2i(-1)}'| < |R_{3,i}|
\end{equation}
now follows  from (\ref{e:ss}) in the same way as (\ref{e:zz}) will
follow from (\ref{equationclaim}).
Finally, for $\zz$ small enough and $i$ not too near to $l_0-1$ we clearly have that
\begin{equation}   \label{e:4}   
|R_{4,2i(-1)}'| < |R_{4,i}|.
\end{equation}
We conclude from (\ref{e:1}), (\ref{e:2}), (\ref{e:3}) and (\ref{e:4})
and (\ref{e:pp}) that
\begin{equation*}   
\frac{|R_{h,2i(-1)}'|}{|P_{2i(-1)}'|} < 3 \frac{|R_{h,i}|}{|P_i|},     
\qquad 1 \le h \le 4.
\end{equation*}
This shows that
\begin{equation}   \label{e:zi'}   
z_i' > 1- 3 \zz.
\end{equation}
Set
\begin{align}   
z_{<i} =
\frac{\big| \coprod_{j=1}^{i-1}Z_j \big|}{\big| \coprod_{j=1}^{i-1}P_j\big|}
\quad &\mbox{ and } \quad
z_{<i}' = 
\frac{\big| \coprod_{j=1}^{i-1}Z_j'\big|}{\big| \coprod_{j=1}^{i-1}P_j'\big|}.
\intertext{By (\ref{e:zi}) and (\ref{e:zi'}), we may assume that for all $i \in
\{1, \dots, l_0-1 \}$}  
z_{<i} > 1- \zz \quad &\mbox{ and } \quad z_{<i}' > 1- 3 \zz. \notag \\
\intertext{In particular,}      
z > 1- \zz \quad&\mbox{ and } \quad z' > 1- 3 \zz   \label{e:z}     \\
\intertext{and}       
z_{<i_0} > 1- \zz \quad &\mbox{ and } \quad z_{<i_0}' > 1- 3 \zz.
\label{e:zi0}
\end{align}
Comparing the two volumes embedded into $\coprod_{j=1}^{i_0-1}P_j$, we get from
(\ref{equationvolumes}) that
\begin{equation}   \label{equationcompare}   
z_{<i_0}' \frac{\pi^{n-1}}{\hat{a}^{n-1}n!}\big[\hat{a}^n - (\hat{a} - u_{i_0})^n\big] 
=
z_{<i_0} \frac{\pi^{n-1}}{\hat{a}_1^{n-1}n!} \Bigg[\bigg(\frac{\hat{a}_1}{\rho}\bigg)^n - 
\bigg(\frac{\hat{a}_1}{\rho} - u_{2i_0-1}'\bigg)^n \Bigg].
\end{equation}
By (\ref{equationcapacities}), $c_{i_0} \le (1-\ee_1) \rho c_{2i_0-1}'$
translates to 
\begin{equation}   \label{equationtranslate}   
u_{2 i_0 -1}' \le \frac{4^n}{(1-\ee_1)\rho} (u_{i_0} - \ee_1 \hat{a}).
\end{equation}
Plugging (\ref{equationtranslate}) into (\ref{equationcompare}), we find 
\[
\Bigg(z_{<i_0}\bigg(\frac{4}{\rho}\bigg)^n - z_{<i_0}'\Bigg) \hat{a}^n
\ge  
\Bigg(z_{<i_0} \bigg(\frac{4}{\rho
(1-\ee_1)}\bigg)^n - z_{<i_0}'\Bigg) (\hat{a} - u_{i_0})^n,
\]
and using (\ref{equationrhon}) and dividing by $z_{<i_0}$ we get
\begin{equation}   \label{equationlast}   
\bigg(\frac{z'}{z} - \frac{z_{<i_0}'}{z_{<i_0}}\bigg) \hat{a}^n 
\ge
\bigg(\frac{z'}{z} \frac{1}{(1-\ee_1)^n}- \frac{z_{<i_0}'}{z_{<i_0}}\bigg)
(\hat{a} -u_{i_0})^n.
\end{equation}
By (\ref{e:z}) and (\ref{e:zi0}), 
$|1- z'/z|$ and $|1- z_{<i_0}'/z_{<i_0}|$ can be made arbitrarily small
by taking $\zz$ small.
(\ref{equationlast}) thus shows that for $\zz$ small enough, 
$1 - u_{i_0} / \hat{a}$ must be small, i.e.\ $i_0$ must be near to $l_0-1$.
This concludes the proof of (\ref{equationclaim}').

\smallskip
Putting everything together, we see that $l_0$ and $a_0$ may be chosen
such that for $i$ not too near to $l_0-1$ 
\begin{equation*} 
\begin{aligned}   
|R_{2,i}| \; 
\overset{(\ref{equationclaim})}{>} \;
 (1- \ee_1) \rho |R_{2,2i(-1)}'| \;            
&  \overset{(\ref{equationrhon}), (\ref{e:z})}{>} \; 
(1-\ee_1) 4 \sqrt[n]{1-\zz}\, |R_{2,2i(-1)}'|            \\
& \overset{\;\;\;\;(\ref{e:x})\;\;\;\;}{>} \; 
4(1-\ee_1)^2 |R_{2,2i(-1)}'|                             \\
& \overset{\;\;\;\;\;\;\;\;\;\;\;\;\;\;}{>} \; (4- \ee) |R_{2,2i(-1)}'|.  
\end{aligned}
\end{equation*}
This proves
(\ref{e:200}.2.0).  

\medskip  
Suppose now that 
(\ref{e:200}.h.m), 
$1 \le h \le 4$, and hence also (\ref{e:201}.h.m)
hold true for $m \le k-1$. (\ref{e:201}.h.m) and (\ref{e:zi}) imply that for
$i$ not too near to $l_k-1$
\begin{equation}   \label{e:1.}
_kz_i > 1 - \zz.
\end{equation}
The reasoning which implied (\ref{e:ss}) thus also shows that for $i$
as in (\ref{e:ss})
\begin{equation}   \label{e:2.}
_kc_{2^ki} > {_kc_{2^{k-1}i}'}.
\end{equation}
Since $l_0$ is large and $\zz$ is small, $_kc_{2^{k-1}i} - {_kc_{2^ki}}$
is small. We thus see that for $i$ not too near to $l_0-1$
\begin{equation}   \label{e:3.}
_kc_{i} > {_kc_{2i(-1)}'}
\end{equation}
almost holds true, and hence also
\begin{equation}   \label{e:4.}
|_kR_{2,2i(-1)}'| | < \left|_kR_{2,i}\right|
\end{equation}
almost holds true.
Next, observe that (\ref{e:zi}) and (\ref{e:1.}) imply that
$A_0(a_k,l_k)/A_0(a_0,l_0)$ is near to $4^k$. This and (\ref{e:Gg}) show that
\begin{equation}   \label{e:5.}
A_0(a_k,l_k) > 12 l_k \pi,
\end{equation}
and in the same way as we derived (\ref{e:3}) from (\ref{e:ss}) and
(\ref{e:Hh}) we may derive from (\ref{e:3.}) and (\ref{e:5.}) that 
\begin{equation}   \label{e:6.}
|_kR_{3,2i(-1)}'| < \left|_kR_{3,i}\right|
\end{equation}
almost holds true.
Finally, by (\ref{e:1.}),
we also have that for $i$ not too near to $l_k-1$
\begin{equation}   \label{e:7.}
|_kR_{4,2i(-1)}'| < \left|_kR_{4,i}\right|.
\end{equation}
We infer from (\ref{e:pp}), (\ref{e:4.}), (\ref{e:6.}) and (\ref{e:7.})
that 
\begin{equation*}   
\frac{|_kR_{h,2i(-1)}'|}{|_kP_{2i(-1)}'|} < 3 \frac{|_kR_{h,i}|}{|_kP_i|},     
\qquad 1 \le h \le 4,
\end{equation*}
i.e.
\[ 
_ky_i' > 1 - 3 \zz.
\]
Proceeding exactly as in the case $k = 0$ we thus
get that for $i$ not too near to $l_k -1$
\begin{equation}   \label{e:8.}
_kc_i > (1- \ee_1)  \rho_k \, _kc_{2i(-1)}',
\end{equation}
from which (\ref{e:200}.2.k)
follows in the same way as for $k = 0$.

\bigskip
{\bf{(R3)}}
Set
\[
D_i(a,l) = \sq^{n-1}(1) \setminus \sq^{n-1}(N_i k_i)
\]
and 
\begin{equation}   \label{equationxi}
W_i(a,l) = \Bigg]\sum_{j=1}^{i-1} b_j(a,l), \sum_{j=1}^i b_j(a,l)\Bigg[ \times
]0,1[ \times \Bigg(1 - \frac{\sum_{j=1}^{i-1} b_j(a,l)}{a}\Bigg) \tr^{n-1}(\pi),
\end{equation}
$1 \le i \le l-1$.
Moreover, let $C_i$ be the cube in the $y$-factor of the fibers which
will be deformed and let $K_i$ be the extra space in $P_i$ needed to move
$C_i$ along the $y_j$-directions, $j \ge 2$. 
Then,
\[
R_3(a,l) = { \ff_{a,l} \Bigg(\coprod_{i=1}^{l-1} W_i(a,l) \Bigg) 
\times D_i(a,l) } \, \cup \, \coprod_{i=1}^{l-2} K_i.
\]

We first of all observe that $K_i \subset \ff_{a,l}(W_i(a,l)) \times
C_i$ and that $|C_i| / |D_i(a,l)|$ is small for $i$ not too near to
$l-1$ and $a$ large, since then $k_i(a,l)$ is large. We thus may forget
about the $K_i$. 
Next, as in ($R_2$), notice that for $i$ not too near to $l_k-1$, 
\[
| {_kR_{3,2i-1}'} \cap {_kT_i}| / | {_kR_{3,2i-1}'}| \ra 0 \quad \mbox{ as }\;
l_k \ra \infty, 
\]
whence we may neglect $_kR_{3,2i-1}' \cap {_kT_i}$ and prove
(\ref{e:200}.3.k)
with $_kR_{3,2i-1}'$ replaced by ${_kR_{3,2i-1}'} \setminus {_kT_i}$
(which we denote again by $_kR_{3,2i-1}'$).

By (\ref{equationsystem1}),
\begin{eqnarray} \label{e:Bb}
N_i(a,l) = 
\left\{ \begin{array}{ll}
           k_i -2,    & \;\; (k_i \mbox{ even}) \\
           k_i -3,    & \;\; (k_i \mbox{ odd})
        \end{array}
   \right. 
\quad \mbox{ for } 1 \le i \le l-2.
\end{eqnarray}
This and \mbox{Figure \ref{figure51}} show that for these $i$, 
\begin{equation}  \label{e:Cc}
\bigg( 1-\frac{3}{k_i(a,l)} \bigg) (n-1) 
\bigg( 1-\frac{N_i(a,l)}{k_i(a,l)}\bigg)        <
|D_i(a,l)| < (n-1)\bigg(1- \frac{N_i(a,l)}{k_i(a,l)}\bigg).
\end{equation}
Observe now that $c_i k_i = c_{2i}' k_{2i}' < c_{2i-1}' k_{2i-1}'$.
Hence, by (\ref{equationclaim}), 
\begin{equation}   \label{e:Dd}
k_{2i(-1)}' > (1-\ee_1) \rho k_i 
\end{equation}
if $i$ is not too near to $l_0-1$.
(\ref{e:Bb}) and (\ref{e:Dd}) imply that for these $i$ 
\begin{equation}   \label{e:Ee}
\frac{1 - N_i/k_i}{1 - N_{2i(-1)}'/k_{2i(-1)}'}
> \frac{2}{3} (1 -\ee_1) \rho. 
\end{equation}
Using again that for $i$ not too near to $l-1$, $k_i(a,l)$ is large whenever $a$
is large, we conclude from (\ref{e:Cc}) and (\ref{e:Ee}) that for $a_0$
large enough and $i$ not too near to $l_0-1$,
\begin{equation}   \label{e:Ff}
\frac{|D_i|}{|D_{2i(-1)}'|}
> \frac{2}{3} (1 - 2 \ee_1) \rho. 
\end{equation}

We conclude that for such $a_0$ and $i$ 
\begin{equation}  \label{e:zz}
\begin{aligned}   
|R_{3,i}|/|R_{3,2i(-1)}'|  \; 
\overset{(\ref{e:Hh}),(\ref{e:Ff})}{>} \; 
2 \frac{5}{6} \frac{2}{3}(1- 2 \ee_1) \rho  \;    
&  > \;                        
\frac{10}{9} (1-2\ee_1)4(1- \ee_1) \;  
> \; 4 -\ee.
\end{aligned}
\end{equation}
This proves (\ref{e:200}.3.0).

\smallskip
Suppose again that (\ref{e:200}.h.m), $1 \le h \le 4$, holds true for $m \le
k-1$. Then (\ref{e:8.}) implies 
\[
_kk_{2i(-1)}' > (1 - \ee_1) \rho_k \, _kk_i
\]
if $i$ is not too near to $l_k -1$, and proceeding as before we obtain
(\ref{e:200}.3.k).

\bigskip
{\bf{(R4)}}
Recall that $R_4(a,l) = \ff_{a,l}(Y(a,l))$ (cf.\ \mbox{Figure \ref{figure57}}).

\smallskip
To any partition $\bar{Z} = \coprod_{i=1}^{l-1} \bar{Z}_i$ of 
$E^{2n}(\pi, \dots, \pi,\bar{a})$ looking as in \mbox{Figure \ref{figure57}}
associate the set $X(\bar{Z}) = \coprod X_i(\bar{Z})$ which is obtained
from $\bar{Z}$ by replacing each fiber in $\bar{Z}_i$ by the maximal
fiber in $\bar{Z}_i$ (see \mbox{Figure \ref{figure57}}). 
Set $Y_i (\bar{Z}) = X_i(\bar{Z}) \setminus \bar{Z}_i$ and $Y(\bar{Z}) =
\coprod Y_i (\bar{Z})$. 
Clearly, if the partitions $E^{2n}(\pi, \dots, \pi, \bar{a}) =
\coprod_{i=1}^{l-1} \bar{Z}_i$ and $E^{2n}(\pi, \dots, \pi, \Bar{\Bar{a}}) =
\coprod_{i=1}^{l-1} \Bar{\Bar{Z}}_i$ are similar to each other, then
\begin{equation}   \label{e:alpha}
\frac{|Y_i (\bar{Z})|}{|\bar{Z}_i|}  =
\frac{|Y_i (\Bar{\Bar{Z}})|}{|\Bar{\Bar{Z}}_i|}.
\end{equation}
Let $B^{2n}(\bar{A}) = \coprod_{i=1}^l \bar{P}_i$ be a partition as in
\mbox{Figure \ref{figure56}} and assume that 
\[
\frac{|\bar{Z}_i|}{|\bar{P}_i|} > 1 - \zz \quad
\mbox{ and }  \quad
\frac{|\Bar{\Bar{Z}}_i|}{|\bar{P}_i|} > 1 - \zz
\quad \mbox{ for } 1
\le i \le i_0.
\]
Clearly, if $\zz$ is small enough and $i_0$ is large enough, $\bar{Z}$ and
$\Bar{\Bar{Z}}$ are almost similar. (\ref{e:alpha}) thus shows that
given $i_1$ not too large we may find $\zz$ and $i_0$ such that for $i
\le i_1$
\begin{equation}   \label{e:beta}
\frac{|Y_i (\bar{Z})|}{|\bar{Z}_i|} 
< (1 + \ee_1)
\frac{|Y_i (\Bar{\Bar{Z}})|}{|\Bar{\Bar{Z}}_i|}.
\end{equation}

Given $\hat{a}_m \in I_m(a_0)$, $m \in \NN_0$, and $1 \le i \le l_0-1$, set 
\[
Z_i(\hat{a}_m) = \coprod_{j = 2^m(i-1)+1}^{2^mi} Z_j(\hat{a}_m, l_m),
\]
$Z(\hat{a}_m) = \coprod Z_i(\hat{a}_m)$, $P(Z_i(\hat{a}_m)) = \coprod_{j
= 2^m(i-1)+1}^{2^mi} P(Z_j(\hat{a}_m,l_m))$ and $z(Z_i(\hat{a}_m)) =
|Z_i(\hat{a}_m)|/|P(Z_i(\hat{a}_m))|$.
For $a_0$ large and $i$ as above we clearly have that for all $m \in
\NN_0$ and $\hat{a}_m \in I_m(a_0)$
\begin{equation}   \label{e:yy}
\frac{\left| \coprod_{j=2^m(i-1)+1}^{2^mi} Y(Z_j(\hat{a}_m,l_m))
\right|}
{|P(Z_i(\hat{a}_m))|} 
\le
\frac{|Y_i (\hat{a}, l_0)|}
{|P_i(\hat{a},l_0)|}.
\end{equation}
Assume now that for some $m$, $i$ not too near to $l_0-1$ and
$2^m(i-1)+1 \le j \le 2^m$
\begin{equation}   \label{e:EEE}
\frac{\big|R_{h,j} (\hat{a}_m, l_m)\big|}
{\big|P_j(\hat{a}_m,l_m)\big|} 
\le
\frac{1}{(2-\ee)^m}
\frac{|R_{h,i} (\hat{a}, l_0)|}{|P_i(\hat{a},l_0)|}, \qquad 1 \le h \le 3.
\end{equation}
(\ref{e:yy}) and (\ref{e:EEE}) in particular imply that for these $i$ 
\begin{equation}   \label{e:nu}
z(Z_i(\hat{a}_m)) \ge z_i.
\end{equation}
(\ref{e:nu}) and (\ref{e:beta}) imply that $l_0$ and $a_0$ may be chosen
such that for all $\hat{a}_m$, $\hat{a}_{m'}$ satisfying
(\ref{e:EEE}) and $i$ not too near to $l_0-1$
\begin{equation}   \label{e:FFF}
\frac{|Y_i (Z(\hat{a}_m))|}{|Z_i(\hat{a}_m)|} 
< (1 + \ee_1)
\frac{|Y_i (Z(\hat{a}_{m'}))|}{|Z_i(\hat{a}_{m'})|}.
\end{equation}
Suppose now that (\ref{e:200}.h.m), $1 \le h \le 4$, holds true for $m \le
k-1$. We then have shown in ($R_h$), $1 \le h \le 3$, that (\ref{e:EEE})
holds true for $m \le k+1.$ (\ref{e:FFF}) thus implies that for $i$ not
too near to $l_0-1$
\begin{equation*}   
\frac{|Y_i (Z(\hat{a}_{k+1}))|}{|Z_i(\hat{a}_{k+1})|} 
< (1 + \ee_1)
\frac{|Y_i (Z(\hat{a}_k))|}{|Z_i(\hat{a}_k)|},
\end{equation*}
and (\ref{e:nu}) with $m=k$ now shows that for these $i$
\begin{equation}   \label{e:a1}
\frac{|Y_i (Z(\hat{a}_{k+1}))|}{|P(Z_i(\hat{a}_{k+1}))|} 
< 
\frac{1 + \ee_1}{1 -\zz}
\frac{|Y_i (Z(\hat{a}_k))|}{|P(Z_i(\hat{a}_k))|}.
\end{equation}
Pick $\ee_2$ so small that 
\begin{equation}   \label{e:ee2}
\left(1-\frac{\ee}{4}\right) \frac{1+\ee_2}{1-\ee_2}
\frac{1+\ee_1}{1-\zz} <1.
\end{equation}
This is possible since
\[
\left(1-\frac{\ee}{4}\right) \frac{1+\ee_1}{1-\zz} 
\overset{(\ref{e:x})}{<}
\left(1-\frac{\ee}{4}\right) \frac{1+\ee_1}{1-\ee_1} <1.
\]
We will show that $l_0$ and $a_0$ can be chosen such
that for any $\hat{a}_m$ satisfying (\ref{e:nu}), $i$ not too near to
$l_0-1$ and $2^m(i-1)+1 \le j \le 2^mi$
\begin{equation}   \label{e:lambda}
(1-\ee_2) |Y(Z_i(\hat{a}_m))| < 4^m | Y_j(\hat{a}_m,l_m)| 
< (1 + \ee_2) |Y(Z_i(\hat{a}_m))|.
\end{equation}
The second inequality in (\ref{e:lambda}) with $m=k+1$,
(\ref{e:a1}), the first inequality in (\ref{e:lambda}) with $m=k$ 
and (\ref{e:ee2}) then imply (\ref{e:200}.4.k).

In order to prove (\ref{e:lambda}), pick some small $\zz_0 = \zz$ and assume
$l_0$ and $a_0$ to be so large that for all
$\hat{a} \in I_0(a_0)$, $z_i(\hat{a}, l_0) >
1 -\zz_0$ whenever $i$ is not too near to $l_0-1$. Write $\bar{a}$ or
$\Bar{\Bar{a}}$ for any $a \ge a_1$ which satisfies (\ref{e:nu}). Then
\begin{equation}   \label{e:dd}
z(Z_i(\bar{a})) > 1- \zz_0 
\end{equation}
if $i$ is not too near to $l_0-1$.
\begin{figure}[h] 
 \begin{center}
  \psfrag{u}{$u$}
  \psfrag{u2}{$u_{2^m(i-1)}$}
  \psfrag{uk}{$u_{2^m(i-1)+2^{m-1}}$}
  \psfrag{uki}{$u_{2^mi}$}
  \psfrag{um}{$u_M$}
  \psfrag{Z1}{$Z_1$}
  \psfrag{Z2}{$Z_2$}
  \psfrag{Y0}{$Y_0$}
  \psfrag{Y1}{$Y_1$}
  \psfrag{Y2}{$Y_2$}
  \psfrag{dd}{$\dd$}
  \leavevmode\epsfbox{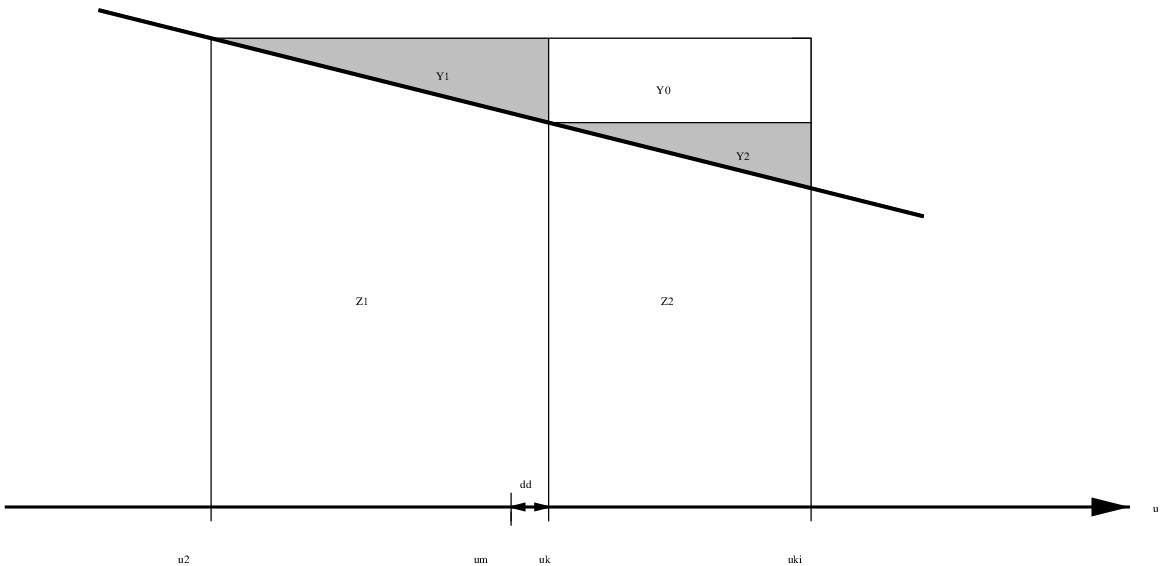}
 \end{center}
 \caption{$X_0$} \label{figure58}
\end{figure}
%
%
%
Fix once and for all such an $i$. Given $\hat{a}_m \in I_m(a_0)$, 
$m \in \NN$,
which satisfies (\ref{e:nu}),
set $d = u_{2^m i} - u_{2^m(i-1)}$, $u_M = u_{2^m(i-1)} + d/2$ and $\dd =
u_{2^m(i-1)+2^{m-1}} -u_M$, and write $Z_0 = Z_i(\hat{a}_m)$, 
$Z_1 = \coprod_{j = 2^m(i-1)+1}^{2^m(i-1)+2^{m-1}} Z_j(\hat{a}_m,l_m)$ and
$Z_2 = Z_0 \setminus Z_1$. Also write $X_j = X(Z_j)$,  $Y_j = Y(Z_j)$ and
$P_j = P(Z_j)$, $j = 0,1,2$ (see \mbox{Figure \ref{figure58}}).
Finally, define $R_h(Z_j)$, $1 \le h \le 4$, in the obvious way.

Define $\aa$, $\bb$ and $\gg_1$ by
\begin{eqnarray}\label{e:A}  
\frac{|X_1|}{|X_2|} = (1+\aa) \frac{d/2 + \dd}{d/2 - \dd}, 
\end{eqnarray}
\begin{eqnarray}\label{e:B}
|X_j| \le (1+\bb) |Z_j|, \qquad j = 1,2,
\end{eqnarray}
and
\begin{equation}   \label{e:C}
|P_1| = (1+\gg_1) |P_2|.
\end{equation}
We assume that $\bb$ is chosen minimal, and we observe that $\gg_1$ is
independent of $\hat{a}_m$ and is small since $i$ is not too near to $l_0-1$ and
$l_0$ is large. By (\ref{e:dd}), $|Z_0| >
(1-\zz_0)|P_0|$. This and (\ref{e:C}) readily imply that
\begin{equation}   \label{e:CC}
|Z_j| > (1- (2+ \gg_1) \zz_0) |P_j|, \qquad j=1,2.
\end{equation}
Thus, since $\gg_1 <1$, 
\begin{equation} \label{e:D}
\begin{aligned}   
(1+\aa) \frac{d/2 + \dd}{d/2-\dd} \; 
\overset{(\ref{e:A})}{=} \; 
\frac{|X_1|}{|X_2|} \;         &
\overset{(\ref{e:B})}{\ge} \;                   
\frac{|Z_1|}{(1+\bb)|Z_2|}          \\
&  \overset{(\ref{e:CC})}{>} \;
\frac{(1-3 \zz_0)|P_1|}{(1+\bb)|P_2|} \;
> \;                                           
\frac{1-3 \zz_0}{1+\bb} 
\end{aligned}
\end{equation}
and
\begin{equation} \label{e:E}
\begin{aligned}   
\frac{d/2 + \dd}{d/2-\dd} \; 
< \; 
\frac{|X_1|}{|X_2|} \; 
&\overset{(\ref{e:B})}{\le} \;                   
\frac{(1+\bb)|Z_1|}{|Z_2|}          \\
&  \overset{(\ref{e:CC})}{<} \;
\frac{(1+\bb)|P_1|}{(1-3 \zz_0)|P_2|} \;
\overset{(\ref{e:C})}{=} \;                                           
\frac{(1+\bb)(1+\gg_1)}{1- 3 \zz_0}. 
\end{aligned}
\end{equation}
If $\dd<0$, by (\ref{e:D}), 
\begin{align*}
d (\aa + \bb +\aa \bb + 3 \zz_0) &> 
|\dd|(4+2 \aa + 2\bb + 2 \aa \bb - 6 \zz_0),
\intertext{and if $\dd \ge 0$, by (\ref{e:E}),}
d(\gg_1 + \bb +\gg_1 \bb + 3 \zz_0) &> 
\dd (4+2 \gg_1 + 2\bb + 2 \gg_1 \bb - 6 \zz_0).
\end{align*}
Set $\mu = \max (\aa, \gg_1)$. Then
\begin{equation} \label{e:F}
|\dd| < \frac{d}{2} (\mu + \bb + 3 \zz_0)
\end{equation}
if $\zz_0$, $\bb$ and $\mu$ are small enough.

Set $c = \hat{a}_m - u_{2^m(i-1)}$.
Observe that, by (\ref{e:dd}), if $\zz_0$ is small,
$d(\bar{a})/d(\Bar{\Bar{a}})$ and
$c(\bar{a})/c(\Bar{\Bar{a}})$ are near to
$\bar{a}/\Bar{\Bar{a}}$ for all $\bar{a}$, $\Bar{\Bar{a}}$. 
Hence, $d(\bar{a})/c(\bar{a})$ is essentially independent of $\bar{a}$. Let
$\nu_1$ be such that $d(\bar{a})/c(\bar{a}) \le \nu_1$ for all $\bar{a}$. 
Since $c(\bar{a})$ is large for $i$ not too near to
$l_0-1$ and since $l_0$ is also large, $\nu_1$ is small.
Moreover, we readily compute
\begin{equation}  \label{e:AAA}
\aa = \frac{n-1}{2} \frac{d+2\dd}{c} + o\bigg(\frac{d}{c}\bigg)
\end{equation}
and
\begin{equation}  \label{e:BBB}
\bb = \frac{n-1}{4} \frac{d+2\dd}{c} + o\bigg(\frac{d}{c}\bigg).
\end{equation}
Thus, $\aa$ and $\bb$ are dominated by $\nu_1$, i.e.\ there are
small constants $\aa_1$ and $\bb_1$ such that $\aa \le
\aa_1$ and $\bb \le \bb_1$ for all $\bar{a}$. Set $\mu_1 = \max
(\aa_1, \gg_1)$.

Next, notice that $|Y_1|/|P_1|$ and $|Y_2|/|P_2|$ are essentially half
as large as $|Y_0|/|P_0|$ and hence also about half as large as
$|Y_i(\hat{a},l_0)|/ |P_i(\hat{a},l_0)|$. Indeed, 
\[
|Y_0| = \frac{1}{(n-1)!}\Big(\frac{\pi c}{\hat{a}_m}\Big)^{n-1}
\Bigg[d -\frac{c}{n} \Bigg(1-\bigg(1-\frac{d}{c}\bigg)^n\Bigg)\Bigg],
\]
and $|Y_1|$ respectively $|Y_2|$ are obtained from this expression by
replacing $d$ by $d/2+\dd$ respectively $c$ by $c-(d/2+\dd)$ and $d$ by
$d/2-\dd$. This yields
\begin{equation} \label{e:GGG}
\begin{aligned}   
\bigg|\frac{|Y_1|}{|Y_0|} - \frac{1}{4}\bigg| 
\; &\;= \;\; 
\frac{1}{4} \bigg|\frac{n-2}{6}\frac{d}{c} + 4 \frac{\dd}{d}\bigg| +
o \bigg(\frac{d}{c}\bigg) + o \bigg(\frac{\dd}{d}\bigg)            \\        
&\overset{(\ref{e:F})}{<}  \;
\frac{n}{2} \nu_1 + \mu_1 + \bb_1 + 3 \zz_0,
\end{aligned}
\end{equation}
and since $\nu_1$ is small, it turns out that the same estimate also holds true for $|Y_2|/|Y_0|$.
Moreover, (\ref{e:C}) implies that
\begin{equation}  \label{e:HHH}
\frac{|P_j|}{|P_0|} \ge \frac{1}{2+\gg_1}, \qquad j=1,2.
\end{equation}
If $\zz_0$, $\bb_1$, $\mu_1$ and $\nu_1$ and also $\ee$ are small
enough, we hence get
\begin{equation} \label{e:xx}
\begin{aligned}   
\frac{|Y_j|}{|P_j|}  \; \overset{(\ref{e:GGG}), (\ref{e:HHH})}{<}  \;
  \frac{3}{5} \frac{|Y_0|}{|P_0|}
& \; \overset{(\ref{e:FFF})}{<}  \;
\frac{3}{5} (1+ \ee_1 ) 
\frac{|Y_i(\hat{a},l_0)|}
              {|Z_i(\hat{a},l_0)|}    \\      
&\;\; < \;\;
\frac{3}{5} \frac{1 + \ee_1}{1-\zz_0} 
\frac{|Y_i(\hat{a},l_0)|}
              {|P_i(\hat{a},l_0)|}      \\
&\;\; < \;\; 
\frac{2}{3} \frac{|Y_i(\hat{a},l_0)|}
              {|P_i(\hat{a},l_0)|}, \qquad j = 1,2.
\end{aligned}
\end{equation}
We conclude that for $j=1,2$
\begin{equation} \label{e:III}
\begin{aligned}   
z(Z_j) \; = \; \frac{|Z_j|}{|P_j|} 
\;\;& = \;\; 1 - \frac{\sum_{h=1}^4 |R_h(Z_j)|}{|P_j|}           \\
\;\;& > \;\; 1 - \frac{\sum_{h=1}^3 |R_h(Z_j)| + |Y_j|}{|P_j|}   \\
& \!\!\!\!\!\!\overset{(\ref{e:EEE}), (\ref{e:xx})}{>}  \;
1 - \frac{2}{3} \frac{\sum_{h=1}^4 |R_{h,i}(\hat{a},
l_0)|}{|P_i(\hat{a}, l_0)|}   \\        
\;\;& > \;\; 1 - \frac{2}{3} \zz_0.
\end{aligned}
\end{equation}
In particular, $\zz_0$ in (\ref{e:dd}) may be replaced by $\zz_1 =
\frac{2}{3} \zz_0$.

We conclude that $l_0$ and $a_0$ may be chosen such that for all $\hat{a}_m$
\begin{equation}  \label{e:CCC}
(1-L_1) |Y_0| < 4 |Y_j| < (1+L_1) |Y_0|, \qquad j=1,2.
\end{equation}
Here, we put 
\[
L_1 = L(\zz_1, \bb_1, \mu_1, \nu_1) = 4(\mu_1 + \bb_1 + 3 \zz_1) + 2n \nu_1. 
\]
Observe that $L$ is linear
in $\zz_1$, $\bb_1$, $\mu_1$ and $\nu_1$.

\smallskip
Assume now that $m \ge 2$ and consider the partition $Z_1 = Z_1^2
\coprod Z_2^2$ whose components consist of $2^{m-2}$ consecutive
components of $Z(\hat{a}_m,l_m)$.
Set $d' = d/2 + \dd$ and define $\dd'$ to be the difference of the
$u$-width of $Z_1^2$ and $d'/2$. If $\aa'$ is defined by 
\begin{equation*} 
\frac{|X_1^2|}{|X_2^2|} = (1+\aa') \frac{d'/2 +\dd'}{d'/2 -\dd'}, 
\end{equation*}
we have
\begin{equation} \label{e:L}
\aa' =  \frac{n-1}{2} \frac{d' +2\dd'}{c} + o\bigg(\frac{d'}{c}\bigg).
\end{equation}
Since $\zz_1$ is small, $\dd'/d'$ is small. (\ref{e:AAA}) and (\ref{e:L}) thus show that
$\aa$ is near to $2\aa'$. In particular,
\begin{equation} \label{e:M}
\aa' < \frac{2}{3} \, \aa.
\end{equation}
Similarly, if $\bb'$ is the minimal constant with 
\begin{equation*} 
|X_j^2| \le (1+\bb')|Z_j^2|, \qquad j =1,2,
\end{equation*}
we have
\begin{equation} \label{e:MM}
 \begin{aligned}   
\bb' &= \frac{n-1}{4} \max \bigg( \frac{d' +2\dd'}{c},
        \frac{d'-2\dd'}{c-d'/2-\dd'}\bigg) + o\bigg(\frac{d'}{c}\bigg) \\
     &= \frac{n-1}{4} \frac{d' +2|\dd'|}{c} + o\bigg(\frac{d'}{c}\bigg),
 \end{aligned}
\end{equation}
and we conclude from (\ref{e:BBB}) and (\ref{e:MM}) as above that 
\begin{equation} \label{e:N}
\bb' < \frac{2}{3} \, \bb.
\end{equation}
A similar but simpler calculation shows that $\gg'$, which is defined by
$P(Z_1^2) = (1+\gg') P(Z_2^2)$, satisfies
\begin{equation} \label{e:O}
\gg' < \frac{2}{3} \, \gg_1.
\end{equation}
Next, since $\dd'/d$ is small, we also have that
\begin{equation} \label{e:P}
\frac{d'}{c} < \frac{2}{3} \, \nu_1.
\end{equation}
Consider now the partition $Z_2 = Z_3^2 \coprod Z_4^2$. While for $Z_1$
we had $c' =c$, now, $c'' = \hat{a}_m -u_{2^m(i-1)+2^{m-1}} = c-d'$.
But $c''/c = 1 - d'/c$ is near to $1$, whence the same arguments as above show
(\ref{e:M}), (\ref{e:N}), (\ref{e:O}) and (\ref{e:P}) with $\aa'$,
$\bb'$, $\gg'$ and $c'$ replaced by $\aa''$, $\bb''$, $\gg''$ and
$c''$. 
Finally, an argument analogous to the one which proved (\ref{e:III}) shows
$z(Z_j^2) > 1 - \frac{2}{3} \zz_1, \; 1 \le j \le 4.$ 
Summing up, we have shown that there are constants 
$\zz_2 = \frac{2}{3} \zz_1$, $\bb_2$,
$\mu_2$ and $\nu_2$ independent of $\hat{a}_m$ such that $L_2 = L(\zz_2, \bb_2,
\mu_2, \nu_2)$ satisfies $L_2 < \frac{2}{3} L_1$ and such that for all
$\hat{a}_m$
\[
(1-L_2)|Y_j| < 4 | Y_{2j(-1)}^2 | < (1+L_2)|Y_j|,
\qquad j=1,2.
\]

In general, let $Z^k(\hat{a}_m)$, $0 \le k \le m$, be the partition of $Z_0$
whose components consist of $2^{m-k}$ consecutive components of
$Z(\hat{a}_m,l_m)$.
Applying the above arguments to the components of $Z^k(\hat{a}_m)$, we see by
finite induction that there are constants $L_k, \; 1 \le k \le m$, with
$L_{k+1} < \frac{2}{3} L_k$ such that for all $\hat{a}_m$
\[
(1-L_{k+1}) |Y_j^k| < 4|Y_{2j(-1)}^{k+1}| < (1+L_{k+1}) |Y_j^k)|,
\]
$1 \le j \le 2^k$, $0 \le k \le m-1.$
Hence, with
\[
\pi_\pm(x) =  \prod_{k=1}^\infty \Bigg( 1 \pm \bigg(\frac{2}{3}\bigg)^k x \Bigg)
\]
we have that for all $j \in \{1, \dots, 2^m \}$
\begin{equation}  \label{e:S}
 \begin{aligned} 
\pi_-(L_1) |Y_0|  & \; < \; \prod_{k=1}^m (1- L_k) |Y_0|     \\
                        & \; < \; 4^m |Y_j^m|                      \\
                        & \; < \; \prod_{k=1}^m (1+ L_k) |Y_0|  
                          \; < \; \pi_+(L_1) |Y_0|.                           
 \end{aligned}
\end{equation}
Let $l_0$ and $a_0$ be so large that for $i$ not too
near to $l_0-1$, $L_1$ is so small that $1-\ee_2 < \pi_-(L_1)$ and
$\pi_+(L_1) < 1+\ee_2$. Then (\ref{e:S}) implies (\ref{e:lambda}).
This completes the proof of Proposition \ref{p:c}.
\proofend

\subsection{Lagrangian folding}  
\label{lagrangianfolding}

As already mentioned at the beginning of this section, 
there is a Lagrangian version of folding developed by Traynor in
\cite{T}. Here, the whole ellipsoid or the whole polydisc is viewed
as a Lagrangian product
of a cube and a simplex or a cube, and folding is then simply achieved
by wrapping the base cube around the base of the cotangent bundle of
the torus via a linear map. This version has thus a more algebraic
flavour. However, it yields good embeddings only for comparable shapes,
while the best embeddings of an ellipsoid into a polydisc respectively of
a polydisc into an ellipsoid via Lagrangian folding pack less than $1/n!$
respectively $n!/n^n$ of the volume.

\medskip 
For the convenience of the reader we review the method briefly.
\\
Write again $\RR^{2n}(x,y) = \RR^n(x) \times \RR^n(y)$ and set
\begin{gather*}
\sq (a_1, \dots, a_n) = \{ 0 < x_i < a_i, \; 1 \le i \le n \} \subset
\RR^n(x), \\
\tr (b_1, \dots, b_n) = \left\{ 0 < y_1, \dots, y_n \, \bigg| \, \sum_{i=1}^n
\frac{y_i}{b_i} <1 \right\} \subset \RR^n(y)
\end{gather*}
and 
\[
T^n = \RR^n(x)/ \pi \ZZ^n.
\]
The embeddings are given by the compositions of maps
\[
\begin{array}{rcl}
E(a_1 - \ee, \dots, a_n - \ee)  
& \xrightarrow{\aa_E}  & \sq^n (1) \times \tr (a_1, \dots, a_n)  \\
& \xrightarrow{\bb}    & \sq (q_1 \pi, \dots, q_n \pi) \times \tr (\frac{a_1}{q_1
  \pi}, \dots, \frac{a_n}{q_n \pi})  \\
& \xrightarrow{\gg}    & T^n \times \tr^n (\frac{A}{\pi}) \\
& \xrightarrow{\dd_E}  & B^{2n}(A)
\end{array}
\]
respectively
\[
\begin{array}{rcl}
P(a_1, \dots, a_n)  
& \xrightarrow{\aa_P}  &  \sq^n (1)
  \times \sq (a_1, \dots, a_n)  \\
& \xrightarrow{\bb}    &  \sq (q_1 \pi, \dots, q_n \pi) \times \sq (\frac{a_1}{q_1
  \pi}, \dots, \frac{a_n}{q_n \pi})  \\
& \xrightarrow{\gg}    &  T^n \times \sq^n (\frac{A}{\pi}) \\
& \xrightarrow{\dd_P}  &  C^{2n}(A),
\end{array}
\]
where $\ee >0$ is arbitrarily small and the $q_i$ are of the form $k_i$ or
$1/k_i$ for some $k_i \in \NN$.

$\aa_E$ and $\aa_P$ are the map $(x_1, y_1, \dots , x_n,y_n) \mapsto
(-y_1,x_1, \dots ,-y_n,x_n)$ followed by the maps described at the
beginning of section \ref{foldinginhigherdimensions}, and
$\bb$ is a diagonal linear map:
\[
\bb = \mbox{diag} \left[ q_1 \pi, \dots, q_n \pi, \frac{1}{q_1\pi}, \dots,
\frac{1}{q_n\pi} \right].
\]
Next, let 
\[
\tilde{\dd}_E \colon \sq^n (\pi) \times \tr^n
\left( \frac{A}{\pi} \right) \hookrightarrow B^{2n}(A)
\]
and
\[
\tilde{\dd}_P \colon \sq^n (\pi) \times \sq^n
\left( \frac{A}{\pi} \right) \hookrightarrow C^{2n}(A)
\]
be given by
\begin{eqnarray*}
(x_1, \dots, x_n, y_1, \dots, y_n)   &  \mapsto   &
                 (\sqrt{y_1} \cos 2x_1, \dots, \sqrt{y_n} \cos 2x_n,  \\
             & &  -\sqrt{y_1} \sin2x_1, \dots, -\sqrt{y_n} \sin 2x_n).
\end{eqnarray*}
Notice that $\tilde{\dd}_E$ respectively $\tilde{\dd}_P$ extend to an 
embedding of $T^n
\times \tr^n(A /\pi)$ respectively $T^n \times \sq^n(A/\pi)$. 
These extensions are the maps $\dd_E$ and $\dd_P$. 
We finally come to the folding map $\gg$.
\begin{lemma}  \label{lemmamatrix}
\begin{itemize}
\item[(i)]
If the natural numbers $k_1, \dots, k_{n-1}$ are relatively prime, then
\[
M(k_1, \dots, k_{n-1}) = \begin{pmatrix}
               1 &    &        &      &      -\frac{1}{k_1}      \\
                 & 1  &        &  0   &      -\frac{1}{k_2}      \\
                 &    & \ddots &      &      \vdots              \\
                 & 0  &        &  1   &      -\frac{1}{k_{n-1}}  \\
                 &    &        &      &          1               \\
           \end{pmatrix}      
\]
embeds $\sq (\pi/k_1, \dots, \pi/k_{n-1}, k_1 \dots k_{n-1}\pi)$ into $T^n$.
\item[(ii)]
For any $k_2, \dots, k_n \in \NN \setminus \{1\}$
\[
N(k_2, \dots, k_n) = \begin{pmatrix}
      1  &-\frac{1}{k_2}&        &        &               &         \\
         &    1    &-\frac{1}{k_3}&        & \Large{0}     &         \\
         &         & \ddots & \ddots &               &         \\
         &         &        & \ddots & -\frac{1}{k_{n-1}}&     \\
         &\Large{0}&        &        &  1            &-\frac{1}{k_n} \\
         &         &        &        &               & 1       \\
           \end{pmatrix}      
\]
embeds $\sq (\pi/(k_2 \dots k_n), k_2 \pi,\dots, k_n \pi)$ into $T^n$.
\end{itemize}
\end{lemma}

\proof
ad (i). Let $Mx=Mx'$ for $x,x' \in \sq (1/k_1, \dots, 1/k_{n-1}, k_1
\dots k_{n-1})$, so
\begin{equation} \label{equationai}
x_i - \frac{x_n}{k_i} = x_i' - \frac{x_n'}{k_i} + l_i, \qquad 1 \le i \le n-1
\end{equation}
for some $l_i \in \ZZ$ and 
\begin{equation} \label{equationbi}
x_n = x_n' + l_n,
\end{equation}
where $l_n \in \ZZ$ satisfies $|l_n| < k_1 \dots k_{n-1}$. Substituting
(\ref{equationbi}) into (\ref{equationai}) we get 
\begin{equation} \label{equationci}
x_i - x_i' = l_i + \frac{l_n}{k_i}, \qquad  1 \le i \le n-1.
\end{equation}
If $l_n=0$, we conclude $x=x'$. Otherwise, $|x_i-x_i'| < 1/k_i$ for $1
\le i \le n-1$ and (\ref{equationci}) imply that $l_n$ is an integral
multiple of all the $k_i$, whence by the assumption on the $k_i$ we have
$|l_n| \ge k_1 \dots k_{n-1}$, a contradiction.

\smallskip
ad (ii). Let $Nx=Nx'$ for $x,x' \in \sq (1/(k_2 \dots k_n), k_2, \dots,
k_n)$, so
\begin{equation} \label{equationaii}
x_i - \frac{x_{i+1}}{k_{i+1}} = x_i' - \frac{x_{i+1}'}{k_{i+1}}+ l_i, 
\qquad 1 \le i \le n-1
\end{equation}
for some $l_i \in \ZZ$ and 
\begin{equation} \label{equationbii}
x_n = x_n' + l_n.
\end{equation}
Substituting (\ref{equationbii}) into the last equation of (\ref{equationaii})
and resubstituting the resulting equations successively into the
preceding ones, we get
\begin{equation} \label{equationcii}
x_1 = x_1' + \frac{l_n}{k_2 \dots k_n} + \frac{l_{n-1}}{k_2 \dots
k_{n-1}} + \frac{l_{n-2}}{k_2 \dots k_{n-2}} + \dots + \frac{l_2}{k_2} +
l_1.
\end{equation}
Since $|x_1-x_1'| < 1/(k_2 \dots k_n)$, equation (\ref{equationcii}) has no
solution for $x_1 \neq x_1'$, hence $x_1 = x_1'$, and substituting
this into (\ref{equationaii}) and using $|x_i-x_i'| < k_i, \; 2 \le i
\le n$, we successively find $x_i = x_i'$.
\proofend

The folding map $\gg$ can thus be taken to be $M \times M^\ast$, where
$M$ is as in (i) or (ii) of the lemma and $M^\ast$ denotes the transpose
of the inverse of $M$.

\begin{remark}  \label{remarkproduct}
{\rm 
For polydiscs, the construction clearly commutes with taking
products. For ellipsoids, a similar compatibility holds: Let $M_1^\ast$
respectively $M_2^\ast$ be linear injections of $\tr (a_1, \dots ,a_m)$ 
into $\tr (a_1', \dots ,a_m')$ respectively $\tr (b_1, \dots ,b_n)$ into $\tr
(b_1', \dots ,b_n')$. Then $M_1^\ast \oplus M_2^\ast$ clearly injects
$\tr (a_1, \dots , a_m, b_1, \dots ,b_n)$ into $\tr (a_1', \dots ,a_m',
b_1', \dots ,b_n')$. Thus, given (possibly trivial) Lagrangian foldings
$\ll_1$ and $\ll_2$ which embed $E(a_1, \dots, a_m)$ into $E(a_1', \dots,
a_m')$ and $E(b_1, \dots, b_n)$ into $E(b_1', \dots, b_n')$, the Lagrangian
folding $\ll_1 \oplus \ll_2$ embeds $E(a_1, \dots, a_m, b_1, \dots
,b_n)$ into $E(a_1', \dots, a_m', b_1', \dots,b_n')$.
\diam
}
\end{remark}

In the following statements, $\ee$ denotes any positive number.
\begin{proposition}  \label{propositionlagrangianfolding}
\begin{itemize}
\item[(i)]
Let $k_1 < \dots < k_{n-1}$ be relatively prime and $a>0$. 
Then
\begin{itemize}
\item[$(i)_E$]
$
E^{2n}(\pi, \dots, \pi,a) \hookrightarrow B^{2n}(\max \{ (k_{n-1}+1) \pi,
                        \frac{a}{k_1 \cdots \, k_{n-1}} \} + \ee)
$
\item[$(i)_P$]
$
P^{2n}(\pi, \dots, \pi, a) \hookrightarrow C^{2n}(\max \{ k_{n-1} \pi,
(n-1)\pi + \frac{a}{k_1 \cdots \, k_{n-1}} \} ).
$
\end{itemize}
\item[(ii)]
Let $n \ge 3$, $k_2, \dots, k_n \in \NN \setminus \{1\}$ and $a_2, \dots, a_n >0$. Then
\begin{itemize}
\item[$(ii)_E$]
$
E(\pi, a_2, \dots, a_n) \hookrightarrow B^{2n}(A+\ee)$, where $A$ is found
as follows: Multiply the first column of $N^\ast$ by $k_2 \cdots k_n$ and
the i th column by $(a_i / \pi)/k_i,\; 2 \le i \le n$. Then add to
every row its smallest entry and add up the entries of each
column. $A/\pi$ is the maximum of these sums.
\item[$(ii)_P$]
$P(\pi, a_2, \dots, a_n) \hookrightarrow P(A_1, \dots , A_n)$, where
the $A_i$ are found as follows: Multiply $N^\ast$ as in
$(ii)_E$. $A_i/\pi$ is the sum of the absolute values of the entries of
the i th row.
\end{itemize}
\end{itemize}
\end{proposition}

\proof
ad (i). Write $y' = M^\ast (k_1, \dots ,k_{n-1})y$. We have 
\[
M^\ast(k_1, \dots, k_{n-1}) = \begin{pmatrix}
               1 &               &        &                   &         \\
                 &            1  &        &  0                &         \\
                 &               & \ddots &                   &         \\
                 &               &        &  1                &         \\
   \frac{1}{k_1} & \frac{1}{k_2} & \dots  & \frac{1}{k_{n-1}} & 1       \\
           \end{pmatrix} .     
\]
Thus, given $y \in \tr (k_1, \dots, k_{n-1}, \frac{a/\pi}{k_1 \cdots \,
k_{n-1}})$,
\begin{eqnarray*}
y_1' + \dots + y_n' & = &  (k_1+1) \frac{y_1}{k_1} + \dots +
         (k_{n-1}+1) \frac{y_{n-1}}{k_{n-1}} + \frac{a / \pi}{k_1 \cdots k_{n-1}}
         \frac{y_n}{\frac{a / \pi}{k_1 \cdots \, k_{n-1}}}                   \\
                    & < &  \max \left\{ k_{n-1} +1, \frac{a / \pi}{k_1 \cdots k_{n-1}} \right\},
\end{eqnarray*}  
and given $y \in \sq (k_1, \dots, k_{n-1}, \frac{a / \pi}{k_1 \cdots \, k_{n-1}})$,
\[
y' \in \sq (k_1, \dots, k_{n-1}, n-1 + \frac{a / \pi}{k_1 \cdots k_{n-1}}).
\]
ad (ii). We have
\[
N^\ast(k_2, \dots, k_n) = \begin{pmatrix}
      1             &             &         &        &           &         \\
\frac{1}{k_2}       &     1       &         &        &    0      &         \\
-\frac{1}{k_2 k_3}  &\frac{1}{k_3}&\ddots   &        &           &         \\
\vdots              &\vdots       &\ddots   & \ddots &           &         \\
\frac{(-1)^{n-1}}{k_2 \cdots \,k_{n-1}}&\frac{(-1)^{n-2}}{k_3 \cdots \,k_{n-1}}&\dots&\frac{1}{k_{n-1}}&1& \\
\frac{(-1)^n}{k_2 \cdots \,k_n}&\frac{(-1)^{n-1}}{k_3 \cdots \,k_n}&\dots&\frac{-1}{k_{n-1}k_n}&\frac{1}{k_n}&1 \\
           \end{pmatrix}.      
\]
Observe that we are free to compose $N^\ast$ with a
translation. Multiplying the columns as prescribed we get the vertices of
the simplex 
\[
N^\ast \tr \left(k_2 \dots k_n, \frac{a_2/\pi}{k_2}, \dots
,\frac{a_n/\pi}{k_n} \right). 
\]
Adding to the rows of this new matrix its
smallest entry corresponds to translating this new simplex into the
positive cone of $\RR^n(y)$. The claim thus follows.
A similar but simpler procedure leads to the last statement. 
\proofend
\\
\\
Proposition \ref{propositionlagrangianfolding} leads to the number
theoretic problem of finding appropriate relatively prime numbers $k_1,
\dots, k_{n-1}$. An effective method which solves this problem for $a$
large is described in the proof of Proposition \ref{p:refined} (i)$_E$.

\begin{corollary}  \label{corollarylagrangianfolding}
\begin{itemize}
\item[$(i)_E$]
$E^{2n} (\pi, l_{EB}(a), \dots , l_{EB}(a), a)  \hookrightarrow
B^{2n}(l_{EB}(a) + \ee)$,
where 
\[
l_{EB}(a) = \min_{k \in \NN} \max \{ (k+1)\pi, a/k \} = 
\left\{ \begin{array}{cll}
           (k+1)\pi, & (k-1)(k+1) \le a/\pi \le k(k+1) \\
           a/k, & k(k+1) \le a/\pi \le k(k+2).
        \end{array}
   \right. 
\]
\item[$(i)_P$]
$P^{2n} (\pi, l_{PC}(a), \dots , l_{PC}(a), a)  \hookrightarrow
C^{2n}(l_{PC}(a))$,
where 
\[
l_{PC}(a) = \min_{k \in \NN} \max \{ k \pi, a/k+\pi \} = 
\left\{ \begin{array}{cll}
           k \pi, & (k-1)^2 \le a/\pi \le k(k-1) \\
           a/k+\pi, & k(k-1) \le a/\pi \le k^2.
        \end{array}
   \right. 
\]
\end{itemize}
For $n \ge 3$ and any $k \in \NN \setminus \{1\}$
\begin{itemize}
\item[$(ii)_E$]
$E^{2n}(\pi, k^n \pi, \dots , k^n \pi) \hookrightarrow
B^{2n}((k^{n-1}+k^{n-2}+(n-2)k^{n-3})\pi + \ee)$,
\item[$(ii)_P$]
$P^{2n}(\pi, (k-1)k^{n-1} \pi, \dots , (k-1)k^{n-1} \pi) \hookrightarrow
C^{2n}(k^{n-1} \pi)$.
\end{itemize}
\end{corollary}

\proof
In (i)$_E$ and (i)$_P$ Remark \ref{remarkproduct} was
applied. For both (ii)$_E$ and (ii)$_P$ choose $k_2 = \dots = k_n =
k$. In (ii)$_E$, the maximal sum is the one of the entries of the $n-$1 st column, and in
(ii)$_P$ all the sums are $k^{n-1}$.
\proofend
\begin{figure}[h] 
 \begin{center}
  \psfrag{1}{$1$}
  \psfrag{2}{$2$}
  \psfrag{3}{$3$}
  \psfrag{4}{$4$}
  \psfrag{5}{$5$}
  \psfrag{6}{$6$}
  \psfrag{9}{$9$}
  \psfrag{A/pi}{$\frac{A}{\pi}$}
  \psfrag{a/pi}{$\frac{a}{\pi}$}
  \psfrag{spc/pi}{$\frac{s_{PC}}{\pi}$}
  \psfrag{lpc/pi}{$\frac{l_{PC}}{\pi}$}
  \psfrag{inclusion}{$\mbox{inclusion}$}
  \psfrag{v}{$\mbox{volume condition}$}
  \leavevmode\epsfbox{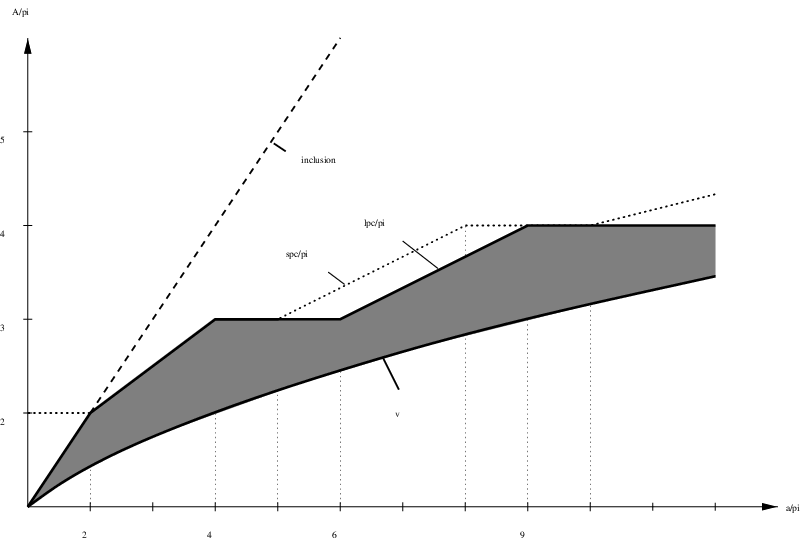}
 \end{center}
 \caption{What is known about $P(\pi,a) \hookrightarrow C^4(A)$} \label{figure21}
\end{figure}
%
\\
{\bf Examples.}
\\

ad (i)$_E$ and (i)$_P$.
Remark \ref{remarkproduct} and Proposition \ref{propositionlagrangianfolding}
$(i)$ applied to opposite entries imply that for any $k \in \NN$
\[
E^{2n} (\pi, k \pi, k^2 \pi, \dots, k^{2l} \pi) \hookrightarrow
B^{2n}((k^l+k^{l-1})\pi + \ee)
\]
and
\[
P^{2n} (\pi, k \pi, k^2 \pi, \dots, k^{2l} \pi) \hookrightarrow
C^{2n}((k^l+k^{l-1})\pi)
\]
if $n=2l+1$ is odd and
\[
E^{2n} (\pi, k^2 \pi, k^4 \pi, \dots, k^{2n-2} \pi) \hookrightarrow
B^{2n}((k^{n-1}+k^{n-2})\pi + \ee)
\]
and
\[
P^{2n} (\pi, k^2 \pi, k^4 \pi, \dots, k^{2n-2} \pi) \hookrightarrow
C^{2n}((k^{n-1}+k^{n-2})\pi)
\]
if $n$ is even.
\\

ad (ii)$_E$.
For $n=3$, Proposition \ref{propositionlagrangianfolding} yields 
\[
E(\pi, a_2, a_3) \hookrightarrow B^6 \left( \max \left\{ k_3(k_2+1)\pi,
\frac{a_2}{k_2k_3}(k_3+1)+\pi, \frac{a_3}{k_3}+\pi \right\} + \ee \right)
\]
for any $k_2, k_3 \in \NN \setminus \{1\}$. With
$(k_2,k_3)=(k,lk-1)$ we thus get for any $k \in \NN \setminus
\{1\}$ and $l \in \NN$
\[
E \left(\pi, \frac{k(lk-1)^2}{l}\pi, k(lk-1)^2 \pi \right) \hookrightarrow
B^6(k(lk-1)\pi+\pi + \ee).
\]

ad (ii)$_P$.
For $n=3$, Proposition \ref{propositionlagrangianfolding} yields 
\[
P(\pi, a_2, a_3) \hookrightarrow C^6 \left( \max \left\{ k_2 k_3 \pi,
k_3 \pi + \frac{a_2}{k_2}, \pi + \frac{a_2}{k_2 k_3}+ \frac{a_3}{k_3} \right\}
\right)
\]
for any $k_2, k_3 \in \NN \setminus \{1\}$. With
$(k_2,k_3)=(k,lk-l+1)$ we thus get for any $k \in \NN \setminus
\{1\}$ and $l \in \NN$
\[
P(\pi, (k-1)k(lk-l+1)\pi, l(k-1)k(lk-l+1) \pi) \hookrightarrow
C^6(k(lk-l+1)\pi).
\]
\hfill $\;\Diamond$

\subsection{Symplectic versus Lagrangian folding}  
                                 \label{symplecticversuslagrangianfolding}

For small $a$, the estimate $s_{EB}$ provides the best result known. 
For example,
we get $\frac{s_{EB}}{\pi} (4 \pi) = 2.6916\dots$, whence we have proved
\\
\\
{\bf Fact.} {\it $\, E(\pi, 4 \pi)$ embeds in $B^4(2.692 \, \pi)$.}
\\
\\
$l_{EB} (a) < s_{EB} (a)$ happens first at $a/\pi = 5.1622\dots$.
In general, computer calculations suggest that $l_{EB}$ and $s_{EB}$
yield alternately better estimates: For all $k \in \NN$ we seem to 
have that $l_{EB} < s_{EB}$ on an interval around $a = k(k+1)\pi$ and $s_{EB}
< l_{EB}$ on an interval around $k(k+2)\pi$; moreover, they suggest that
\[
\lim_{k \ra \infty} (s_{EB}(k(k+2)\pi) - l_{EB}(k(k+2)\pi)) =0,
\]
i.e.\ $l_{EB}$ and $s_{EB}$ seem to be asymptotically equivalent. We
checked the above statements for $k \le 5\,000$.

\begin{remark}  \label{remarkdeb}
{\rm
The difference $d_{EB}(a) = l_{EB}(a) - \sqrt{\pi a}$ between $l_{EB}$
and the volume condition attains local maxima at $a_k = k(k+2)\pi$,
where $d_{EB}(a) = (k+2)\pi - \sqrt{k(k+2)}\,\pi$. This is a decreasing
sequence converging to $\pi$.
\diam
}
\end{remark}

\smallskip
\mbox{Figure \ref{figure18}} summarizes the results. 
The non trivial estimates from below are provided by Ekeland-Hofer capacities,
which yield $A(a) \ge a$ for $a \in [\pi, 2 \pi ]$ and $ A(a) \ge 2 \pi$
for $a>2 \pi$.

\subsection{Summary} \label{summary}

Given $U \in \co(n)$ and $\aa >0$, set $\aa U = \{\aa z \in \CC^n \, |
\, z \in U\}$.

For $U,V \in \co(n)$ define squeezing constants
\[
s(U,V) = \inf \{ \aa \, | \, \mbox{there is a symplectic embedding } \ff
\colon U \hookrightarrow \aa V \}.
\]
Specializing, we define {\it squeezing numbers}
\[
s_{q_2 \dots q_n}^E (U) = s(U,E(1,q_2, \dots, q_n))
\]
and
\[
s_{q_2 \dots q_n}^P (U) = s(U,P(1,q_2, \dots, q_n)),
\]
and we write $s^B(U)$ for $s_{1\dots1}^E(U)$ and $s^C(U)$ for
$s_{1 \dots 1}^P (U)$.

With this notation, the main results of this section read
\begin{eqnarray} \label{equationinequalities}
s^B(E(\pi,a))  & \le & \min (s_{EB}(a), l_{EB}(a)) \\
s^B(P(\pi,a))  & \le & s_{PB}(a)                   \\
s^C(E(\pi,a))  & \le & s_{EC}(a)                   \\
s^C(P(\pi,a))  & \le & \min (s_{PC}(a), l_{PC}(a))  
\end{eqnarray}
and 
\[
\begin{array}{lcl}
s^C(P^{2n}(\pi,\dots ,\pi,a)) &\le &s_{PC}^{2n}(a)       
\end{array}
\]

\section{Packings}  \label{packings}

In the previous section we tried to squeeze a given simple shape into a
minimal ball and a minimal cube. This problem may be reformulated as follows:
\\
\\
{\it ``Given a ball $B$ respectively a cube $C$ and a simple shape $S$,
what is the largest simple shape similar to $S$ which fits into $B$
respectively $C$?''}

\medskip
\noindent
or equivalently:

\medskip
\noindent
{\it ``Given a ball or a cube, how much of its volume may be
symplectically packed by a simple shape of a given shape?''}
\\
\\
More generally, given $U \in \co (n)$ and any connected symplectic
manifold $(M^{2n},\oo)$, define the {\it $U$-width} of $(M,\oo)$ by
\[
w(U, (M, \oo)) = \sup \{ \aa \, | \; \mbox{there is a symplectic
embedding } \ff \colon \aa U \hookrightarrow (M, \oo) \},
\]
and if the volume $\mbox{Vol}(M, \oo) = \frac{1}{n!} \int_M \oo^n$ is
finite, set
\[
p(U, (M,\oo)) = \frac{|w(U,(M,\oo)) U|}{\mbox{Vol}(M,\oo)}.
\]
In this case, the two invariants determine each other,
$p(U,(M,\oo)) >0$ by Darboux's theorem, and if in
addition $n=1$, $p(U,(M,\oo)) =1$ by \nobreak{Theorem
\ref{t:volume}}.

Given real numbers $1 \le q_2 \le \dots \le q_n$, we define 
{\it weighted widths}
\begin{eqnarray*}
w_{q_2 \ldots q_n}^E (M, \oo) &=& w(E(1, q_2, \dots, q_n),(M,\oo)),
\\
w_{q_2 \ldots q_n}^P (M, \oo) &=& w(P(1, q_2, \dots, q_n),(M,\oo))
\end{eqnarray*}
and {\it packing numbers}
\begin{eqnarray*}
p_{q_2 \ldots q_n}^E (M,\oo) &=& p(E(1, q_2, \dots, q_n),(M,\oo)) 
= \frac{(w_{q_2 \ldots q_n}^E (M,\oo))^n q_2 \dots q_n}{n! \, \mbox{Vol}(M,\oo)},
\\
p_{q_2 \ldots q_n}^P (M,\oo) &=& p(P(1, q_2, \dots, q_n),(M,\oo)) 
= \frac{(w_{q_2 \ldots q_n}^P (M,\oo))^n q_2 \dots q_n}{\mbox{Vol}(M,\oo)}.
\end{eqnarray*}
Write $w(M, \oo)$ for the Gromov width $w_{1 \dots 1}^E(M, \oo)$ and
$p(M, \oo)$ for $p_{1 \dots 1}^E(M, \oo)$.
\begin{example}           \label{examplepacking}
{\rm
Assume that $(M, \oo) = (V, \oo_0) \in \co(n).$ 
By the very definitions of squeezing constants and
widths we have 
\[
w(U,V) = \frac{1}{s(U,V)}.
\]
In particular, we see that squeezing numbers and weighted widths of
simple shapes determine each other via
\begin{eqnarray} 
w_{q_2 \ldots q_n}^E(E(\pi, p_2 \pi, \dots , p_n \pi))  & = &
\frac{\pi^2}{s_{p_2 \ldots p_n}^E(E(\pi, q_2 \pi, \dots , q_n \pi))}, 
\label{equationwidth1}               \\
w_{q_2 \ldots q_n}^P(P(\pi, p_2 \pi, \dots , p_n \pi))  & = &
\frac{\pi^2}{s_{p_2 \ldots p_n}^P(P(\pi, q_2 \pi, \dots , q_n \pi))}, 
\label{equationwidth2}               \\
w_{q_2 \ldots q_n}^E(P(\pi, p_2 \pi, \dots , p_n \pi))  & = &
\frac{\pi^2}{s_{p_2 \ldots p_n}^P(E(\pi, q_2 \pi, \dots , q_n \pi))}, 
\label{equationwidth3}               \\
w_{q_2 \ldots q_n}^P(E(\pi, p_2 \pi, \dots , p_n \pi))  & = &
\frac{\pi^2}{s_{p_2 \ldots p_n}^E(P(\pi, q_2 \pi, \dots , q_n \pi))}. 
\label{equationwidth4}
\end{eqnarray}
Combined with the estimates stated in subsection \ref{summary},
these equations
provide estimates of weighted widths and packing numbers of simple
shapes from below.
\diam
}
\end{example}
If $(M,\oo)$ is an arbitrary symplectic manifold whose Gromov width 
is known to be large, these results may be used to
estimate $w_{q_2 \ldots q_n}^E(M,\oo)$ and $p_{q_2 \ldots q_n}^E(M,\oo)$ reasonably
well from below.
\\
\\
{\bf Example.}
Let $T^2(\pi)$ be the 2-torus of volume $\pi$ and $S^2(2\pi)$ the
sphere of volume $2\pi$ and endow $M = T^2(\pi) \times S^2(2\pi)$
with the split symplectic structure. Theorem \ref{t:gromov}(ii) shows
that $p(M) =1$. Thus, by (\ref{equationwidth1}) and 
(\ref{equationinequalities})
\begin{eqnarray*}
p_q^E(M) \ge p_q^E(B^4(2\pi)) & = & \frac{(w_q^E(B^4(2\pi)))^2 q}{4 \pi^2}\\
                       & = & \frac{q \pi^2}{(s^B(E(\pi, q \pi)))^2} \ge 
\frac{q\pi^2}{(\min(s_{EB}(q\pi),l_{EB}(q\pi)))^2}.
\end{eqnarray*}
In particular, $\lim_{q \ra \infty} p_q^E(M) =1$.
\diam
\\

On the other hand, $w(U, (M,\oo)) \ge w(V,(M,\oo))$ whenever $U \le_3
V$; 
in particular, $w \ge w_{q_2 \dots q_n}^E \ge w_{q_2 \dots q_n}^P$
for all $1 \le q_2 \le \dots \le q_n$.
Thus, if $w(M,\oo)$ and the weights are small, we get good estimates of
weighted widths and packing numbers from above.
\\
\\
{\bf Example.}
Let $r \ge 1$ and $M = S^2(\pi) \times S^2(r \pi)$ with the split symplectic
structure. By the Non-Squeezing Theorem stated at the beginning of
Appendix B we have 
$w(M) \le \pi$, whence $w_q^E(M) \le \pi$ and $p_q^E(M) \le
\frac{q}{2r}$. For $q \le r$ the obvious embedding $E(\pi, q\pi)
\hookrightarrow P(\pi, r\pi) \hookrightarrow M$ shows that these
inequalities are actually equalities. 
\diam
\\

The knowledge of the Gromov width is thus of
particular importance to us. Recently considerable progress has been
made in computing or estimating the Gromov width of closed 4-manifolds. 
An overview on these results is given in Appendix B. 
\\
\\
{\bf Remark.}
Since the Gromov width is the smallest symplectic capacity we might try
to estimate it from above by using other symplectic capacities. However,
other capacities (like the Hofer-Zehnder capacity or the first
Ekeland-Hofer capacity, Viterbo's capacity and the capacity arising 
from symplectic homology in the case of subsets of $\RR^{2n}$) 
are usually even harder to
compute. In fact, we do not know of any space for which a capacity
other than the Gromov width is known and finite while its Gromov width
is unknown.
\diam

\subsection{Asymptotic packings}  \label{asymptoticpackings}

%
%
\begin{theorem}  \label{t:volume}
Let $M^n$ be a connected manifold endowed with a volume form $\OO$ and
let $U \subset \RR^n$ be diffeomorphic to a standard ball. Then $U$
embeds in $M$ by a volume preserving map if and only if $|U| \le \Vol
(M, \OO)$.
\end{theorem}
\proof
Endow $\overline{\RR}_{>0} = \RR_{>0}\cup \{ \infty \}$ with the topology
whose base of open sets is given by joining the open intervals $]a,b[ \,
\subset
\RR_{>0}$ with the subsets of the form $]a, \infty] = \, ]a, \infty[ \,
\cup  \, \{\infty \}$. Denote the Euclidean norm on $\RR^n$ by $\| \cdot \|$
and let $S_1$ be the unit sphere in $\RR^n$. 
\begin{lemma}  \label{l:star}
Let $\RR^n$ be endowed with its standard smooth structure, let $\mu
\colon S_1 \ra \overline{\RR}_{>0}$ be a continuous function and let
\[
S = \left\{ x \in \RR^n \, \Big|\, x =0 \mbox{ or } \, 0 < \| x \| < 
\mu \left( \frac{x}{\| x \|} \right) \right\}
\]
be the starlike domain associated to $\mu$.
Then $S$ is diffeomorphic to $\RR^n$.
\end{lemma}
{\bf Remark.} The diffeomorphism guaranteed by the lemma may be chosen
such that the rays emanating from the origin are preserved.
\\
\\
{\it Proof of the lemma.} If $\mu (S_1) = \{ \infty \}$, there is nothing
to prove. For $\mu$ bounded, the lemma was proved by Ozols \cite{O}. 
If $\mu$ is neither bounded nor $\mu (S_1) = \{ \infty \}$,
Ozols's proof readily extends to our situation. Using his notation, the
only modifications needed are: Require in addition that $r_0 <1$ and
that $\ee_1 <2$, and define continuous functions $\tilde{\mu}_i \colon
S_1 \ra \RR_{>0}$ by 
\[
\tilde{\mu}_i = \min \{ i , \mu - \ee_i + \dd_i/2 \}.
\]
With these minor adaptations the proof in \cite{O} applies word by word.
\proofend
\\
Next, pick a complete Riemannian metric $g$ on $M$. (We refer to \cite{KN} for
basic notions and results in Riemannian geometry.) The existence of such
a metric is guaranteed by a theorem of Whitney \cite{W}, according to
which $M$ can be embedded as a closed submanifold in some $\RR^m$. We
may thus take the induced Riemannian metric. A direct and elementary
proof of the existence of a complete Riemannian metric is given in
\cite{NO}.
Fix a point $p \in M$, let $\exp_p \colon T_pM \ra M$ be the
exponential map at $p$ with respect to $g$, let $C(p)$ be the cut locus
at $p$ and set $\widetilde{C}(p) = \exp_p^{-1}(C(p))$.
Let $S_1$ be the unit sphere in $T_pM$, let $\mu_p \colon S_1 \ra
\overline{\RR}_{>0}$ be the function defining $\widetilde{C}(p)$ and let
$S_p \subset T_pM$ be the starlike domain defined by
$\widetilde{C}(p)$. Since $g$ is complete, $\mu_p$ is continuous \cite[p.\
98]{KN}. We are thus in the situation of Lemma \ref{l:star}, and since
$\exp_p(S_p) = M \setminus C(p)$ \cite[p.\ 100]{KN}, we obtain
\begin{corollary}  \label{c:full}
Let $(M^n,g)$ be a complete Riemannian manifold. Then the maximal normal
neighbourhood $M \setminus C(p)$ of any point $p$ in $M$ is
diffeomorphic to the standard $\RR^n$. 
\end{corollary}
Using polar coordinates on $T_pM$ we see from Fubini's Theorem that
$\widetilde{C}(p)$ has zero measure; thus the same holds true for $C(p)$,
whence 
\[
\Vol (S_p, \exp_p^* \OO) = \Vol (M \setminus C(p) ,\OO)= \Vol (M,\OO).
\]
Theorem \ref{t:volume} now follows from Lemma \ref{l:star} and
\begin{proposition} {\rm (Greene-Shiohama, \cite{GS}) \rm}
\label{propositiongreeneshiohama}
Two volume forms $\OO_1$ and $\OO_2$ on an open manifold are diffeomorphic
if and only if the total volume and the set of ends of infinite volume
are the same for both forms. 
\end {proposition}
\proofend
{\bf Remark.}
The existence of a volume preserving embedding of a set $U$ as above
with $|U| < \Vol(M,\OO)$ immediately follows from Moser's deformation
technique if $M$ is closed and from 
Proposition \ref{propositiongreeneshiohama}, which is itself an
extension of that technique to open manifolds, if $M$ is open. The main
point in Theorem \ref{t:volume}, however, is that {\it all} of the
volume of $M$ can be filled. This is in contrast to the full
symplectic packings by $k$ balls established in \cite{MP}, \cite{B1} and
\cite{B2}.
\diam
\\
In view of the Non-Squeezing Theorem and the existence of symplectic
capacities, very much in contrast to the volume-preserving case, there
exist strong obstructions to full packings by ``round'' simple shapes in the
symplectic category. (We refer to the previous sections for related
results on embeddings into simple shapes and to Appendix B
for an overview on known results on the Gromov width of closed four
manifolds.)

However, the results of section \ref{flexibility} show for example that
for embeddings into four dimensional simple shapes packing obstructions
more and more disappear if we pass to skinny domains.
The main goal of this section is to show that in the limit rigidity
indeed disappears.
\begin{theorem} \label{theoremasymptotics}
Let $(M,\oo)$ be a connected symplectic manifold of finite
volume. Then
\[
p_\infty^E(M,\oo) = \lim_{q \ra \infty} p_{1 \dots 1 q}^E(M,\oo)
\quad \mbox{ and } \quad 
p_\infty^P(M,\oo) = \lim_{q \ra \infty} p_{1 \dots 1 q}^P(M,\oo)
\]
exist and equal $1$.
\end{theorem}
{\bf Remark.}
Remark \ref{remarkproduct}, Proposition
\ref{propositionlagrangianfolding}(i) and the theorem immediately
imply that for any $(M, \oo)$ as in the theorem 
\[
\lim_{q \ra \infty} p_{qq^2 \dots \, q^{n-1}}^E (M, \oo)
\quad \mbox{ and } \quad
\lim_{q \ra \infty} p_{qq^2 \dots \, q^{n-1}}^P (M, \oo)
\]
exist and equal $1$.                                     \hfill \diam
\\

\begin{figure}[h] 
 \begin{center}
  \psfrag{aP}{$\aa P^{2n}(\pi, \dots ,\pi, a)$}
  \psfrag{L1}{$L_1$}
  \psfrag{L2}{$L_2$}
  \psfrag{L3}{$L_3$}
  \psfrag{C1}{$\overline{C}_1(\ll)$}
  \psfrag{C2}{$\overline{C}_2(\ll)$}
  \psfrag{C3}{$\overline{C}_3(\ll)$}
  \psfrag{C4}{$\overline{C}_4(\ll)$}
  \psfrag{wl}{$\sqrt{\ll}$}
  \psfrag{e2}{$\ee_2$}
  \psfrag{x1}{$x_1$}
  \psfrag{psi}{$\psi$}
  \psfrag{M}{$M$}
  \leavevmode\epsfbox{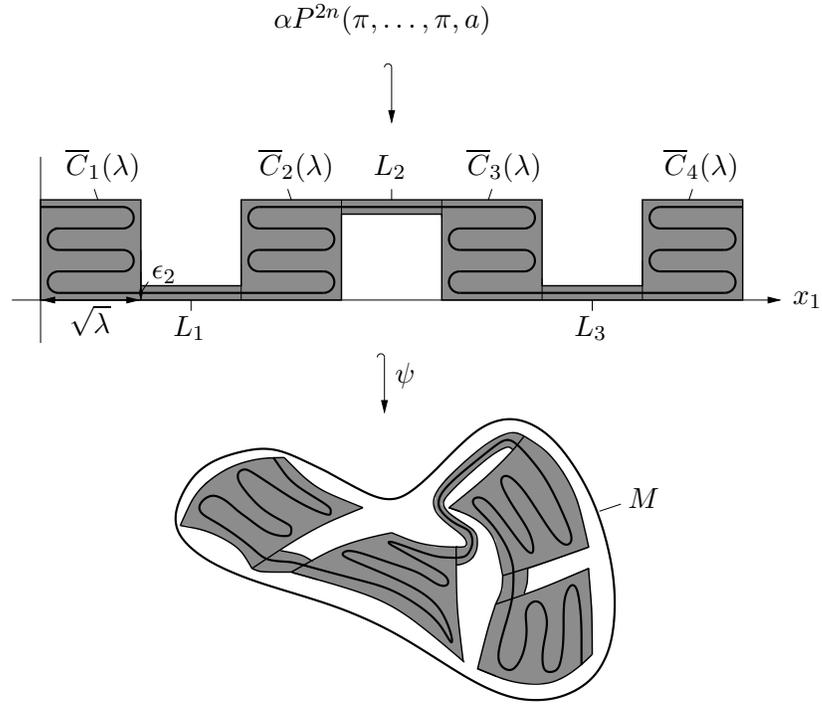}
 \end{center}
 \caption{Asymptotic filling by polydiscs} \label{figure27}
\end{figure}

The proof of the statement for polydiscs  proceeds along the following
lines: We first fill $M$ up to
some $\ee$ with small disjoint closed cubes, which we connect by 
lines. We already know how to asymptotically fill these cubes with thin
polydiscs, and we may use neighbourhoods of the lines to pass from one
cube to another (cf.\ \mbox{Figure \ref{figure27}}). 

The case of ellipsoids is less elementary.
For $n \le 3$, the statement for ellipsoids follows from the one for
polydiscs and the fact that a polydisc may be asymptotically filled by
skinny ellipsoids. This is proved in the same way as
(\ref{equationbadlimes}). 
In higher dimensions, however, symplectic folding alone is not powerful
enough to fill a polydisc by thin ellipsoids, since there is no
elementary way of filling a cube by balls.
However, algebro-geometric methods imply that in any dimension cubes can
indeed be filled by balls. Using this, we may almost fill $(M,\oo)$ by
equal balls, which we connect again by thin lines. The claim then
readily follows from the proof of Proposition \ref{propositionlim=1}.

\medskip

We begin with the following

\begin{lemma}  \label{lemmamcduffpolterovich}
{\rm (McDuff-Polterovich, \cite{MP})}
Let $(M, \OO)$ be a symplectic manifold of finite volume. Then, given
$\ee >0$, there is an embedding of a disjoint union of closed equal
cubes $\coprod \overline{C}(\ll)$ into $M$ such that $|\coprod
C(\ll)|
> \Vol\,(M) -2\ee$.
\end{lemma}

\proof
Assume first that $M$ is compact and cover $M$ with Darboux charts $V_i
= \ff_i(U_i), \; i=1,\dots,m$.
Pick closed cubes $\overline{C}_1, \dots, \overline{C}_{j_1} \subset U_1$
of possibly varying size such that 
\[
\sum_{j=1}^{j_1} |C_j| > \Vol (V_1) - \frac{\ee}{m}.
\] 
Proceeding by finite induction, for $i>1$, set $k_i = \sum_{l=1}^{i-1}
j_l$ and pick closed cubes $\overline{C}_{k_i+1}, \dots,
\overline{C}_{k_i+j_i} \subset U_i \setminus \ff_i^{-1}
(\bigcup_{j=1}^{i-1} V_j)$ such that
\[
\sum_{j=1}^{j_i} |C_{k_i+j}| > \Vol (V_i \setminus
\bigcup_{j=1}^{i-1} V_j) - \frac{\ee}{m}.
\]
Choose now $\ll$ so small that all the cubes $\overline{C}_k, \; 1 \le k
\le k_{m+1}$, admit an embedding of a disjoint union $\coprod_{j=1}^{n_k}
\overline{C}(\ll)$ such that 
$n_k |C(\ll)| > |C_k| - \ee/k_{m+1}$.
In this way, we get an embedding of $\sum_{k=1}^{k_{m+1}} n_k$ closed
cubes into $M$ filling more than $\Vol (M) -2\ee$.

If $M$ is not compact, choose a volume-preserving embedding $\ff \colon$
\newline
$\overline{B^{2n}} (\Vol (M) -\ee) \hookrightarrow M$ (cf.\ Theorem
\ref{t:volume}) and apply the already proved part to
$\mbox{$(\overline{B^{2n}} (\Vol (M) -\ee), \ff^* \oo)$}$.
\proofend
%
%

We next connect the cubes by thin lines.

Pick $\ee_1 >0$ and let $\ff = \coprod_{i=1}^k \ff_i \colon \coprod_{i=1}^k
\overline{C}_i(\ll) \hookrightarrow M$ be a corresponding embedding
guaranteed by Lemma \ref{lemmamcduffpolterovich}. Extensions of the $\ff_i$
to small neighbourhoods of $\overline{C}_i(\ll)$ are still denoted by
$\ff_i$. We may assume that the faces of the $\overline{C}_i(\ll)$ are
cubes and that all the $\overline{C}_i(\ll)$ lie in the positive cone of
$\RR^{2n}$ and touch the $x_1$-axis. Join these cubes by straight lines
$L_i$ as described in \mbox{Figure \ref{figure27}}, i.e.\ fixing regular
parameterizations $L_i(t) \colon [0,1] \ra L_i$ we have 
$L_i(0) \in \partial \overline{C}_i(\ll), \; L_i(1) \in \partial \overline{C}_{i+1}(\ll)$
and
\[
L_i(t) =
\left\{ \begin{array}{ll}
         (x_1(L_i(t)), 0,\dots,0)     &   \mbox{for $i$ odd,}    \\
         (x_1(L_i(t)), \sqrt{\ll},\dots,\sqrt{\ll})   &  \mbox{for $i$ even.} 
        \end{array}
   \right. 
\]
Let now $\coprod_{i=1}^{k-1} \ll_i \colon \coprod L_i \ra M \setminus \coprod
\ff_i (C_i(\ll))$
be a disjoint family of embedded curves in $M$ which touches $\coprod
\ff_i (\overline{C}_i(\ll)$ only at the points $\ll_i(0)$ and $\ll_i(1)$
and coincides with $\ff_{i|L_i}$ respectively $\ff_{i+1|L_{i+1}}$ on a
small neighbourhood of $\overline{C}_i(\ll)$ respectively
$\overline{C}_{i+1}(\ll)$.
Choose 1-parameter families of symplectic frames $\{ e_{j,i}(t)
\}_{j=1}^{2n}$ respectively $\{ e_{j,i}'(t)\}_{j=1}^{2n}$ along $L_i(t)$
respectively 
$\ll_i(L_i(t))$ such that $e_{1,i}(t) = \frac{d}{dt} \ll_i(t)$ and
$e_{1,i}'(t) = \frac{d}{dt} \ll_i(L_i(t))$. 
Let $\tilde{\psi}_i$ be an extension of $\ll_i$ to a
neighbourhood of $L_i$ which coincides with $\ff_i$ respectively $\ff_{i+1}$ on
a neighbourhood of $\ll_i(0)$ respectively $\ll_i(1)$ and which sends the
symplectic frame along $L_i(t)$ to the one along $\ll_i(L_i(t))$, i.e.
\[ \left(T_{L_i(t)} \tilde{\psi}_i\right) (e_{j,i}(t)) = e_{j,i}'(t). \]
$\tilde{\psi}_i$ is thus a diffeomorphism on a neighbourhood of $L_i$
which is symplectic along $L_i$.
Using a variant of Mosers's method (see \cite[Lemma 3.14 and its
proof\,]{MS}) we see that $\tilde{\psi}_i$ may be deformed to an embedding
$\psi_i$ of a possibly smaller neighbourhood of $L_i$ which still
coincides with $\ll_i$ on $L_i$ and $\ff_i$ respectively $\ff_{i+1}$ on a
neighbourhood of $L_i(0)$ respectively $L_i(1)$, but is symplectic
everywhere. Choose $\ee_2 >0$ so small that for all $i$, $\psi_i$ is defined
on $N_i(\ee_2) = \{x_1(L_i(t))\} \times [0,\ee_2]^{2n-1}$ if $i$ is odd
and on $N_i(\ee_2) = \{x_1(L_i(t))\} \times [\sqrt{\ll} -\ee_2,
\sqrt{\ll}]^{2n-1}$ if $i$ is even.

Summing up, we see that there exists $\ee_2 >0$ such that  
\[ \cn (\ee_2) = \coprod C_i(\ll) \coprod N_i(\ee_2) \] 
symplectically embeds in $M$. 

It remains to show that
$\cn (\ee_2)$ may be asymptotically filled by skinny polydiscs. We try
to fill $\cn(\ee_2)$
by $\aa P^{2n}(\pi, \dots,\pi, a)$ with $\aa$ small and $a$ large by
packing the $C_i(\ll)$ as described in subsection
\ref{embeddingsofpolydiscs} and using $N_i(\ee_2)$ to pass from
$C_i(\ll)$ to $C_{i+1}(\ll)$. Here we think of $\aa P^{2n}(\pi,
\dots,\pi, a)$ as $\frac{\aa^2}{\ee_2} \sq (a, \pi, \dots, \pi) \times
\sq (\ee_2, \dots, \ee_2)$ and of $C^{2n}(\ll)$ as $\frac{1}{\ee_2} \sq
(\ll, \dots, \ll) \times \sq (\ee_2, \dots, \ee_2)$.
Write $P_i$ for the restriction of the image of $\aa P^{2n}(\pi,
\dots,\pi, a)$ to $C_i(\ll)$.
In order to guarantee that the ``right'' face of $P_i$ and the ``left''
face of $N_i(\ee_2)$ fit, we require that the number of folds in each
$z_1$-$z_2$-layer is even
and that the component of $P_i$ between its right face and the last
stairs touches $\partial \overline{C}_i(\ll)$ wherever possible. This second point
may be achieved by making $n-$1 of the stairs in $P_i$ a little bit
higher than necessary. The part of the image of $\aa P^{2n}(\pi,
\dots,\pi, a)$ between $P_i$ and $P_{i+1}$ will thus be contained in
$N_i(\ee_2)$ whenever $\aa^2 \pi < \ee_2^2$.

Now, in Proposition \ref{propositionhigherintocubes} we have
\[
\lim_{a \ra \infty} \frac{a \pi^{n-1}}{(s_{PC}^{2n}(a))^n} =1,
\]
and hence, by duality,
\begin{equation} \label{equationlim}
\lim_{q \ra \infty} p_{1\dots 1q}^P(C^{2n}(\ll)) =1.
\end{equation}
(\ref{equationlim}) is clearly not affected by the two minor
modifications which we required above for the packing of $C_i(\ll)$.
Thus half of the theorem follows.

\medskip

As explained above, in order to prove the statement for ellipsoids we
need the following non-elementary result.

\begin{proposition}    \label{propositionMP}
{\rm (McDuff-Polterovich, \cite[Corollary 1.5.F]{MP})}
For each positive integer $k$, arbitrarily much of the volume of
$\overline{C^{2n}}(\pi)$ may be filled by $n!\,k^n$ equal closed balls.
\end{proposition}

This proposition may be proved in two different ways, either via symplectic
blowing up and fibrations or via symplectic branched
coverings. Combining it with Lemma \ref{lemmamcduffpolterovich}, we see
that we may fill as much of the volume of $(M, \oo)$ by disjoint equal
closed balls as we want.

\begin{figure}[h] 
 \begin{center}
  \psfrag{0}{$0$}
  \psfrag{l}{$\ll$}
  \psfrag{l+1}{$\ll+1$}  
  \psfrag{m+1l+m}{$(m+1) \ll+m$}
  \psfrag{x1}{$x_1$}
  \psfrag{B0}{$B_0(\ll)$}
  \psfrag{B1}{$B_1(\ll)$}
  \psfrag{Bm}{$B_m(\ll)$}
  \psfrag{e}{$\ee$}
  \psfrag{dots}{$\dots$}
  \leavevmode\epsfbox{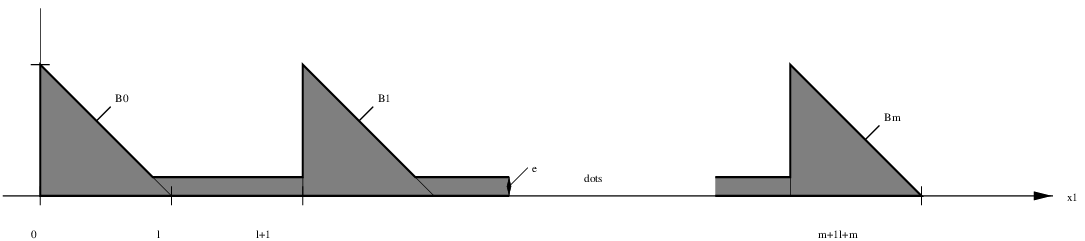}
 \end{center}
 \caption{$\cn (\ee)$} \label{figure59}
\end{figure}

So assume that $(M, \oo)$ is almost filled by $m+1$ disjoint equal
closed balls
$\left( \overline{B}_i(\ll), \ff_i \right)$,
$0 \le i \le m$. By Lemma
\ref{lemmafibres}(ii) we may think of $B_i(\ll)$ as fibered over 
$]i \ll+i, (i+1)\ll+i[ \, \times \, ]0,1[$ with fibers 
$\gg \tr^{n-1}(\ll) \times \sq^{n-1}(1)$, $1 \ge \gg  >0$ 
(cf.\ \mbox{Figure \ref{figure59}}). 
%
%
Exactly as in the case of cubes we find an $\ee >0$ such that $\ff =
\coprod_{i=0}^m \ff_i$ extends to a symplectic embedding $\psi$ of a
small neighbourhood of 
\[ \cn (\ee) = \coprod B_i (\ll) \bigcup \; ]0, m\ll +m[ \, \times \,
]0, \ee[^{2n-1}. \]
Let $\tau_i \colon \RR^{2n} \ra \RR^{2n}$, $z \mapsto z + i(\ee -1, 0,
\dots, 0)$ and set
\[ \widetilde{N}_i (\ee)  = \; ]i \ll + (i-2)\ee, i \ll +i \ee[ \, 
\times \, ]0,1[ \, \times \, ]0, \ee[^{2n-2} \]
and
\[ \widetilde{\cn} (\ee) = \coprod_{i=0}^m \tau_i(B_i(\ll)) \bigcup
\coprod_{i=1}^m \widetilde{N}_i(\ee). \]
It is a simple matter to find a symplectomorphism $\ss$ of $\RR^2$ such
that $\ss \times id_{2n-2}$ embeds $\widetilde{\cn} (\ee)$ into an
arbitrarily small neighbourhood of $\cn (\ee)$.
It thus remains to show that $\widetilde{\cn} (\ee)$ may be asymptotically
filled by skinny ellipsoids. We try to fill $\widetilde{\cn} (\ee)$ by $\aa
E^{2n}(\pi, \dots, \pi,a)$ with $\aa$ small and $a$ large by packing the
$B_i(\ll)$ as in the proof of Proposition \ref{propositionlim=1} and
using $\widetilde{N}_i(\ee)$ to pass from $B_i(\ll)$ to $B_{i+1}(\ll)$. To
this end, think of $\aa E^{2n}(\pi, \dots, \pi, a)$ as fibered over
$\sq(\aa^2 a,1)$ with fibers $\frac{\bb^2}{\ee} \tr^{n-1}(\pi) \times
\sq^{n-1}(\ee)$, $\aa \ge \bb >0$.

We observe that the present packing problem is easier then the
one treated in Proposition \ref{propositionlim=1} inasmuch as now only a
part of $\aa E^{2n}(\pi, \dots, \pi, a)$ is embedded into a $B_i(\ll)$,
whence the ellipsoid fibres decrease slowlier.

\begin{figure}[h] 
 \begin{center}
  \psfrag{1}{$1$}
  \psfrag{y2}{$y_2$}
  \psfrag{y3}{$y_3$}
  \psfrag{m1}{$\dd_1$}
  \psfrag{m2}{$\dd_2$}
  \leavevmode\epsfbox{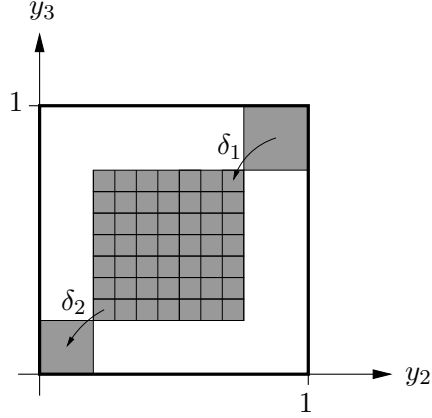}
 \end{center}
 \caption{The two deformations in $P_1(\tau_2(B_2(\ll)))$} \label{figure60}
\end{figure}

Let $\coprod_{i=1}^l P_i$ be a partition of $\tau_1(B_1(\ll))$ as in the
proof of Proposition \ref{propositionlim=1} and let $\gg \tr^{n-1}(\ll)
\, \times \, \sq^{n-1}(1)$ be the smallest fiber of $P_{l-1}$. Assume
that $l$ is so large that $\gg \ll < \ee$ and that $\aa$ is so small that
$\aa^2 \pi < \ee {\gg \ll}$. The image of the last ellipsoid fiber
mapped to $P_{l-1}$ is then contained in $\widetilde{N}_1(\ee)$, and we
may pass to $\tau_2(B_2(\ll))$. 
Having reached $P_1(\tau_2(B_2(\ll)))$, we first of all move the
ellipsoid fiber out of the connecting floor and then deform the fiber of
the second floor to a fiber with maximal $\tr^{n-1}$-factor ($\mu_1$ in
\mbox{Figure \ref{figure60}}). We then fill the remaining room in
$P_1(\tau_2(B_2(\ll)))$ as well as possible 
(cf.\ \mbox{Figure \ref{figure60}})
and proceed filling $\tau_2(B_2(\ll))$ as before. 
%
%
The above modification in the filling of $\tau_2(B_2(\ll))$ clearly does
not affect the result in Proposition \ref{propositionlim=1}. Going on in
the same way, we fill almost all of $\widetilde{\cn}(\ee)$.
This concludes the proof of \nobreak{Theorem \ref{theoremasymptotics}}. 
\proofend

\subsection{Refined asymptotic invariants}
\label{refinedasymptoticinvariants}

Theorem \ref{theoremasymptotics} shows that the asymptotic packing
numbers $p_\infty^E$ and $p_\infty^P$ are uninteresting invariants.
However, we may try to recapture some symplectic information on the
target space
by looking at the convergence speed. Given $(M, \oo)$ with $\Vol (M,
\oo) < \infty$ consider the function
\[ [1, \infty [ \ra \RR, \quad q \mapsto 1 - p_{1 \dots 1 q}^E(M, \oo)
\] and define a refined asymptotic invariant by 
\[
\aa_E (M,\oo) = \sup \{ \bb \,|\, 1 - p_{1\dots1q}^E(M,\oo) =
O(q^{-\bb}) \}.
\]
Define $\aa_P(M,\oo)$ in a similar way.

Let $U \in \co(n)$ with piecewise smooth boundary $\partial U$. 
Given a subset $S \subset \partial U$, let 
\[ S_s = \{ x \in U \, | \, d(x,S) < s \} \]
be the $s$-neighbourhood of $S$ in $U$. 
We say that $U$ is {\it admissible}, if there exists $\ee>0$
such that $U \setminus \partial U_\ee$ is connected.

\begin{example}   \label{ex:camel}
{\rm
Let $K(h,k) \subset \RR^{2n}$ be a camel space:
\[ K(h,k) = \{ x_1 <0 \} \cup \{x_1 >0 \} \cup H(h,k), \]
where
\[ H(h,k) = \left\{ \sum_{i=2}^n x_i^2 + \sum_{i=1}^n y_i^2 < h^2, \;
x_1=k \right\}.  \] 
Pick sequences $(h_i)_{i \in \NN}$ and $(k_i)_{i \in \NN}$ with
$h_1>h_2>\dots$, $h_i \ra 0$ and $0=k_1<k_2 < \dots$, $k_i \ra 1$,
let $C = \{ -1 < x_1, \dots, x_n, y_1, \dots, y_n <1 \}$ be a cube
and set
\[ U = C \cap \bigcap_{i=1}^\infty K(h_i,k_i).            \]
Then $C$ is not admissible. 
Thickening the walls and smoothing the
boundary, we obtain non admissible sets with smooth boundaries. \diam
}
\end{example}

\begin{proposition} \label{p:refined}
Let $U \in \co(n)$ be admissible and let $(M^{2n}, \oo)$ be a closed
symplectic manifold. Then
\begin{itemize}
\item[$(i)_E$]
\quad $\aa_E(U) \ge \frac{1}{n}$ \quad if $n \le 3$ or if $U \in \ce(n)$
\item[$(ii)_E$]
\quad $\aa_E(M,\oo) \ge \frac{1}{n}$ \quad if $n \le 3$
\item[$(i)_P$]
\quad $\aa_P(U) \ge \frac{1}{n}$
\item[$(ii)_P$]
\quad $\aa_P(M,\oo) \ge \frac{1}{n}$.
\end{itemize}
\end{proposition}

\noindent
{\bf Question.}
Given $\gg \in \, ]0,\frac{1}{2}[$, are there sets $U,V \in \co(2)$ with
$\aa_E(U) = \aa_P(V) = \gg \,$?
Candidates for such necessarily non admissible sets are the sets
described in Example \ref{ex:camel} with $(h_i), (k_i)$ chosen appropriately.
\diam
\\
\\
{\it Proof of Proposition \ref{p:refined}.}
ad $(i)_P$. If $U$ is a cube, the claim follows at once from
\nobreak{Proposition \ref{propositionhigherintocubes}}. 
If $U$ is an arbitrary admissible set, let 
\[ N_d = \{ (x_1, \dots, x_{2n}) \in \RR^{2n} \,|\, x_i \in d \ZZ,\,
1 \le i \le n \} \]
be the $d$-net in $\RR^{2n}$, and let $\cc_d$ be the union of all those
open cubes in $\RR^{2n} \setminus N_d$ which lie entirely in
$U$. Observe that $U \setminus \cc_d \subset {\partial U}_s$ 
$\partial U_s$
whenever 
$d \sqrt{2n}<s$. Let
$s_0 <\ee$ and $d_0 < s_0 / \sqrt{2n}$. Pick $\aa_0$ much smaller than
$d_0$ and exhaust $\cc_{d_0}$ with $\frac{\aa_0}{2} 
P^{2n}(\pi, \dots, \pi, a_0)$ by successively filling the cubes in 
$\cc_{d_0}$. More generally,
let $k \in \NN_0$, suppose that we almost exhausted $\cc_{d_0 / 2^k}$
by $\frac{\aa_0}{2^k} P^{2n}(\pi, \dots,\pi,a_k)$ and consider 
$\cc_{d_0 / 2^{k+1}}$. Then
\begin{equation} \label{e:r1}
U \setminus \cc_{d_0/2^{k+1}} \subset {\partial U}_{s_0/2^{k+1}}.
\end{equation}
We fill the cubes in $\cc_{d_0/2^k}$ by $\frac{\aa_0}{2^{k+1}} 
P^{2n}(\pi,\dots,\pi,a_{k+1})$ in the same order as we filled them by 
$\frac{\aa_0}{2^k} P^{2n}(\pi,\dots,\pi,a_k)$, but in between also fill the
cubes in $\cc_{d_0/2^{k+1}} \setminus \cc_{d_0/2^k}$. Observe that in
order to come back from a cube $C_{k+1} \in \cc_{d_0/2^{k+1}}$ to its
``mother-cube'' $C_k \in \cc_{d_0/2^k}$, we possibly have to use some
extra space in $C_k$, but that for the subsequent filling by 
$\frac{\aa_0}{2^{k+2}} P^{2n}(\pi, \dots,\pi,a_{k+2})$ this extra space will
be halved.

Since the $a_k$ were chosen maximal and since we exhaust more and more
of $U$,
\begin{equation}  \label{e:r2}
\lim_{k \ra \infty} \frac{a_{k+1}}{a_k} = 2^n.
\end{equation}
(\ref{e:r2}), the preceding remark and the case of a cube show that for
any $\dd >0$ there is a constant $C_1(\dd)$ such that for any $k$, any
$k' \le k$ and any $C_{k'} \in \cc_{d_0 / 2^{k'}}$
\begin{equation}  \label{e:r3}
\frac{\left|C_{k'} \setminus \mbox{image}\left(\frac{\aa_0}{2^{k'}} 
P^{2n}(\pi,\dots,\pi,a_k)\right)\right|}{|C_{k'}|} 
< C_1(\dd) a_k^{-\frac{1}{n}+\dd}.
\end{equation}
Let $\partial_kU$ be the $k$-dimensional components of $\partial U$, $0
\le k \le 2n-1$, and let $|\partial_kU|$ be their $k$-dimensional
volume. Then there are constants $c_k$ depending only on $U$ such that 
\[
\lim_{s \ra 0^+} \frac{|{\partial_kU}_s|}{s^{2n-k}} = c_k,
\]
whence
\begin{equation} \label{e:r4}
\lim_{s \ra 0^+} \frac{\big|{\partial U}_{s/2}\big|}
                      {\left|{\partial U}_s\right|}               =
\lim_{s \ra 0^+} \frac{\big|{\partial_{2n-1} U}_{s/2}\big|}
                      {\left|{\partial_{2n-1} U}_s\right|}
= \frac{1}{2}.
\end{equation}
(\ref{e:r1}), (\ref{e:r4}) and (\ref{e:r3}) imply that for any $\dd >0$
there is a constant $C_2(\dd)$ such that for any $k$ 
\begin{equation}  \label{e:r5}
\left|\cc_{d_0/2^k} \setminus
\mbox{image} \left( \frac{\aa_0}{2^k} P^{2n}(\pi,\dots,\pi,a_k)\right) \right|
< C_2(\dd) a_k^{-\frac{1}{n}+ \dd}.
\end{equation}
Next, (\ref{e:r1}), (\ref{e:r2}) and (\ref{e:r4}) show that for any $\dd
>0$ there is a constant $C_3(\dd)$ such that for any $k$
\begin{equation}  \label{e:r6}
|U \setminus \cc_{d_0/2^k}|
\le |B_{s_0/2^k}| < C_3(\dd) a_k^{-\frac{1}{n}+ \dd}.
\end{equation}
(i)$_P$ now follows from (\ref{e:r5}) and (\ref{e:r6}).

\smallskip
ad (ii)$_P$.
Cover $M$ with Darboux charts $(U_i, \ff_i)$, $i = 1, \dots ,m$, and
choose admissible subsets $V_i$ of $U_i$ such that the sets $W_i =
\ff_i(V_i)$ are disjoint and $\bigcup_{i=1}^m \overline{W}_i =M$. 
Choose
different points $p_i, q_i \in V_i$, set $\tilde{p}_i  = \ff_i(p_i)$, 
$\tilde{q}_i  = \ff_i(q_i)$, let $\widetilde{\ll}_i \colon [0,1] \ra M$
be a family of smooth, embedded and disjoint curves connecting
$\tilde{q}_i$ with $\tilde{p}_{i+1}$, and set $\ll_{i,j} = \ff_j^{-1}
(\widetilde{\ll}_i)$, $1 \le i \le m-1$, $1 \le j \le m$. We may assume that
near $q_i$ respectively $p_{i+1}$, $\ll_{i,i}$ respectively $\ll_{i,i+1}$ are
linear paths parallel to the $x_1$-axis.
As in the proof of {Theorem \ref{theoremasymptotics}} we find $\ee >0$
such that the $\widetilde{\ll}_i$ extend to disjoint symplectic embeddings
\[
\psi_i \colon [0,1] \times [-\ee, \ee]^{2n-1} \ra M
\]
whose compositions
$\psi_{i,i} = \ff_i^{-1} \circ \psi_i$ respectively $\psi_{i,i+1} =
\ff_{i+1} \circ \psi_i$ restrict to translations near $\{ 0 \} \times
[-\ee, \ee]^{2n-1}$ respectively $\{1 \} \times [-\ee,
\ee]^{2n-1}$. More generally, set $\psi_{i,i} = \ff_j^{-1} \circ
\psi_i$, and given $\dd \le \ee$, set 
\[
\psi_i^\dd = {\psi_i}|_{[0,1] \times [-\dd, \dd]^{2n-1}}
\quad \mbox{ and } \quad
\psi_{i,j}^\dd = {\psi_{i,j}}|_{[0,1] \times [-\dd, \dd]^{2n-1}}. 
\]
Let $\aa$ be so small that $\aa^2 \pi < 4 \dd^2$. We may then fill $M$
with $\aa P^{2n}(\pi, \dots , \pi, a)$ by successively filling $W_i \setminus
\coprod_{k=1}^{m-1} \mbox{image} \, \psi_k^\dd$ and passing from $W_i$
to $W_{i+1}$ with the help of $\psi_i^\dd$. 

In order to estimate the convergence speed of the filling of $W_i$, let
us look at the corresponding filling of $V_i$ instead. Set
\[
\ll_{i,j}^\dd = \{ x \in V_i \, | \, d(x, \, \mbox{image} \, \ll_{i,j}) 
< \dd \} 
\quad \mbox{ and } \quad 
V_i^\dd = V_i \setminus \coprod_j \ll_{i,j}^\dd.
\]
Let $L$ be a Lipschitz-constant for $\coprod_{i,j} \psi_{i,j}$. 
Then
\begin{equation}  \label{e:AAAA}
\mbox{image} \, \psi_{i,j}^\dd \; \subset \; \ll_{i,j}^{L \dd}.
\end{equation}
With $V_i$ also $V_i^0$ is admissible, and so there is $\dd_0 >0$ such
that $V_i^{L \dd_0}$ is connected. This and (\ref{e:AAAA}) show that we
may fill $V_i$ with a part of $\aa_0 P^{2n}(\pi, \dots , \pi, a_0)$ by
entering $V_i$ through $\ll_{i,i}^{L \dd_0}$, filling as much of $V_i^{L
\dd_0}$ as possible and leaving $V_i$ through $\ll_{i,i+1}^{L \dd_0}$.
Let $_i\cc_d^\dd$ be the union of those open cubes in $\RR^{2n}
\setminus N_d$ which lie entirely in $V_i^{L\dd}$. Then
\begin{equation}  \label{e:BBBB}
V_i \setminus {_i\cc_{d_0}^{\dd_0}} \;  \subset \; \coprod_{j=1}^{m-1}
\ll_{i,j}^{2 L \dd_0} \, \bigcup \; {(\partial V_i)}_{L \dd_0}
\end{equation}
whenever $d_0 \sqrt{2n} < L \dd_0$. Finally, 
\begin{equation}  \label{e:CCCC}
\lim_{s \ra 0^+} 
\frac{\left|\ll_{i,j}^{s/2}\right|}{\left|\ll_{i,j}^s \right|} 
= \frac{1}{2^{2n-1}} < \frac{1}{2}.
\end{equation}
(ii)$_P$ now follows from (\ref{e:BBBB}), (\ref{e:CCCC}) and the proof
of (i)$_P$.

\smallskip
ad (i)$_E$ and (ii)$_E$.
By the Folding Lemma, $E(\pi, a) \hookrightarrow P(\pi, (a+\pi)/2)$,
whence the case $n=2$ follows from (i)$_P$ and (ii)$_P$.

Let $n=3$, and let $U$ be a cube. We fill $U$ as described in 
\ref{highembeddingellipsoidsintocubes}.
This asymptotic packing problem resembles the one in the proof of
Proposition \ref{p:c}. Again, for given $a$, the region in $U$ not
covered by the image of the maximal ellipsoid $\aa E(\pi, \pi, a)$
fitting into $U$ decomposes into several disjoint regions $R_h(a)$, $2
\le h \le 4$. 
\begin{itemize}
\item[]
$R_2(a)$ is the space needed for folding.
\item[]
$R_3(a)$ is the union of the space needed to deform the ellipsoid
fibers and the space caused by the fact that 
the sum of the sizes of the ellipsoid fibres embedded into a column of
the cube fibre and the $x_3$-width of the space needed to deform one of
these ellipsoid fibres might be smaller than the size of the cube fibre.
\item[]
$R_4(a)$ is the space caused by the fact that the size of the ellipsoid
fibres decreases during the filling of a column of the cube fibre.
\end{itemize}
We compare $R_h(a)$ with $R_h(2^na) = R_h(8a)$. 
Let $\aa' E(\pi, \pi, 8a)$ be the maximal ellipsoid fitting into $U$. 
A volume comparison shows that for $a$ large $\aa '$ is very close to
$\aa /2$.
A similar but simpler analysis than in the proof of 
\nobreak{Proposition \ref{p:c}} 
now shows that given $\ee >0$ there is $a_0$ such that for any
$a \ge a_0$
\[ (2-\ee) \left| R_h(8a) \right| < \left| R_h(a) \right|, 
\qquad 2 \le h \le 4. \]
This implies the claim in case of a cube. The general case follows from
this case in the same way as (i)$_P$ and (ii)$_P$ followed from the case
of a cube.

Finally, let $E = E(b_1, \dots, b_n)$.
It follows from the description of Lagrangian folding in subsection
\ref{lagrangianfolding} and from Lemma \ref{lemmamatrix}(i) 
that given $n-1$ relatively prime numbers $k_1, \dots , k_{n-1}$ there
is an embedding $E^{2n}(\pi, \dots, \pi, a) \hookrightarrow \bb E(b_1,
\dots, b_n)$ whenever
\begin{equation} 
\left. \begin{array}{rcl} \label{e:r11}
    \frac{\pi}{\bb b_i} + \frac{\pi}{k_i \bb b_n} & < &   \frac{1}{k_i},
    \qquad 1 \le i \le n-1                           \\     
    \frac{\pi}{\bb b_n}        & < &   \frac{k_1 \cdots k_{n-1} \pi}{a}.
       \end{array}
   \right\}.
\end{equation}
W.l.o.g.\ we may set $b_n = 1$. (\ref{e:r11}) then reads 
\begin{equation}  \label{e:r12}
\left. \begin{array}{rcl}
    k_i \pi  & < &   (\bb -1) b_i, \qquad 1 \le i \le n-1              \\     
    a        & < &   k_1 \cdots k_{n-1} \bb.
       \end{array}
   \right\}.
\end{equation}
Pick some (large) constant $C$ and define $\bb$ by
\[
b_1 \cdots b_{n-1} \bb^n = \pi^{n-1} \left(a+Ca^\frac{n-1}{n}\right).
\]
Moreover, pick $n-1$ prime numbers $p_1, \dots, p_{n-1}$, let $l$ be the
least common multiple of $\{ p_i-p_j \,|\, 1 \le i < j \le n-1 \}$, define $m_i$, $1 \le i
\le n-1$, by 
\[
m_i = \max \{ m \in \NN \,|\, m_i l-p_i < (\bb-1) b_i/\pi \}
\]
and set $k_i = m_il-p_i$. We claim that the $k_i$ are relatively
prime. Indeed, assume that for some $i \neq j$
\begin{equation}  \label{e:r13}
d \,|\, m_il-p_i \quad \mbox{ and } \quad d \,|\, m_jl-p_j.
\end{equation}
Then $d$ divides $(m_i l - p_i) - (m_j l - p_j) = p_i - p_j$, and hence,
by the definition of $l$, $d$ divides $l$. But then, by (\ref{e:r13}), 
$d$ divides $p_i$ and $p_j$, whence $d=1$.

The first $n-1$ inequalities in (\ref{e:r12}) hold true by the
definition of the $k_i$,
and since $b_i \le 1$,
\begin{eqnarray*}
\pi^{n-1} k_1 \cdots k_{n-1} \bb &>& (\bb b_1 -l-1) \cdots (\bb_{n-1}
-l-1) \bb \\
& = & b_1 \cdots b_{n-1} \bb^n + \sum_{i=1}^{n-1} (-1)^ic_i \bb^{n-i}, 
\end{eqnarray*}
where the $c_i$ are positive constants depending
only on $b_1, \dots, b_{n-1}$ and $l$. For $a$ large enough the last
expression is larger than $b_1 \cdots b_{n-1} \bb^n - c_1 \bb^{n-1}$,
which equals
\[
\pi^{n-1} \left(a+Ca^\frac{n-1}{n}\right) 
- c_1 \left(\frac{\pi^{n-1}}{b_1 \cdots b_{n-1}}\right)^\frac{n-1}{n}
\left( a+C a^\frac{n-1}{n} \right)^\frac{n-1}{n}
\]
and this is larger than $\pi^{n-1}a$ whenever $a$ and $C$ are large
enough.

Finally, we have that
\[
\frac{|E^{2n}(\pi, \dots, \pi, a)|}{|\bb E(b_1, \dots, b_n)|} 
= \frac{\pi^{n-1}a}{\bb^n b_1 \cdots b_{n-1}} 
= \frac{1}{1+ C a^{-\frac{1}{n}}} 
= 1- Ca^{-\frac{1}{n}} + o\left(a^{-\frac{1}{n}}\right),
\]
from which the second claim in (i)$_E$ follows.
\proofend

\noindent
{\bf Remark.}
Suppose that we knew that there is a natural number $k$ such that 
the cube $C^{2n}$ admits a full symplectic packing by $k$ equal balls
and such that the space of symplectic embeddings of $k$ equal balls into
$C^{2n}$ is unknotted.
Combining such a result with {Proposition \ref{p:c}} and the
techniques used in the proof of {Theorem \ref{theoremasymptotics}} and
{Proposition \ref{p:refined}} we may derive that 
\[
\aa_E (U) \ge \frac{1}{2n} 
\quad \mbox{ and } \quad \aa_E(M, \oo) \ge \frac{1}{2n}
\]
for any admissible $U \in \co(n)$ and any closed symplectic manifold
$(M^{2n}, \oo)$.

\subsection{Higher order symplectic invariants}  \label{higherorder}

The construction of good higher order invariants for subsets of
$\RR^{2n}$ has turned out to be a
difficult problem in symplectic topology. The known such
invariants are Ekeland-Hofer capacities \cite{EH1, EH2} 
and symplectic homology \cite{FH, FHW}, which 
both rely on the variational study of periodic orbits of certain
Hamiltonian systems, and the symplectic homology constructed via
generating functions \cite{T2}. 
We propose here some higher order invariants which are based on an
embedding approach. 

Let $(M^{2n}, \oo)$ be a symplectic manifold and let
\[
e_1(M,\oo) = \sup \{ A\,|\, B^{2n}(A) \mbox{ symplectically embeds in }
(M,\oo) \}
\]
be the Gromov-width of $(M,\oo)$. We inductively define $n-1$ other
invariants by 
\begin{eqnarray*}
e_i(M,\oo) = 
\sup \{ A \,| & \!\!\!\! E^{2n}(e_1(M,\oo), \dots ,e_{i-1}(M,\oo),  A, \dots, A) \\
              & \quad \mbox{ symplectically embeds in } (M,\oo) \}.
\end{eqnarray*}
Similarly, given $U \in \co(n)$, let 
\[
e^n(U) = \inf \{ A\,|\, U \mbox{ symplectically embeds in }
B^{2n}(A) \}
\]
and inductively define $n-1$ other invariants $e^i(U)$ by
\[
e^i(U) = \inf \{ A \,|\, U \mbox{ symplectically embeds in } E^{2n}(A,
\dots , A, e^{i+1}(U), \dots ,e^n(U)  \}.
\]
Clearly, 
\[ e_1(M, \oo) \le e_2(M, \oo) \le \dots \le e_n(M, \oo) \]
and 
\[ e^1(M, \oo) \le e^2(M, \oo) \le \dots \le e^n(M, \oo). \]
Moreover,
$
e_i(M,\aa \oo) = |\aa| \, e_i(M,\oo)$ and $e^i(U,\aa \oo_0) = |\aa| \,
e^i(U,\oo_0)$ for all $\aa \in \RR \setminus \{0\}$, 
and $e_i$ and $e^i$ are indeed invariants, that is 
$e_i(M,\oo) = e_i(N,\tau)$ and $e^i(U,\oo_0) = e^i(V,\oo_0)$
if there are symplectomorphisms $\ff \colon (M,\oo) \ra (N,\tau)$ and
$\psi \colon (U,\oo_0) \ra (V,\oo_0)$.

\begin{example}   \label{ex:order}
{\rm
Ekeland-Hofer capacities show that
\[
e_i(E(a_1, \dots, a_n)) = a_i, \qquad 1 \le i \le n,
\]
and 
\[
e^i(E(a_1, \dots, a_n)) = a_i, \qquad 1 \le i \le n, \qquad \mbox{if } 2 a_1
\ge a_n. 
\]
\hfill{\diam}
}
\end{example}
$e_1$ and $e^n$ are also monotone and nontrivial, and are hence
symplectic capacities (see \cite{HZ} for the axioms of a symplectic
capacity). This, however, does not hold true for any of the
higher invariants.
Indeed, let $Z(\pi) = D(\pi) \times \RR^{2n-2}$ be the standard
symplectic cylinder. Then
\[
e_i(Z(\pi)) = \infty \quad \mbox{ for all } i \ge 2.
\]
Moreover, Example \ref{ex:order} and Theorem 2A show that none of the
$e_i$, $i \ge 2$, is monotone, and the same holds true for $e^i$, $i \le
n-1$. For instance, set
$U_\ll = \frac{4}{3} E(\ll^{-1} \pi, \ll \pi)$ 
and 
$V = E(\pi, 2 \pi)$.
By Theorem \ref{theoremasymptotics}, $U_\ll$ symplectically embeds in
$V$ and $e^2(U_\ll)$ is near to $\frac{4}{3} \pi$ if $\ll$ is large. 
Then also $e^1(U_\ll)$ is near to $\frac{4}{3} \pi$; but $e^1(V) = \pi$.

Similar invariants may be constructed by looking at polydiscs instead of
ellipsoids. 

These considerations indicate that it should be difficult to construct
higher order symplectic capacities via an embedding approach.

\section{Appendix}

\subsection*{A. Computer programs}

All the Mathematica programs of this appendix may be found under
\begin{center}
ftp://ftp.math.ethz.ch/pub/papers/schlenk/folding.m
\end{center}
For convenience, in the
programs (but not in the text) both the $u$-axis and the capacity-axis
are rescaled by a factor $1/\pi$.

\subsubsection*{A1. The estimate $s_{EB}$}

As said at the beginning of \ref{embeddingellipsoidsintoballs} we fix
$a$ and $u_1$ and try to embed $E(\pi, a)$ into
$B^4(2\pi +(1-2\pi/a)u_1)$ by multiple folding. If this works,
we set $A(a, u_1) = 2\pi + (1-2\pi/a) u_1$ and $A(a, u_1) = a$ otherwise.
\begin{verbatim}
A[a_, u1_] :=
  Block[{A=2+(1-2/a)u1},
     j  = 2;         
     uj = (a+1)/(a-1)u1-a/(a-1);
     rj = a-u1-uj;
     lj = rj/a;
     While[True,
           Which[EvenQ[j],                      
                   If[rj <= uj,                  
                      Return[A],
                      If[uj <= 2lj,                 
                         Return[a],
                         j++;                      
                         uj = a/(a-2)(uj-2lj);   
                         rj = rj-uj;             
                         li = lj;              
                         lj = rj/a             
                        ]
                     ],
                 OddQ[j],                       
                   If[rj <= uj+li,                 
                      Return[A],             
                      j++;                       
                      uj = (a+1)/(a-1)uj;          
                      rj = rj-uj;           
                      lj = rj/a                
                     ]
                ]                     
            ]                    
         ]              
\end{verbatim}
This program just does what we proposed to do in 
\ref{embeddingellipsoidsintoballs} in order to decide if the embedding
attempt associated with $u_1$ succeeds or fails.
Note, however, that in the {\tt Oddq[j]}-part, we did not check whether
the upper
left corner of $F_{j+2}$ is contained in $T(A,A)$. However, this negligence
does not cause troubles, since if the left edge of $F_{j+2}$ indeed
exceeds $T(A,A)$, the embedding attempt will fail in the subsequent
{\tt EvenQ[j+1]}-part. In fact, that the left edge of $F_{j+2}$ exceeds
$T(A,A)$
means that $l_{j+1} > u_{j+1}$; hence $r_{j+1} > u_{j+1}$ (since
otherwise the embedding attempt would have succeeded in the preceding
{\tt OddQ[j]}-part), but $u_{j+1} \le 2 l_{j+1}$.

\medskip
Writing again $u_0$ for the minimal $u_1$ which leads to an embedding, $A(a,
u_1)$ is equal to $a$ for $u_1<u_0$ and it is a linear increasing function for
$u_1 \ge u_0$. Since, by (\ref{equation23}), we may assume that $u_0 \le a/2$, we have
$A(a, u_0) \le \pi + a/2 <a$, whence $u_0$ is found up to
accuracy $acc/2$  
by the following bisectional algorithm.
\begin{verbatim}
u0[a_, acc_] := 
  Block[{},     
     b = a/(a+1);
     c = a/2;
     u1 = (b+c)/2;
     While[(c-b)/2 > acc/2,
       If[A[a,u1] < a, c=u1, b=u1];
       u1 = (b+c)/2
          ];
     Return[u1]
       ]
\end{verbatim}
Here the choice $b=a \pi /(a+\pi)$ is also based on (\ref{equation23}).
Up to accuracy $acc$, the resulting estimate $s_{EB}(a)$ is given by 
\begin{verbatim}
sEB[a_, acc_] := 2 + (1-2/a)u0[a,acc].
\end{verbatim}

\subsubsection*{A2. The estimate $s_{EC}$}

Given $a$ and $u_1$, we first calculate the height of the image of the
corresponding embedding. The following program is easily understood by
looking at \mbox{Figure \ref{figure18a}}.
\begin{verbatim}
h[a_, u1_] :=
  Block[{l1=1-u1/a},          
     j  = 2;             
     uj = (a+1)/(a-1)u1-a/(a-1);         
     rj = a-u1-uj;                       
     lj = rj/a;                          
     hj = 2l1;                            
     While[rj > u1+l1 - lj,         
           j++;                     
           uj = (a+1)/(a-1)uj;        
           rj = rj-uj;               
           li = lj;                 
           lj = rj/a;                
           If[EvenQ[j], hj = hj+2li]
          ];
     Which[EvenQ[j], 
           hj = hj+lj,
           OddQ[j],  
           hj = hj+Max[li,2lj]
          ];
     Return[hj]                
       ]
\end{verbatim}
As explained in \ref{embeddingellipsoidsintocubes}, the optimal folding
point $u_1$ is the $u$-coordinate of the unique intersection point of
$h(a,u_1)$ and $w(a,u_1)$. It may thus be found again by a bisectional
algorithm.

\begin{verbatim}
u0[a_, acc_] := 
  Block[{},
     b = a/(a+1);
     c = a/2;
     u1 = (b+c)/2;
     While[(c-b)/2 > acc/2,
       If[h[a,u1] > 1+(1-1/a)u1, b=u1, c=u1];
       u1 = (b+c)/2
          ];
     Return[u1]
      ]
\end{verbatim}
Again, the choices $b = a \pi / (a+\pi)$ and $c = a/2$ reflect that we
fold at least twice in which case $u_1 \ge l_1$ must hold true. 
Up to accuracy $acc$, the resulting estimate $s_{EC}(a)$ is given by
\begin{verbatim}
sEC[a_, acc_] := 1+(1-1/a)u0[a,acc].
\end{verbatim}

\subsection*{B. Report on the Gromov width of closed symplectic manifolds}

Recall that given any symplectic manifold $(M^{2n}, \oo)$ its
Gromov width is defined by
\[
w(M, \oo) = \sup \{ c \; | \mbox{ there is a symplectic embedding } 
(B^{2n}(c), \oo_0) \hookrightarrow (M,\oo) \}.
\]
Historically, the width provided the first example of a symplectic
capacity. Giving the size of the largest Darboux chart of $(M,\oo)$, the
width is always positive, and in the closed case it is finite. We now
restrict to closed manifolds and define an equivalent packing invariant
by
\[
p(M^{2n}, \oo) = \frac{|B^{2n}(w(M,\oo))|}{\mbox{Vol}(M,\oo)}
= \frac{w(M, \oo)^n}{n! \, \mbox{Vol} (M,\oo)}.
\]
In two dimensions the width is the volume and $p=1$ 
(see Theorem \ref{t:volume}). The basic
result to discover rigidity in higher dimensions is a version of
Gromov's Non-Squeezing Theorem \cite{LMP}. 
\\
\\
{\bf Non-Squeezing Theorem (compact version)}     \label{theoremnonsqueezing}
{\it Let $(M^{2n}, \oo)$ be closed, let $\ss$ be an area form on $S^2$
such that $\int_{S^2} \ss = 1$ and
assume that there is a symplectic embedding
$
B^{2n+2}(c) \hookrightarrow (M \times S^2, \oo \oplus a \ss).
$
Then $a \ge c$.
}
\\
\\
{\bf Remark.}
More generally, let $S^2 \hookrightarrow M \ltimes S^2 \xrightarrow{\pi} M$
be an oriented $S^2$-bundle over a closed manifold $M$ and let $\oo$ be
a symplectic form on $M \ltimes S^2$ whose restriction to the fibers is
nondegenerate and induces the given orientation. 
In particular, $a = \langle [\oo], [pt \times S^2] \rangle >0$. Then the
proof of the above Non-Squeezing Theorem also implies that $c \le a$
whenever $B^{2n+2} (c)$ symplectically embeds in $(M \ltimes S^2,
\oo)$. We will verify this below in the case where $M$ is 2-dimensional.
\diam

Since the theory of $J$-holomorphic curves works best in dimension four,
the deepest results on the Gromov-width have been
proved for 4-manifolds. 
Given a symplectic 4-manifold $(M,\oo)$, let $c_1$ be the first
Chern class of $(M,\oo)$ with respect to the contractible set of almost
complex structures compatible with $\oo$.
Let $\cc$ be the class of symplectic 4-manifolds $(M,\oo)$ for which
there exists a class $A \in H_2(M; \ZZ)$ with non-zero Gromov invariant
and $c_1(A) + A^2 \neq 0$.  
Recall that a symplectic 4-manifold is called {\it
rational} if it is the symplectic blow-up of $\CC \PP^2$ 
and that it is said to be {\it ruled} if it is an $S^2$-bundle over
a Riemann surface. 
The class $\cc$ consists of symplectic blow-ups of 
\begin{itemize}
\item[$\bullet$] rational and ruled manifolds;
\item[$\bullet$] manifolds with $b_1 = 0$ and $b_2^+ =1$;
\item[$\bullet$] manifolds with $b_1 =2$ and $(H^1(M;\ZZ))^2 \neq 0$.
\end{itemize}
We refer to \cite{M2} for more information on the class $\cc$.

Recall that by definition an {\it exceptional sphere} in a symplectic
4-manifold $(M,\oo)$ is a symplectically embedded 2-sphere $S$ of
self-intersection number $S \cdot S =-1$,
and that $(M,\oo)$ is said to be {\it minimal} 
if it contains no exceptional spheres.
Combining the technique of symplectic blowing-up with Taubes theory 
of Gromov invariants, Biran \cite[Theorem 6.A]{B1}
showed that for the symplectic 4-manifolds $(M,\oo)$ in class $\cc$ all packing
obstructions come from exceptional spheres in the symplectic blow-up of
$(M, \oo)$ and from the volume constraint.
His result suffices to compute the Gromov-width of all minimal manifolds
in the class $\cc$.
\begin{theorem}  \label{theorembiran2f}
{\rm (Biran \cite[Theorem 2.F]{B1})}
Let $(M, \oo)$ be a closed symplectic 4-manifold in the class $\cc$
which is minimal and neither rational nor ruled. Then $p(M, \oo) = 1$.
\end{theorem}
Examples of manifolds satisfying the conditions of the above theorem are
hyper-elliptic surfaces and the surfaces of Barlow, Dolgachev and
Enriques, all viewed as K\"ahler surfaces.
\\

We next look at minimal manifolds which are rational or ruled. 

Let $\oo_{SF}$ be the unique $\mbox{U}(3)$-invariant K\"ahler form on
$\CC \PP^2$ whose integral over $\CC \PP^1$ equals $\pi$. In the
rational case, by a theorem of Taubes \cite{Ta}, $(M, \oo)$ is
symplectomorphic to $(\CC \PP^2, a \oo_{SF})$ for some $a >0$, thus $p(M,
\oo) =1$.

Denote by $\Sigma_g$ the Riemann surface of genus $g$. There are
exactly two orientable $S^2$-bundles with base $\Sigma_g$, namely the
trivial bundle $\pi \colon \Sigma_g \times S^2 \ra \Sigma_g$ and the
nontrivial bundle $\pi \colon \Sigma_g \ltimes S^2 \ra \Sigma_g$
\cite[Lemma 6.25]{MS}. Such a manifold is called a ruled surface. 
$\Sigma_g \ltimes S^2$ is the projectivization $\PP(L_1 \oplus \CC)$
of the complex rank two bundle $L_1 \oplus \CC$ over $\Sigma_g$, where
$L_1$ is a holomorphic line bundle of Chern index $1$.
A symplectic form $\oo$ on a ruled surface is called {\it compatible}
with the given ruling $\pi$ if it restricts on each fiber to a
symplectic form. Such a symplectic manifold is then called a {\it ruled
symplectic manifold}. It is known that every symplectic structure on a
ruled surface is diffeomorphic to a form compatible with the given
ruling $\pi$ via a diffeomorphism which acts trivially on homology, and
that two cohomologous symplectic forms compatible with the same ruling
are isotopic \cite{LM4}.
A symplectic form $\oo$ on a ruled surface is thus determined up to
diffeomorphism by the class $[\oo] \in H^2(M; \RR)$.

Fix now an orientation of the fibers of the given ruled symplectic
manifold. We say that $\oo$ is {\it admissible} if its restriction to
each fiber induces the given orientation.

Consider first the trivial bundle $\Sigma_g \times S^2$ with its given
orientation, and let $\{ B = [\Sigma_g \times pt], F = [pt \times S^2]
\}$ be a basis of
$H^2(M;\ZZ)$ (here and henceforth we identify homology and cohomology
via Poincar\'e duality). Then a cohomology class $c = bB + aF$ can
be represented by an admissible form if and only if $c(B) = a >0$ and
$c(F) = b >0$. We write $\Sigma_g(a) \times S^2(b)$ for this ruled
symplectic manifold. 

In case of the nontrivial bundle $\Sigma_g \ltimes S^2$ a basis of
$H^2(\Sigma_g \ltimes S^2; \ZZ)$ is given by $\{ A, F \}$, where $A$ is
the class of a section with selfintersection number $-1$ and $F$ is the
fiber class. Set $B= A+ \dfrac{F}{2}$. $\{ B, F\}$ is then a basis of
$H^2(\Sigma_g \ltimes S^2; \RR)$ with $B \cdot B = F \cdot F = 0$ and $B
\cdot F =1$. It turns out that in case $g=0$ a form $c = bB + a F$ can
be represented by an admissible form if and only if $a > \frac{b}{2}
>0$, while in case $g \ge 1$ this is possible if and
only if $a>0$ and $b>0$ \cite[Theorem 6.27]{MS}.
We write $(\Sigma_g \ltimes S^2, \oo_{ab})$ for this ruled symplectic
manifold.

Finally note that each admissible form is cohomologous to a standard
K\"ahler form. For the trivial bundles these are just the split forms,
and for the non-trivial bundles we refer to \cite[p.\ 276]{L}.
\begin{theorem}  \label{t:gromov}
Let $(M^4, \oo)$ be a ruled symplectic manifold, i.e.\ either 
$(M, \oo) = \Sigma_g(a) \times S^2(b)$ or 
$(M, \oo) = (\Sigma_g \ltimes S^2, \oo_{ab})$.
If $(M, \oo) = S^2(a) \times S^2(b)$ we may assume that $a \ge b$. Then 
\begin{itemize}
\item[(i)] 
  $ p(S^2(a) \times S^2(b)) \, = \, p(S^2 \ltimes S^2, \oo_{ab}) \, = \,
    \frac{b}{2a} $
\item[(ii)]
  $ p(\Sigma_g(a) \times S^2(b)) \, = \, p(\Sigma_g \ltimes S^2,
    \oo_{ab}) \, = \, \min \{ 1, \frac{b}{2a} \} \mbox{ if } g \ge 1 $
\end{itemize}
\end{theorem}
The statements for the trivial bundles are proved in \cite[Theorem
6.1.A]{B1}, and the ones for the non-trivial bundles
are calculated in \cite{S}. Observe that the upper bounds predicted by
the Non-Squeezing Theorem and the volume condition are sharp in all
cases.
Explicit maximal embeddings are easily found for $g=0$ and for $g \ge 1$
if $a \ge b$ \cite{S}, but no explicit maximal embeddings are
known for $g \ge 1$ if $a<b$.

Also notice that 
$p(S^2(b) \times \Sigma_g(a)) = \min \{ 1, \frac{b}{2a} \}$ if $g \ge 1$ 
implies that the Non-Squeezing Theorem does not remain
valid if the sphere is replaced by any other closed surface.
\\

If $(M^4, \oo)$ does not belong to the class $\cc$ only very few is
known about $p(M, \oo)$. Indeed, no obstructions to full packings are known.
Some flexibility results for products of higher genus surfaces were
found by Jiang.

\begin{theorem} 
{\rm (Jiang \cite[Corollary 3.3 and 3.4]{J})}   \label{theoremjiang}
Let $\SS$ be any closed surface of area $a>1$.
\begin{itemize}
\item[(i)]
Let $T^2$ be the 2-torus. There is a constant $C>0$ such that $p(T^2(1)
\times \SS(a)) \ge C$.
\item[(ii)]
Let $g \ge 2$. There is a constant $C(g) >0$ depending only on $g$ such
that $w(\SS_g (1) \times \SS(a)) \ge C(g) \log a$.
\end{itemize}
\end{theorem}

\noindent
{\bf Remark.}
If $\SS = S^2$ Birans sharp result in Theorem
\ref{t:gromov} is of course much better.   \diam
\begin{example}    \label{examplejiang}
{\rm
Set $R(a) = \{ (x,y) \in \RR^2 \, | \, 0<x<1, \; 0<y<a \}$, and consider the
linear symplectic map
\begin{eqnarray*}
\ff \colon  (R(a) \times R(a), dx_1 \wedge dy_1 + dx_2 \wedge dy_2)  &  \ra   &
                (\RR^2 \times \RR^2, dx_1 \wedge dy_1 + dx_2 \wedge dy_2)   \\
      (x_1,y_1,x_2,y_2) & \mapsto & (x_1+y_2,y_1, -y_2, y_1+x_2).
\end{eqnarray*}
Let $p \colon \RR^2 \ra T^2 = \RR / \ZZ \times \RR
/ \ZZ$ be the projection onto the standard symplectic torus. Then $p
\circ \ff \colon R(a) \times R(a) \ra T^2 \times \RR^2$ is an embedding;
indeed, given $(x_1,y_1,x_2,y_2)$ and $(x_1',y_1',x_2',y_2')$ with 
\begin{eqnarray}
x_1 + y_2 & \equiv & x_1' + y_2' \mod \ZZ  \label{equationI}  \\
      y_1 & \equiv &        y_1' \mod \ZZ  \label{equationII} \\
     -y_2 & =      & -y_2'                 \label{equationIII}\\
y_1 + x_2 & =      & y_1' + x_2'           \label{equationIV} 
\end{eqnarray}
(\ref{equationIII}) gives $y_2 = y_2'$ and thus (\ref{equationI})
implies $x_1 \equiv x_1' \mod \ZZ$ whence $x_1 = x_1'$. Moreover,
(\ref{equationII}) and (\ref{equationIV}) show that $y_1-y_1' = x_2'-x_2
\equiv 0 \mod \ZZ$, hence $x_2 = x_2'$ and $y_1 = y_1'$.

Next observe that $p \circ \ff (R(a) \times R(a)) \subset T^2 \times
]-a,0[ \times ]-a-1, a+1[$. Thus $R(a) \times R(a)$ embeds in $T^2(1)
\times \SS(2a(a+1))$, and since $B^4(a)$ embeds in $R(a) \times R(a)$
and $B^4(1)$ embeds in $T^2(1) \times \SS(a)$ for any $a \ge 1$, we have
shown

\begin{proposition}  \label{propositionjiang}
Let $a\ge 1$. Then
\[
p(T^2(1) \times \SS(a)) \ge \frac{\max \{ a+1 - \sqrt{2a+1},2 \}}{4a}.
\]
\end{proposition}
In particular, the constant $C$ in Theorem \ref{theoremjiang}(i) can be
chosen to be $C=1/8$.
\diam
}
\end{example} 

It would be interesting to have a complete list of those symplectic
4-manifolds with $p(M,\oo) =1$.
As we have seen above, the minimal such manifolds in class $\cc$ are
those which are not ruled, the trivial bundles $\Sigma (a) \times S^2(b)$
with $g(\Sigma) \ge 1$ and $a \ge 2b$ and the nontrivial bundles
$(\Sigma \ltimes S^2, \oo_{ab})$ with $g(\Sigma) \ge 1$ and $a \le 0$.
Combining the techniques of \cite{B1} with Donaldson's existence result
for symplectic submanifolds, Biran \cite{B2} found examples with
$p(M,\oo) =1$ which do not belong to $\cc$.

\medskip
In higher dimensions almost no flexibility results are known. Note
however that for the standard K\"ahler form $\oo_{SF}$ on $\CC P^n$
we have $p(\CC P^n, \oo_{SF}) =1$ (see
e.g.\ \cite{MP}), and that the technique used in Example 
\ref{examplejiang} shows
that given any constant form $\oo$ on $T^{2n}$ and an area form $\ss$ on $\SS$
with $\int_\SS \ss =1$ there is a constant $C>0$
such that $p(T^{2n} \times \SS, \oo \oplus a\ss) \ge C$ 
(\cite[Theorem 3.1]{J}).

\medskip

Felix Schlenk, Mathematik, ETH Zentrum, 8092 Z\"urich, Switzerland

\medskip
{\it E-mail address:} felix@math.ethz.ch

\end{document}